\documentclass[12pt,psamsfonts]{article} 
\usepackage{a4,fullpage,amssymb,epsf,psfrag,times}
\usepackage{graphicx}
\usepackage{amsmath}
\usepackage{amsfonts}
\usepackage{enumerate}
\usepackage{latexsym}
\usepackage{epsfig}
\usepackage{amsthm}  
\usepackage{amsmath}
\usepackage{amssymb} 
\usepackage{latexsym}
\usepackage{epsfig} 
\usepackage{amssymb,latexsym}
\usepackage[all]{xy}

\makeatletter\@addtoreset{equation}{section}\makeatother
\makeatletter\@addtoreset{figure}{section}\makeatother
\makeatletter\@addtoreset{table}{section}\makeatother

\newtheorem{theorem}{Theorem}[section]

\newtheorem{prop}[theorem]{Proposition}
\newtheorem{lemma}[theorem]{Lemma}

\newtheorem{cor}[theorem]{Corollary}

\newtheorem{roughclassification}[theorem]{Rough Classification}
\theoremstyle{plain}

\newcommand{\R}{{\mathbb R}}
\newcommand{\C}{{\mathbb C}}
\newcommand{\Z}{{\mathbb Z}}

\newcommand{\T}{{\mathbb T}}

\newcommand{\op}[1]{\!\!\mathop{\rm ~#1}\nolimits}

\newcommand{\scriptop}[1]{\!\!\mathop{\mbox{\rm \scriptsize ~#1}}\nolimits}

\newenvironment{remark}{\refstepcounter{theorem}\par\medskip\noindent{\bf Remark~\thetheorem~~}}{\unskip\nobreak\hfill\hbox{ $\oslash$}\par\bigskip}

\newenvironment{example}{\refstepcounter{theorem}\par\medskip\noindent{\bf Example~\thetheorem~~}}{\unskip\nobreak\hfill\hbox{ $\oslash$}\par\bigskip}

\newenvironment{definition}{\refstepcounter{theorem}\par\medskip\noindent{\bf Definition~\thetheorem~~}}{\unskip\nobreak\hfill\hbox{ $\oslash$}\par\bigskip}

\newcommand{\got}[1]{\mathfrak{#1}}

\begin{document}

\title{Symplectic actions of two tori on four manifolds}
\author{Alvaro Pelayo\thanks{Partly
funded by a \emph{Rackham Predoctoral Fellowship} (2006-2007) and a \emph{Rackham Dissertation Fellowship} (2005-2006).}}

\maketitle

\begin{abstract}
We classify
symplectic actions of $2$\--tori on compact, connected symplectic 
$4$\--manifolds, up to equivariant symplectomorphisms. 
This extends results of Atiyah, Guillemin--Sternberg, Delzant and Benoist. 
The classification is in terms of a collection of invariants,
 which are invariants of the topology of the manifold,  
of the torus action and of the symplectic form. 
We construct explicit models of such symplectic manifolds with torus actions, defined in terms of these invariants. 

We also classify, up to equivariant symplectomorphisms,
symplectic actions of $(2n-2)$\--dimensional tori
on $2n$\--dimensional symplectic manifolds,
when at least one orbit is a $(2n-2)$\--dimensional 
symplectic submanifold. Then
we show that a $2n$\--dimensional symplectic manifold $(M, \, \sigma)$ equipped
with a free symplectic action of a $(2n-2)$\--dimensional torus with at least one $(2n-2)$\--dimensional
symplectic orbit is equivariantly diffeomorphic to $M/T \times T$ equipped with
the translational action of $T$. Thus
two such symplectic manifolds are equivariantly diffeomorphic if and only if their
orbit spaces are surfaces of the same genus. 
\end{abstract}

%\tableofcontents

\section{Introduction}

We extend
the theory of Atiyah \cite{A}, Guillemin \cite{G}, 
Guillemin--Sternberg \cite{GS}, Delzant \cite{Delzant}, and 
Benoist \cite{Be} to symplectic actions of tori which are not necessarily Hamiltonian. 
Although Hamiltonian actions of $n$\--dimensional tori
on $2n$\--dimensional manifolds are present in many integrable systems in classical mechanics,
non--Hamiltonian actions occur also in physics, c.f. Novikov's article \cite{novikov}. 
Interest on non--Hamiltonian motions may be found in the physics literature, for example: Sergi--Ferrario \cite{SF}, Tarasov \cite{Tarasov} and Tuckerman's articles \cite{Tucker1},
\cite{Tucker2} and the references therein.

In this paper we give a classification
of symplectic actions of $2$\--tori on compact connected symplectic $4$\--manifolds in terms of 
a collection of invariants, some of which are algebraic while others are topological or geometric.
A consequence of our classifications is
the following.

\begin{theorem}
The only $4$\--dimensional symplectic manifold equipped with
a non--locally--free and non--Hamiltonian effective symplectic action of a $2$\--torus is,
up to equivariant symplectomorphisms, the product $\T^2 \times S^2$,
where
$\T^2=(\R/\Z)^2$ and
the first factor of $\mathbb{T}^2$ acts on the left factor by translations
on one component, and the second factor acts on $S^2$ by rotations about the
vertical axis of $S^2$. The symplectic form is
a positive linear combination of the
standard translation invariant form on $\mathbb{T}^2$
and the standard rotation invariant form on $S^2$.
\end{theorem}

Duistermaat
and I showed in \cite{DuPe} that a compact connected symplectic manifold $(M, \, \sigma)$
with a symplectic torus action with at least one coisotropic principal orbit  is
an associated $G$\--bundle $G\times _HM_{\textup{h}}$ 
whose fiber is a symplectic toric manifold $M_{\textup{h}}$ 
with $T_{\textup{h}}$\--action and whose base $G/H$ is a torus bundle over a torus.
Here $T_{\textup{h}}$ is the unique maximal subtorus 
of $T$ which acts in a Hamiltonian fashion on $M$,
$G$ is a two\--step nilpotent Lie group
which is an extension of the torus $T$, and $H$ 
is a commutative closed Lie subgroup of $G$ 
which acts on $M_{\textup{h}}$ via $T_{\textup{h}}$ and is
defined in terms of the holonomy of a certain
connection for the principal bundle $M_{\textup{reg}} \to M_{\textup{reg}}/T$,
where $M_{\textup{reg}}$ is the set of points where the action is free.
Precisely, $G=T \times (\mathfrak{l}/\mathfrak{t}_{\textup{h}})^*$ where
$\mathfrak{l}$ is the kernel of the antisymmetric bilinear form $\sigma^{\mathfrak{t}}$ on $\mathfrak{t}$ which gives the restriction of $\sigma$ to the orbits, and $\mathfrak{t}_{\textup{h}}$ is
the Lie algebra of the torus $T_{\textup{h}}$. The additive
group $(\mathfrak{l}/\mathfrak{t}_{\textup{h}})^* \subset \mathfrak{l}^*$, 
viewed as the set of linear forms on $\mathfrak{l}$ which vanish on $\mathfrak{t}_{\textup{h}}$,
is the maximal subgroup of $\mathfrak{t}^*$ which acts on the orbit space
$M/T$. We then proved a precise version of the following.

\begin{roughclassification}[Duistermaat--P.-, 2005] \label{rcl1}
Symplectic actions of tori on compact connected symplectic manifolds 
with coisotropic principal orbits 
are classified by an antisymmetric bilinear form $\sigma^{\mathfrak{t}}$ on $\mathfrak{t}$ (the restriction of $\sigma$ to the orbits), the Hamiltonian torus $T_{\textup{h}}$, 
the momentum polytope associated to $T_{\textup{h}}$ by the Atiyah--Guillemin--Sternberg theorem,
a discrete cocompact subgroup in $(\mathfrak{l}/\mathfrak{t}_{\textup{h}})^* \subset \mathfrak{t}^*$ (the
period lattice of $(\mathfrak{l}/\mathfrak{t}_{\textup{h}})^*$),
an antisymmetric bilinear form $c \colon (\mathfrak{l}/\mathfrak{t}_{\textup{h}})^* 
\times (\mathfrak{l}/\mathfrak{t}_{\textup{h}})^* \to \mathfrak{l}$ 
with certain integrality properties (represents the Chern class of  
the principal $T$\--bundle $M_{\textup{reg}} \to M_{\textup{reg}}/T$), and 
the holonomy invariant of a so called admissible connection for the 
principal bundle $M_{\textup{reg}} \to M_{\textup{reg}}/T$.
\end{roughclassification}

On the other hand suppose that $(M, \, \sigma)$ is 
a compact, connected $2n$\--dimensional symplectic manifold equipped with an effective 
action of a $(2n-2)$\--dimensional 
torus $T$ for which at least one $T$\--orbit is a $(2n-2)$\--dimensional 
symplectic submanifold of $(M, \, \sigma)$. Then the orbit space $M/T$ is a compact, connected,
smooth, orientable orbisurface ($2$\--dimensional orbifold) and the projection mapping
$\pi \colon M \to M/T$ is a smooth principal $T$\--bundle for which the collection 
$\{(\op{T}_x(T \cdot x))^{\sigma_x}\}_{x \in M}$ 
of symplectic orthogonal complements to the tangent
spaces to the $T$\--orbits is a flat connection.
Let $p_0$ be any regular point in $M/T$, 
$\pi^{\textup{orb}}_1(M/T, \, p_0)$ be the orbifold fundamental group, 
and $\widetilde{M/T}$ be the orbifold
universal cover of $M/T$.
Then the symplectic manifold $(M, \, \sigma)$ is $T$\--equivariantly symplectomorphic to the
associated bundle
$
\widetilde{M/T} \times_{\pi^{\textup{orb}}_1(M/T, \, p_0)} T,
$
where $\pi^{\textup{orb}}_1(M/T, \, p_0)$ acts on $T$ by means of the monodromy
homomorphism of the aforementioned flat connection.
We will describe the symplectic form on this space, as well as the torus action, in 
Definition \ref{themodel}. The $T$\--action on this bundle comes from the $T$\--action
on $\widetilde{M/T} \times T$ by translations on the right factor.
Then we prove a precise version of the following, c.f. Theorem \ref{TT}:

\begin{roughclassification} \label{rcl2}
Symplectic actions of $(2n-2)$\--dimensional tori on compact connected symplectic $2n$\--manifolds 
for which at least one $T$\--orbit is a $(2n-2)$\--dimensional symplectic submanifold
are classified by a non--degenerate  antisymmetric 
bilinear form $\sigma^{\mathfrak{t}}$ on $\mathfrak{t}$ (the restriction of $\sigma$ to the orbits),
the Fuchsian signature of the orbit
space $M/T$, which is a compact, connected orbisurface (and by this we mean the
genus $g$ of the underlying surface and the tuple of orders $\vec{o}$ of the orbifold singularities
of $M/T$),
the total symplectic area of $M/T$, and 
an element in $T^{2g+n}/\mathcal{G}$ which encodes the holonomy
of the aforementioned flat connection for $\pi \colon M \to M/T$, where $n$ is the number of orbifold singular points of $M/T$, and where $\mathcal{G}$ is the group of matrices
\begin{eqnarray} \label{G0}
\mathcal{G}:=\{  
\left( \begin{array}{cc}
A & \textup{0}  \\
C & D 
\end{array} \right) \in \textup{GL}(2g+n, \, \Z)
 \, | \, A \in \textup{Sp}(2g, \, \Z), \, \,  D \in \mathcal{M}\textup{S}^{\vec{o}}_n \}.
\end{eqnarray}
Here $\op{Sp}(2g, \, \Z)$ stands for symplectic matrices and $\mathcal{M}\textup{S}^{\vec{o}}_n$
for permutation matrices which preserve the tuple $\vec{o}$ of orbifold singularities of $M/T$.
\end{roughclassification}

Moreover we show that if the $T$\--action is free, then
$M$ is $T$\--equivariantly diffeomorphic to the product
$
M/T \times T,
$
equipped with the action of $T$ by translations on the right factor. 
Thus two such symplectic manifolds are equivariantly diffeomorphic if and only if their
corresponding orbit spaces are surfaces of the same genus. 
Using the rough classifications \ref{rcl1}, \ref{rcl2} as a stepping stone we obtain the following classification,
a precise and explicit statement of which is Theorem \ref{520} (pp. 59).

\begin{roughclassification}
Let $(M, \, \sigma)$ be a compact connected $4$\--dimensional symplectic manifold equipped with
an effective symplectic action of a $2$\--torus $T$. Then one and only one of the following
cases occurs:
\begin{itemize}
\item[\textup{1)}]
$(M, \, \sigma)$ is a $4$\--dimensional symplectic toric manifold, determined by
its Delzant polygon.
\item[\textup{2)}]
$(M, \, \sigma)$ is equivariantly symplectomorphic to
a product $\T^2 \times S^2$, where
$\T^2=(\R/\Z)^2$ and
the first factor of $\mathbb{T}^2$ acts on the left factor by translations
on one component, and the second factor acts on $S^2$ by rotations about the
vertical axis of $S^2$. The symplectic form is a positive linear combination
of the
standard translation invariant form on $\mathbb{T}^2$ and
the standard rotation invariant form on $S^2$.
\item[\textup{3)}]
$T$ acts freely on $(M, \, \sigma)$ with all $T$\--orbits being Lagrangian $2$\--tori, 
and it is 
classified in (rough) classification \ref{rcl1} by a 
discrete cocompact subgroup $P$ of 
$\got{t}^*$, an antisymmetric bilinear mapping 
$c \colon \got{t}^*\times \got{t}^* \to\got{t}$ which satisfies certain symmetry and integrality
properties, and the so called holonomy invariant of a so called admissible connection
for $M_{\textup{reg}} \to M_{\textup{reg}}/T$. 
\item[\textup{4)}]
$T$ acts on $(M, \, \sigma)$ with all $T$\--orbits being symplectic $2$\--tori, and
it is classified in (rough) classification \ref{rcl2} by an antisymmetric bilinear form 
$\sigma^{\mathfrak{t}}$ on 
$\mathfrak{t}$, the Fuchsian signature of $M/T$, the total symplectic area of $M/T$,
and an element in $T^{2n+g}/\mathcal{G}$, where $g$ is the genus of $M/T$,
$n$ is the number of singular points of $M/T$, and $\mathcal{G}$ is the group given
by $(\ref{G0})$.
\end{itemize}
\end{roughclassification}

The paper is organized as follows. In Section  2 we describe models of
$(M, \, \sigma)$ up to $T$\--equivariant symplectomorphisms and up to
$T$\--equivariant diffeomorphisms. 

In Section 3
we classify free
symplectic 
torus actions of a $(2n-2)$\--dimensional torus $T$ on 
$2n$\--dimensional symplectic manifolds, when at least one
$T$\--orbit is a $(2n-2)$\--dimensional
symplectic submanifold of $(M, \, \sigma)$, up to $T$\--equivariant 
symplectomorphisms, c.f. Theorem \ref{FTT}, and 
also up to $T$\--equivariant diffeomorphisms, c.f. Corollary \ref{FTTT}.

In Section 4 we extend the results in Section 3 to non--free actions. We only present those parts
of the proofs which are different from the proofs in Section 3, and 
hence we 
suggest that the paper be read linearly between Section 3 and Section 5, both included. In turn
this approach has the benefit that with virtually no repetition we are able to present
the classification in the free case in terms of invariants which are easier to describe than in the general case. 

In Section 5 we provide a classification of symplectic actions of $2$\--dimensional tori on
compact, connected, symplectic $4$\--dimensional manifolds. This generalizes the
$4$\--dimensional Delzant's theorem \cite{Delzant} to non--Hamiltonian actions.

There is extensive literature on the classification of Hamiltonian, symplectic and/or smooth 
torus actions. The papers 
closest to our paper in spirit are the paper \cite{Delzant} by Delzant on the classification of 
symplectic--toric manifolds, and the paper
by Duistermaat and the author \cite{DuPe} on the classification of symplectic torus actions
with coisotropic principal orbits.
The following are other contributions related to our work. 
The paper \cite{K} by Karshon on the classification
of Hamiltonian circle actions on compact connected symplectic $4$\--manifolds.
The book of Audin's \cite{audin}  on Hamiltonian torus 
actions, and Orlik--Raymond's \cite{OR} and Pao's \cite{pao} papers, on the classification of actions 
of $2$\--dimensional tori on $4$\--dimensional compact 
connected smooth manifolds -- they do not assume 
an invariant symplectic structure.
Kogan \cite{Ko}  studied completely integrable 
systems with local torus actions. Karshon and Tolman
studied centered complexity one Hamiltonian 
torus actions in arbitrary dimensions in their article \cite{KT}
and Hamiltonian torus actions with $2$\--dimensional symplectic quotients
in \cite{KT0}.
McDuff \cite{M} and McDuff and Salamon \cite{MS} studied 
non--Hamiltonian circle actions, and Ginzburg 
\cite{ginzburg} non--Hamiltonian symplectic actions of 
compact groups under the assumption of a ``Lefschetz 
condition''.  Symington \cite{s} and Leung and Symington \cite{ls} 
classified $4$\--dimensional compact connected symplectic 
manifolds which are fibered by Lagrangian tori where 
the fibration may have have elliptic or 
focus\--focus singularities. 
\\
\\
{\bf Acknowledgements.} 
The author is grateful to professor Y. Karshon for moral and intellectual support throughout this project, and for comments on preliminary versions, which
have enhanced the clarity and accuracy.

He is grateful to professor J.J. Duistermaat for discussions, specifically on sections 2.4.4, 2.5, 2.6, for hospitality on three visits to Utrecht, and for comments on a preliminary version.

He thanks professor A. Uribe for conversations on 
symplectic normal forms, and professor P. Scott for discussions on orbifold theory,
and helpful feedback and remarks on Section 4.2.
He also has benefited from conversations with
professors D. Auroux, D. Burns, P. Deligne, V. Guillemin, A. Hatcher, D. McDuff, M. Pinsonnault,
R. Spatzier and E. Zupunski. In particular, 
E. Zupunski sat through several talks of the author on the paper and offered feedback.
He thanks professor M. Symington for the suggestion of the problem\footnote{The question treated 
in this paper was raised by professor M. Symington and was communicated to the 
author by professor Y. Karshon.}. Additionally, he thanks professor D. McDuff and professor P. Deligne for the hospitality during visits to Stony Brook and to IAS in the Winter of 2006, to discuss the content of the article \cite{DuPe}, which influenced the presentation of some topics in the current article. He thanks
Oberlin College for the hospitality during the author's visit (October 2006-May 2007), while
supported by a Rackham Predoctoral Fellowship.

\section{Global model}

Unless otherwise stated we assume throughout the section that
$(M, \, \sigma)$ is a smooth compact and connected symplectic manifold and 
$T$ is a torus which acts effectively on $(M, \, \sigma)$ by means of
symplectomorphisms. We furthermore assume that at least one $T$\--orbit is a
$\op{dim}T$\--dimensional symplectic submanifold of $(M, \, \sigma)$.

\subsection{The orbit space $M/T$}

We describe the structure of the orbit space $M/T$, c.f. Definition \ref{maximal}.
\subsubsection{Symplectic form on the $T$\--orbits}

We prove that the symplectic form on every $T$\--orbit of $M$ is given by 
 the same non--degenerate
 antisymmetric bilinear form.

Let $X$ be an element of the Lie algebra $\mathfrak{t}$ 
of $T$, and denote by $X_M$ the smooth vector field on $M$ obtained as the infinitesimal 
action of $X$ on $M$. Let $\omega$ be a smooth differential form, let $\op{L}_v$ denote the Lie derivative with respect to a vector field $v$, and let $\op{i}_v\omega$ denote the usual inner product 
of  $\omega$ with $v$. Since the symplectic form $\sigma$ is $T$\--invariant, we have that
$
\op{d}(\textup{i}_{X_M}\sigma )=\op{L}_{X_M}\sigma =0,
$
where the first equality
follows by combining $\op{d} \sigma =0$ and the homotopy 
identity $\op{L}_v =\op{d}\circ\op{i}_v +\op{i}_v\circ\op{d}$.
The following result follows from \cite[Lem.\,2.1]{DuPe}.

\begin{lemma} \label{sigmat}
Let $(M, \,\sigma)$ be a compact, connected, symplectic 
manifold equipped with an effective symplectic action of a torus $T$ for which there is at least
one  $T$\--orbit which is a $\op{dim}T$\--dimensional symplectic submanifold of $(M, \, \sigma)$. 
Then there exists a unique non--degenerate antisymmetric bilinear form
$\sigma^{\mathfrak{t}} \colon \mathfrak{t} \times \mathfrak{t} \to \R$
on the Lie algebra $\mathfrak{t}$ of $T$
such that 
\begin{eqnarray} \label{abf}
\sigma_x(X_M(x), \,Y_M(x))=\sigma^{\mathfrak{t}}(X,Y),
\end{eqnarray}
for every $X, \,Y \in \mathfrak{t}$, and every $x \in M$.
\end{lemma}

\begin{proof}
In \cite[Lem.\,2.1]{DuPe} it was shown that there is a unique 
antisymmetric bilinear form $\sigma^{\mathfrak{t}} \colon \mathfrak{t} \times 
\mathfrak{t} \to \R$ on the Lie algebra $\mathfrak{t}$ of $T$ such that
expression (\ref{abf}) holds for every $X, \,Y \in \mathfrak{t}$, and every $x \in M$. 
We recall the proof which was given in \cite{DuPe}. It follows from 
Benoist\footnote{Although the observation was probably first done by Kostant. I thank Y. Karshon for making me aware of this.},
\cite[Lem.\,2.1]{Be} that if $u$ and $v$ are smooth vector fields on $M$ such that 
$\textup{L}_u\sigma=0$
and $\textup{L}_{v}\sigma=0$, then $[u,\,v]=\textup{Ham}_{\sigma(u,\,v)}$. Indeed,
observe that
\begin{eqnarray} \label{eqa}
\textup{i}_{[u,\,v]}\sigma=\textup{L}_u(\textup{i}_v 
\sigma)=\textup{i}_u(\op{d}(\textup{i}_v\sigma))
+\op{d}(\textup{i}_u(\textup{i}_v \sigma))=\op{d}(\sigma(u,\,v)).
\end{eqnarray}
Here we used $\textup{L}_u\sigma=0$ in the first equality, the homotopy formula for the 
Lie derivative
in the second identity, and finally $\op{d}\sigma=0$, the homotopy identity and 
$\textup{L}_v\sigma=0$
in the third equality. 
Applying (\ref{eqa}) to
$u=X_M$, $v=Y_M$, where $X,\,Y \in \mathfrak{t}$, we obtain that
\begin{eqnarray} \label{eqqq}
\textup{i}_{[X_M,\,Y_M]}\sigma=\op{d}(\sigma(X_M,\,Y_M)).
\end{eqnarray}
On the other hand, since $T$ is commutative, 
\begin{eqnarray} \label{eqb}
[X_M,\,Y_M]=-[X,\,Y]_M=0,
\end{eqnarray}
and hence $[X,\,Y]=0$. Combining expression (\ref{eqqq}) with expression (\ref{eqb}),
we obtain that the derivative of the real valued function $x \mapsto \sigma_x(X_M(x),
\,Y_M(x))$ identically vanishes on $M$, which in virtue of the connectedness of $M$ and the fact
that $\sigma$ is a symplectic form, implies 
expression (\ref{abf}) for a certain antisymmetric bilinear form $\sigma^{\mathfrak{t}}$
in $\mathfrak{t}$. 
Since there is a $T$\--orbit of dimension $\op{dim}T$ which
is a symplectic submanifold of $(M, \, \sigma)$, 
the form $\sigma^{\mathfrak{t}}$ must be non--degenerate.
\end{proof}

\begin{remark} \label{disremark}
Each tangent space $\op{T}_x(T \cdot x)$ equals the linear span of the vectors $X_M(x)$, $X \in \mathfrak{t}$.
The collection of tangent spaces $\op{T}_x(T \cdot x)$ to the $T$\--orbits $T \cdot x$ forms
a smooth $\op{dim}T$\--dimensional distribution\footnote{Since $T$ is a commutative group, the Lie
brackets of $X_M$ and $Y_M$ are zero for all $X, \, Y \in \mathfrak{t}$, which implies that 
$\mathcal{D}$ is integrable,  which in particular verifies the integrability theorem of Frobenius 
\cite[Th.\,1.60]{warner}
for our particular assumptions (although we do \emph{not} need it).}, which is integrable, where the integral manifold through $x$ is precisely the $T$\--orbit $T \cdot x$.
Since the $X_M$, $X \in \mathfrak{t}$, are $T$\--invariant vector fields, the distribution 
$H=\{\op{T}_x(T \cdot x)\}_{x \in M}$ is $T$\--invariant. Each element of 
$H$ is a symplectic vector space.
\end{remark}

\subsubsection{Stabilizer subgroup classification}

Recall that if $M$ is an arbitrary smooth manifold equipped with a smooth action of a torus $T$,
for each $x\in M$ we write $T_x:=\{ t\in T \, | \, t\cdot x=x\}$ 
for the \emph{stabilizer subgroup of the action of  $T$ on $M$ at the point $x$}. 
$T_x$ is a closed Lie subgroup of $T$. 
In this section we study the stabilizer subgroups
of the action of $T$ on $(M, \, \sigma)$.

In \cite[Sec.\,2]{DuPe}, Duistermaat and the author pointed out that for general
symplectic torus actions the stabilizer subgroups of the action
need not be connected, which is in contrast with the symplectic actions whose principal orbits are
Lagrangian submanifolds, where the stabilizer
subgroups are subtori of $T$; such fact also 
may be found as statement (1)(a) in Benoist's article 
\cite[Lem.\,6.7]{Be}. 

\begin{lemma} \label{sigmat1}
Let $T$ be a torus. Let $(M, \, \sigma)$ be a compact, connected, symplectic 
manifold equipped with an effective, symplectic action of the torus $T$, such that 
at least one $T$\--orbit is a $\op{dim}T$\--dimensional 
symplectic submanifold of $(M, \, \sigma)$. Then
the stabilizer subgroup of the $T$\--action at every point in $M$ is
a finite abelian group.
\end{lemma}

\begin{proof}
Let $\mathfrak{t}_x$ denote the Lie algebra of the stabilizer subgroup $T_x$ of the
action of $T$ on $M$ at the point $x$.
In the article of Duistermaat and the author \cite[Lem.\,2.2]{DuPe}, we observed that
for every $x \in M$ there is an inclusion
$\mathfrak{t}_x \subset \textup{ker}\, \sigma^{\mathfrak{t}}$, and since 
by Lemma \ref{sigmat} $\sigma^{\mathfrak{t}}$
is non--degenerate, its kernel $\textup{ker}\, \sigma^{\mathfrak{t}}$ is trivial, which
in turn implies that $\mathfrak{t}_x$ is the trivial vector space, and hence $T_x$, which is
a closed and hence compact subgroup of $T$, must be a finite group. 
\end{proof}

\begin{lemma} \label{finitestabilizers00}
Let $T$ be a torus. Let $(M, \, \sigma)$ be a compact, connected, symplectic 
manifold equipped with an effective, symplectic action of the torus $T$ such that 
at least one $T$\--orbit is a $\op{dim}T$\--dimensional 
symplectic submanifold of $(M, \, \sigma)$. Then
every $T$\--orbit is a $\op{dim}T$\--dimensional 
symplectic submanifold of $(M, \, \sigma)$.
\end{lemma}

\begin{proof}
Let $x \in M$. Since the torus $T$ is a compact group, the action of $T$ on the smooth manifold
$M$ is proper, and the mapping
\begin{eqnarray} \label{diffeomorphism}
t \mapsto t \cdot x \colon T/T_x \to T \cdot x
\end{eqnarray}
is a diffeomorphism, c.f. \cite[Appendix B]{GGK} or \cite[Sec.\,23.2]{AC}, and in particular, the dimension of quotient group $T/T_x$ equals the dimension of the $T$\--orbit $T \cdot x$. Since by
Lemma \ref{sigmat1} each
stabilizer subgroup $T_x$ is finite, the dimension of $T/T_x$ equals $\op{dim}T$, and
hence every $T$\--orbit is $\op{dim}T$\--dimensional. By Lemma \ref{sigmat} 
the symplectic form $\sigma$ restricted to any $T$\--orbit of the $T$\--action 
is non--degenerate and hence
$T\cdot x$ is a symplectic submanifold of $(M, \, \sigma)$.
\end{proof}

\begin{cor}
Let $(M, \,\sigma)$ be a compact connected symplectic manifold equipped
with an effective symplectic action of a torus $T$ for which at least one, and hence
every $T$\--orbit is 
a $\op{dim}T$\--dimensional symplectic submanifold of $(M, \, \sigma)$. 
Then there exists only finitely many different subgroups of $T$  which occur
as stabilizer subgroups of the action of $T$ on $M$, and each of them is a finite group.
\end{cor}

\begin{proof}
It follows from Lemma \ref{sigmat1} that every stabilizer subgroup of the
action of $T$ on $M$ is a finite group. 
It follows from  the tube theorem, c.f. \cite[Th.\,2.4.1]{DuKo} or \cite[Th.\,B24]{GGK} 
that in the case of a compact smooth manifold equipped
with an effective action of a torus $T$, there exists only 
finitely many different subgroups of $T$ which occur
as stabilizer subgroups.
\end{proof}

If $M$ is $4$\--dimensional and $T$ is $2$\--dimensional we have
the following stronger statement, 
which follows from the tube theorem, since a finite group acting linearly on a disk must
be a cyclic group acting by rotations.

\begin{lemma} \label{finitestabilizers}
Let $T$ be a $2$\--torus. Let $(M, \, \sigma)$ be a compact, connected, symplectic 
$4$\--manifold equipped with an effective, symplectic action of $T$, such that 
at least one, and hence every $T$\--orbit is a $2$\--dimensional symplectic
submanifold of $(M, \, \sigma)$. Then the stabilizer subgroup of the $T$\--action at every point in $M$ is
a cyclic abelian group.
\end{lemma}

\subsubsection{Orbifold structure of $M/T$}

We denote the space of all orbits in $M$ of the $T$\--action by $M/T$, and by 
$\pi \colon M \to M/T$ the canonical projection which assigns to each $x \in M$ the 
orbit $T \cdot x$ through the point $x$. The orbit space is provided with the maximal topology for which the canonical projection $\pi$ is continuous; this topology is Hausdorff.
 Because $M$ is compact and connected, the orbit space $M/T$ is compact and
 connected. For each connected component $C$ of an orbit type $M^H:=\{x \in M \, |\, T_x=H\}$ in $M$ 
of the subgroup $H$ of $T$, the action of $T$ on $C$ induces a 
proper and free action of the torus $T/H$ on $C$, and $\pi(C)$ 
has a unique structure of a smooth manifold such that 
$\pi  \colon C\to\pi(C)$ is a principal $T/H$\--bundle. 
The orbit 
space $M/T$ is not in general a smooth manifold, c.f. Example \ref{ex2}. 
Our next goal is to show that
the orbit space $M/T$ has a natural structure of
smooth orbifold.

\begin{example}[Free action] \label{ex1}
\normalfont
Let $(M, \, \sigma)$ be the Cartesian product 
 $(\R/\Z)^2 \times S^2$
equipped with the product symplectic form of the standard symplectic (area) form
on the torus $(\R/\Z)^2$ and the standard area form on the sphere $S^2$. 
Let $T$ be the $2$\--torus $(\R/\Z)^2$,
and let $T$ act on $M$ by translations on the left factor of the product. Such
action of $T$ on $M$ is free, it has symplectic $2$\--tori as $T$\--orbits,
and the orbit space $M/T$ is equal to the $2$\--sphere 
$S^2$. Probably the simplest example of a $4$\--dimensional symplectic manifold
equipped with a symplectic action of a $2$\--torus for which the torus orbits
are symplectic $2$\--dimensional tori is the $4$\--dimensional torus $(\R/\Z)^2 \times (\R/\Z)^2$ with the
standard symplectic form, on which the $2$\--dimensional torus $(\R/\Z)^2$
acts by multiplications on two of the copies of $\R/\Z$ inside
of $(\R/\Z)^4$. The orbit space is a $2$\--dimensional torus, so
a smooth manifold.
\end{example}

\begin{example}[Non--free action] \label{ex2}
\normalfont
Consider the Cartesian product $S^2 \times (\R/\Z)^2$ of the $2$\--sphere
and the $2$\--torus equipped with the product symplectic form of the standard symplectic (area) form
on the torus $(\R/\Z)^2$ and the standard area form on the sphere $S^2$. The $2$\--torus $(\R/\Z)^2$ acts freely by translations on the right factor
of the product $S^2 \times (\R/\Z)^2$.
Consider the action of the finite group $\Z/2\,\Z$ on $S^2$ which rotates
each point horizontally by $180$ degrees, and the action of  $\Z/2\,\Z$
on the $2$\--torus $(\R/\Z)^2$ given by the antipodal action on the first circle. 
The diagonal action of $\Z/2\,\Z$ on $S^2 \times (\R/\Z)^2$ is free
and hence the quotient space $S^2 \times_{\Z/2\,\Z} (\R/\Z)^2$ is a smooth manifold.
Let $(M, \, \sigma)$ be this associated bundle $S^2 
\times_{\Z/2\,\Z} (\R/\Z)^2$ with the symplectic form and $T$\--actions inherited from
the ones given in the product $S^2 \times (\R/\Z)^2$, where $T=(\R/\Z)^2$. The action
of $T$ on $M$ is not free 
and the $T$\--orbits are symplectic $2$\--dimensional tori.
The orbit space $M/T$ is equal to $S^2/(\Z/ 2\, \Z)$, which is a 
smooth orbifold with two singular points of order $2$, the South and North poles of $S^2$. 
\end{example}

There is no standard definition of orbifold, so
we do not use terminology in orbifolds without clarifying or introducing it.
Following largely but not entirely Boileau--Maillot--Porti \cite[Sec.\,2.1.1]{BMP}, we define orbifold.
We also borrow from ideas in Satake \cite{satake1}, \cite{satake2} and Thurston \cite{Thurston}.
Our definition of orbifold is 
close to that in Haefliger's paper \cite[Sec. 4]{H}. Unfortunately Haefliger's paper does
not appear to be so well known and is frequently not given proper credit. I thank
Y. Karshon for making me aware of this.

\begin{definition}  \label{orb}
A \emph{smooth $n$\--dimensional orbifold} $\mathcal{O}$ is a metrizable topological space $|\mathcal{O}|$
endowed with an equivalence class of orbifold atlases for $\mathcal{O}$.
An \emph{orbifold atlas for $\mathcal{O}$} is 
a collection $\{(U_i, \, \widetilde{U}_i, \, \phi_i, \, \Gamma_i)\}_{i \in I}$
where for each $i \in I$, $U_i$ is an open subset of $\mathcal{O}$,
$\widetilde{U}_i$ is an open and connected subset of $\R^{n}$,
$\phi_i \colon \widetilde{U}_i \to U_i$ is a continuous map, called an \emph{orbifold chart},
and $\Gamma_i$ is a finite group of diffeomorphisms of $\widetilde{U}_i$, satisfying:
\begin{itemize}
\item[i)]
the $U_i$'s cover $\mathcal{O}$,
\item[ii)]
each $\phi_i$ factors through a homeomorphism between $\widetilde{U}_i/\Gamma_i$ and $U_i$, and
\item[iii)]
the charts are compatible. This means that for each $x \in \widetilde{U}_i$
and $y \in \widetilde{U}_j$ with $\phi_i(x)=\phi_j(y)$, there is a diffeomorphism
$\psi$ between a neighborhood of $x$ and a neighborhood of $y$
such that $\phi_j(\psi(z))=\phi_i(z)$ for all $z$ in such neighborhood.
\end{itemize}
Two orbifold atlases are \emph{equivalent} if their union is an orbifold atlas. 
If $x \in U_i$, the \emph{local group $\Gamma_x$ of $\mathcal{O}$ at a point $x \in \mathcal{O}$} is the
isomorphism class of the stabilizer  of the action of $\Gamma_i$ on $\widetilde{U}_i$ at the point $\phi_i^{-1}(x)$. 
A point $x \in \mathcal{O}$ is \emph{regular} if $\Gamma_x$ is trivial, and \emph{singular} otherwise. The \emph{singular locus}
is the set $\Sigma_{\mathcal{O}}$ of singular points of $\mathcal{O}$.
We say that the orbifold $\mathcal{O}$
is compact (resp. connected) if the topological space
$|\mathcal{O}|$ is compact (resp. connected). An \emph{orientation for an orbifold atlas for $\mathcal{O}$}
is given by an orientation on each $\widetilde{U}_i$ which is preserved by every change
of chart map $\psi$ as in part iii) above. The orbifold $\mathcal{O}$ is \emph{orientable} if it has
an orientation. 
\end{definition}

\begin{remark}
One can replace metrizable in Definition \ref{orb} by Hausdorff which although is a weaker condition, it suffices for our purposes.
\end{remark}

By the tube theorem, c.f. \cite[Th.\,2.4.1]{DuKo}, \cite[Th.\,B24]{GGK}
for each $x \in M$ there exists a $T$\--invariant open neighborhood $U_x$ of the $T$\--orbit
 $T \cdot x$ and a $T$\--equivariant diffeomorphism $\Phi_x$ from $U_x$ onto the
 associated bundle $T \times_{T_x} D_x$, where $D_x$ is an open disk centered at the origin in
 $\R^k=\C^{k/2}$, where $k:=\op{dim}M-\op{dim}T$,
 and $T_x$ acts by linear transformations on $D_x$.
The action of $T$ on $T \times_{T_x} D_x$ is induced by the action of $T$ by translations on
the left factor of $T \times D_x$. 
 Because $\Phi_x$
is a $T$\--equivariant diffeomorphism, it induces a homeomorphism $\widehat{\Phi}_x$
 on the quotient $\widehat{\Phi}_x \colon D_x/T_x \to \pi(U_x)$, and there
 is a commutative diagram of the form
 \begin{eqnarray} \label{keydiagram}
\xymatrix{ 
T \times D_x \ar[r]^{\pi_x}      &  T  \times_{T_x} D_x \ar[d]^{p_x}  \ar[r]^{\Phi_x}  & U_x  \ar[d]^{\pi|_{U_x}} \\
 D_x \ar[u]^{i_x} \ar[r]^{\pi'_x} &  D_x/T_x   \ar[r]^{\widehat{\Phi}_x} &      \pi(U_x)},
\end{eqnarray} 
where $\pi_x$, $\pi'_x$, $p_x$ are the canonical projection maps, and $i_x$ is the inclusion map.
Let 
\begin{eqnarray}
\phi_x:= \widehat{\Phi}_x \circ \pi'_x. \label{orbimap}
\end{eqnarray}

 The word lift
 in Lemma \ref{mine} below is used in the sense that 
 \begin{eqnarray} \label{lift}
 \pi_{\Gamma'} \circ i' \circ F=
 f \circ \pi_{\Gamma} \circ i
 \end{eqnarray}
where $i \colon D \to T \times D$, $i \colon D' \to T \times D'$
are the inclusion maps and
$\pi_{\Gamma} \colon T \times D \to T \times_{\Gamma} D$,
$\pi_{\Gamma'} \colon T \times D' \to T \times_{\Gamma'} D'$
are the canonical projections.

 \begin{lemma} \label{mine}
 Let $T$ be a torus.
 Let $\Gamma,\,\Gamma'$ be finite subgroups of $T$ respectively acting linearly on $\Gamma, \, \Gamma'$\--invariant open subsets
 $D, \, D' \subset \mathbb{R}^m$. 
 Let $z \in D$, $z' \in D'$. Let $\Gamma, \, \Gamma'$ act on $T \times D, \, T \times D'$, respectively,
 by the diagonal action, giving rise to smooth manifolds $T \times_{\Gamma}D$ and $T\times_{\Gamma'}D'$ equipped with the $T$\--actions induced by the action of 
 $T$ by left translations on $T \times D$, $T \times D'$, respectively.
 Let $f \colon T \times_{\Gamma}D \to T \times_{\Gamma'}D'$
 be a $T$\--equivariant diffeomorphism such that $f(T \cdot [1,\,z]_{\Gamma})=T \cdot [1,\,z']_{\Gamma'}$. 
 Then there exist open
 neighborhoods $U \subset D$ of $z$, and $U' \subset D'$ of $z'$, 
 and a diffeomorphism
 $F \colon U \to U'$ which lifts $f$ as in (\ref{lift}) and such that $F(z)=z'$.
 \end{lemma}

 One can obtain Lemma \ref{mine} applying the idea of
 the proof method of \cite[Lem. 23]{KIZ} by replacing
 the mapping $f\colon \R^n/\Gamma \to U'/\Gamma'$ therein by 
 $f \colon T \times_{\Gamma}D \to T \times_{\Gamma'}D'$.

 \begin{prop} \label{atlaslemma}
 Let $T$ be a torus. Let $(M, \, \sigma)$ be a compact, connected,
 symplectic manifold equipped with an effective, symplectic
 action of the torus $T$ for which at least one, and hence every $T$\--orbit, 
 is a $\op{dim}T$\--dimensional
 symplectic submanifold of $(M, \, \sigma)$. Then the collection of charts
 \begin{eqnarray} \label{mtatlas}
\widehat{\mathcal{A}}:=\{(\pi(U_x), \, D_x, \, \phi_x, \, T_x)\}_{x \in M} 
\end{eqnarray} 
is an orbifold atlas for the orbit space $M/T$, where for each 
$x \in M$ the mapping $\phi_x$ is defined by expression (\ref{orbimap}),
and the mappings $\widehat{\Phi}_x, \, \pi'_x$ 
are defined by diagram (\ref{keydiagram}).
 \end{prop} 
 
 \begin{proof}
Because $x \in U_x$, we have that $\bigcup_{x \in M}U_x=M$, so the collection
 $\{\pi(U_x)\}_{x \in M}$ covers $M/T$. Since $U_x$ is open,
 $\pi(U_x)$ is open, for each $x \in M$. Let $k:={\op{dim}M-\op{dim}T}$.
By Lemma \ref{sigmat1} the stabilizer group $T_x$ is a finite group of diffeomorphisms.
The disks $D_x$, $x \in M$, given by the tube theorem, are open subsets of $\R^k$,
since $T_x$ is a $0$\--dimensional subgroup of $T$.
Because it is obtained as a composite of continuous maps, 
$\phi_x$ in (\ref{orbimap}) is continuous and it factors 
through $\widehat{\Phi}_x \colon D_x/T_x \to \pi(U_x)$, which is the
homeomorphism on the bottom part of the right square of diagram (\ref{keydiagram}).

It is left to show that the mappings $\phi_x, \, \phi_y$, where $x, \,y \in M$, are compatible on
their overlaps. Indeed, pick $z \in D_x,\, z' \in D_y$ and assume that 
\begin{eqnarray} \label{phixz}
\phi_x(z)=\phi_y(z').
\end{eqnarray}
Let $U_z$ be an open neighborhood of $z$ such that $\widehat{\Phi}_x(U_z/T_x)$
is contained in $\pi(U_x) \cap \pi(U_y)$. Let $U_{z'} \subset D_y$ be the unique
open subset of $D_y$ such that $(\widehat{\Phi}_y)^{-1}(\widehat{\Phi}_x(U_z/T_x))=U_{z'}/T_y$.
Then the composite map
\begin{eqnarray} \label{xy}
\widehat{\Psi}_{xy}:=(\widehat{\Phi}_y)^{-1} \circ \widehat{\Phi}_x \colon U_z/T_x
\to U_{z'}/T_y
\end{eqnarray}
is a homeomorphism which by (\ref{phixz}) satisfies
\begin{eqnarray} \label{zz'}
\widehat{\Psi}_{xy}([z]_{T_x})=[z']_{T_y}. 
\end{eqnarray}
Since by the tube theorem the mappings
$\Phi_x \colon T \times_{T_x} D_x \to U_x$
and $\Phi_y \colon T \times_{T_y} D_y \to U_y$ are $T$\--equivariant
diffeomorphisms, and
by definition of $p_x, \, p_y$ we have that $(p_x)^{-1}(U_z/T_x) =T \times _{T_x} U_z$ and 
$(p_y)^{-1}(U_{z'}/T_y)=T \times_{T_y}U_{z'}$, the
composite map 
\begin{eqnarray} \label{lemmamap}
\Psi_{xy}:=(\Phi_y)^{-1} \circ \Phi_x \colon T \times_{T_x}U_z  \to T \times_{T_y}U_{z'},
\end{eqnarray}
is a $T$\--equivariant diffeomorphism. 
By the commutativity of diagram (\ref{keydiagram}), the map
$\Psi_{xy}$ lifts the map $\widehat{\Psi}_{xy}$.
Then by (\ref{zz'}), we have that
\begin{eqnarray} \label{zz'2}
\Psi_{xy}( T \cdot [1,\,z]_{T_x})=T \cdot [1,\,z']_{T_y},
\end{eqnarray}
where $\pi_x \colon T \times D_x \to T \times_{T_x} D_x$, $\pi_y \colon T \times D_y \to
T \times_{T_y} D_y$ are the canonical projection maps.
Then
the map in (\ref{lemmamap}) is of the form in Lemma \ref{mine} and
we can use this lemma to conclude that $\Psi_{xy}$ lifts to a diffeomorphism 
$
\psi_{xy} \colon W_z \to W_{z'},
$
where $W_z \subset U_z \subset D_x$,
and $W_{z'} \subset U_{z'} \subset D_y$ are open neighborhoods of $z,\,z'$ respectively, and
\begin{eqnarray} \label{pl}
\psi_{xy}(z)=z'.
\end{eqnarray}
Because the map (\ref{lemmamap}) lifts the map (\ref{xy}),
the diffeomorphism $\psi_{xy}$ lifts the
restricted homeomorphism 
$
\widehat{\Psi}_{xy} \colon W_z/T_x \to W_{z'}/T_y
$
induced by (\ref{xy}).
Then by (\ref{orbimap}), 
\begin{eqnarray} \label{exchange}
\phi_x \circ \psi_{xy}=\phi_y
\end{eqnarray}
on $W_z$. Expressions (\ref{pl}) and (\ref{exchange}) and
precisely describe the 
compatibility condition of the charts $\phi_x, \, \phi_y$ on their overlaps.
\end{proof}

 \begin{definition}  \label{maximal}
 Let $T$ be a torus. Let $(M, \, \sigma)$ be a compact, connected,
 symplectic manifold equipped with an effective, symplectic
 action of $T$ for which at least one, and hence every
 $T$\--orbit is a $\op{dim}T$\--dimensional symplectic submanifold of $(M, \, \sigma)$. 
 We call $\mathcal{A}$ the class of atlases
 equivalent to the orbifold atlas  $\widehat{\mathcal{A}}$ 
 defined by expression (\ref{mtatlas}) in Proposition \ref{atlaslemma}. 
 We denote the orbifold $M/T$ endowed with 
the class $\mathcal{A}$
by $M/T$, and the class $\mathcal{A}$ is assumed. 
 \end{definition}

 \begin{remark} \label{manifoldcase}
Since $M$ is compact and connected, $M/T$ is compact and connected. 
If $\op{dim}T=2n-2$, $M/T$ is a compact, connected orbisurface.
If the $T$ action is free, then the local 
groups $T_x$ in Definition \ref{maximal} are all trivial, and $M/T$
is a compact connected surface determined up to diffeomorphism by a non--negative integer, 
its topological genus.
 
Moreover, since every 
$T$\--orbit is a $\op{dim}T$\--dimensional symplectic submanifold of $(M, \sigma)$, 
the action of 
the torus $T$ on $M$ is a locally free and non--Hamiltonian action. Indeed,
by Lemma \ref{sigmat1}, the stabilizers are finite, hence discrete 
groups. On the other hand, if an action is Hamiltonian
the $T$\--orbits are isotropic submanifolds, so the $T$\--action cannot be Hamiltonian.
\end{remark}

\subsection{A flat connection for the projection $M \to M/T$}

In this section we prove that the projection $\pi \colon M \to M/T$
onto the orbifold $M/T$ (c.f. Proposition \ref{atlaslemma} and Definition \ref{maximal})
is a smooth principal $T$\--orbibundle. For such orbibundle there
are notions of connection, and of flat
connection, which extend the classical definition for bundles.
Then we show that $\pi \colon M \to M/T$
comes endowed with a flat connection (c.f. Proposition \ref{atlaslemma}). 

Let $\mathcal{O}, \, \mathcal{O'}$ be smooth orbifolds.
An \emph{orbifold diffeomorphism
$f \colon \mathcal{O} \to \mathcal{O}'$} is a is a homeomorphism
at the level of topological spaces $|\mathcal{O}|
\to |\mathcal{O}'|$ such that for every $x \in \mathcal{O}$
there are charts $\phi_i \colon \widetilde{U}_i \to U_i$,
$x \in U_i$, and $\phi'_j \colon \widetilde{U}'_j
\to U'_j$ such that $f(U_i) \subset U'_j$ and the
restriction $f|_{U_i}$ may be
lifted to a diffeomorphism $\widetilde{f} \colon
\widetilde{U}_i \to \widetilde{U}'_j$ which
is equivariant with respect to some homomorphism
$\Gamma_i \to \Gamma'_j$.

\begin{definition} \label{ic}
Let $T$ be a torus, let $(X, \, \omega)$ be a compact, connected symplectic manifold 
endowed with an effective symplectic action of $T$ for which at least one, and hence
every $T$\--orbit is a $\op{dim}T$\--dimensional symplectic submanifold of $(X, \, \omega)$.
Equip $X/T$ with the orbifold structure in 
Definition \ref{maximal}. Let $Y$ be a smooth manifold equipped with a smooth effective action of $T$.

A continuous surjective map $p \colon Y \to X/T$
is a \emph{smooth principal $T$\--orbibundle} if
the action of $T$ on $Y$ preserves the fibers of $p$ and acts freely and 
transitively on them, and for every $z \in X/T$ the following holds.
Let $\pi(U_x), \, D_x, \, \phi_x, \,T_x$ as in (\ref{mtatlas}). 
For each $x \in X$ for which $z \in \pi(U_x)$
there exists a $T_x$\--invariant open
subset $\widetilde{U}_x$ of $D_x$ and
a continuous map $\psi_x \colon T \times \widetilde{U}_x \to Y$ 
which induces a diffeomorphism between
$T  \times_{T_x} \widetilde{U}_x$ and $p^{-1}(\phi_x(\widetilde{U}_x))$
such that $p \circ \psi_x = \phi_x \circ \pi_1$, where $\pi_1 \colon
T \times \widetilde{U}_x\to \widetilde{U}_x$ is the canonical
projection. Here $T_x$ acts on $T \times \widetilde{U}_x$ by the diagonal action.

A \emph{connection for $p \colon Y \to X/T$} is a smooth vector subbundle $H$ of the tangent bundle
$\op{T} Y$ with the property that for each $y \in Y$,
$H_y$ is a direct complement in $\op{T}_y Y$ of the tangent space
to the fiber of $p$ that passes through $y$.
We say that the connection $H$ is \emph{flat with respect to $p \colon Y
\to X/T$} if the subbundle $H \subset \op{T}Y$ is integrable
considered as a smooth distribution on $Y$. 
\end{definition}

\begin{remark} \label{pullback}
Let $p \colon Y \to X/T$
be a smooth principal $T$\--orbibundle as in Definition \ref{ic}.
Then for every $y \in Y$
there are charts $\phi_i \colon \widetilde{U}_i \to U_i$ of $Y$,
$y \in U_i$, and $\phi'_j \colon \widetilde{U}'_j
\to U'_j$ of $X/T$, such that $p(U_i) \subset U'_j$ and the
restriction $p|_{U_i}$ may be
lifted to a smooth map $\widetilde{p} \colon
\widetilde{U}_i \to \widetilde{U}'_j$.
 
If the local group $T_x$ is trivial
for all $x \in X/T$, the action of $T$ on $X$ is trivial and $X/T$ is a smooth manifold.
Then the mapping $p \colon X
\to X/T$ is a smooth principal $T$\--bundle, in the usual sense of
the theory of fiber bundles on manifolds, for example c.f. \cite{Steenrod}. 
We also use the term $T$\--bundle instead of  $T$\--orbibundle.
\end{remark}

\begin{prop} \label{easynow}
Let $(M, \,\sigma)$ be a compact connected symplectic manifold equipped
with an effective symplectic action of a torus $T$ for which at least one, and hence every 
$T$\--orbit is a $\op{dim}T$\--dimensional symplectic submanifold of $(M, \, \sigma)$. 
Then the collection $\Omega=\{\Omega_x\}_{x \in M}$ of subspaces 
$\Omega_x \subset \op{T}_xM$,
where $\Omega_x$ is the $\sigma_x$\--orthogonal complement to $\op{T}_x(T \cdot x)$ in 
$\op{T}_xM$,
for every $x \in M$, is a smooth distribution on $M$. The projection
mapping $\pi \colon M \to M/T$ is a smooth principal $T$\--bundle
of which $\Omega$ is a $T$\--invariant flat connection.
\end{prop}

\begin{proof}
For each $x  \in M$ let $\widetilde{U}_x=D_x$, where $\widetilde{U_x}, \,D_x$ are
as in Definition Definition \ref{ic}.  Let 
$$
\psi_x(t,\,z):=\Phi_x([t,\,z]_{T_x}),
$$
where $\psi_x$ is also the map in the definition.
With these choices of $\widetilde{U_x}, \psi_x$, by diagram (\ref{keydiagram}) and
by construction of the orbifold atlas on $M/T$, c.f. 
Proposition \ref{atlaslemma} and Definition \ref{maximal}, the projection mapping 
$\pi \colon M \to M/T$ is a smooth principal $T$\--orbibundle. 
Because the tangent space to each $T$\--orbit $(\op{T}_x(T \cdot x), \, \sigma|_{\op{T}_x(T \cdot x)})$
is a symplectic vector space, its symplectic orthogonal complement 
$(\Omega_x, \, \sigma|_{\Omega_x})$ 
is a symplectic vector space.
Here $\sigma|_{\Omega_x}, \, \sigma|_{\op{T}_x(T \cdot x)}$ 
are the symplectic forms respectively induced by $\sigma$ 
on $\Omega_x, \, \op{T}_x(T \cdot x)$. Consider the disk bundle
$T \times_{T_x} \Omega_x$ where $T_x$ 
acts by induced linearized action on $\Omega_x$, and on $T \times \Omega_x$
by the diagonal action. The translational
action of $T$ on the left factor of $T \times \Omega_x$ 
descends to an action of $T$ on $T \times_{T_x} \Omega_x$. 
There exists a unique $T$\--invariant symplectic form $\sigma'$ on $T \times_{T_x} \Omega_x$
such that if $\pi' \colon T \times \Omega_x \to T \times_{T_x} \Omega_x$ is the canonical projection,
\begin{eqnarray} \label{cr}
\pi'^* \sigma'=\sigma|_{\op{T}_x(T \cdot x)} \oplus \sigma|_{\Omega_x},
\end{eqnarray}
where $\sigma|_{\op{T}_x(T \cdot x)} \oplus \sigma|_{\Omega_x}$ denotes the product 
symplectic form on $T \times \Omega_x$.
Then by the symplectic tube theorem 
of Benoist \cite[Prop.\,1.9]{Be}, Ortega and Ratiu \cite{ortegaratiu}
and Weinstein, which we use as it was formulated
in \cite[Sec.\,11]{DuPe}, there exists an open 
neighborhood $D$ of $0 \in \Omega_x$,  
an open $T$\--invariant neighborhood $U$ of $x$ in $M$,
and a $T$\--equivariant symplectomorphism $\Phi \colon T \times_{T_x}D \to U$
with $\Phi ([1,\, 0]_{T_x})=x$.
By $T$\--equivariance, $\Phi$ maps
the zero section of $T \times_{T_x} \Omega_x$ to
the $T$\--orbit $T \cdot x$ through $x$. It follows from (\ref{cr})
that the symplectic--orthogonal complement to such section in $T \times_{T_x} \Omega_x$ is
precisely the $(\op{dim}M-\op{dim}T)$\--dimensional manifold $\pi'(\{1\} \times \Omega_x)$, and hence, 
since $\Phi$ is a $T$\--equivariant symplectomorphism, the image 
$\Phi(\pi'( \{1\} \times \Omega_x))$ is an integral manifold through $x$ of dimension 
$\op{dim}M-\op{dim}T$, so 
 $\{\Omega_x\}_{x \in M}$ is a $T$\--invariant
 integrable distribution.
\end{proof}

\normalfont

We can formulate Proposition \ref{easynow} in the language of foliations as follows.

\begin{cor} \label{easynow2}
Let $(M, \,\sigma)$ be a compact connected symplectic manifold equipped
with an effective symplectic action of a torus $T$ for which at least one, and hence every 
$T$\--orbit is a $\op{dim}T$\--dimensional symplectic submanifold of $(M, \, \sigma)$. 
Then the collection of integral manifolds to the symplectic orthogonal complements
to the tangent spaces to the $T$\--orbits, c.f. Proposition \ref{easynow}, is a 
smooth $T$\--invariant $(\op{dim}M-\op{dim}T)$\--dimensional foliation of $M$.
\end{cor}

\subsection{Symplectic Tube Theorem}

It is convenient to state the the symplectic tube theorem used in the proof of Proposition
\ref{easynow} in a more standard way by replacing $\Omega_x$ by
$\C^m$ and $\sigma^{\Omega_x}$ by $\sigma^{\C^m}$, respectively, where
$m=1/2(\op{dim}M-\op{dim}T)$.
Indeed,
write $\op{i}:=\sqrt{-1} \in \mathbb{C}$ and let $\sigma ^{\C ^m}$ be the symplectic form on 
$\C ^m$
\begin{equation}
\sigma ^{\C ^m}
:=\sum_{j=1}^m\,\op{d}\overline{z^j}\wedge
\op{d} z^j/2\textup{i}. 
\label{sigmaW}
\end{equation}

Because $\Omega_x$ is a symplectic vector space, it has a symplectic basis with $2m$ elements
which induces a direct sum decomposition of $\Omega_x$ into $m$ mutually 
$\sigma|_{\Omega_x}$\--orthogonal 
two--dimensional linear subspaces $E_j$. The stabilizer group $T_x$ acts 
by means of symplectic linear transformations on the symplectic vector space $\Omega_x$, and
the aforementioned symplectic basis can be chosen so that $E_j$ is $T_x$\--invariant.
Averaging any inner product 
on each $E_j$ over $T_x$, we obtain an $T_x$\--invariant inner 
product $\beta _j$ on $E_j$, which
is unique if we also require that the symplectic 
inner product of any orthonormal basis with respect to 
$\sigma_{\Omega_x}$ is equal to $\pm 1$. This leads to 
the existence of a unique complex structure on 
$E_j$ such that, for any unit vector 
$e_j$ in $(E_j,\,\beta _j)$, we have that 
$e_j$, $\textup{i}\, e_j$ is an orthonormal basis 
in $(E_j,\,\beta _j )$ and $\sigma|_{\Omega_x}(e_j,\,\textup{i}\, e_j)=1$. 
This leads to an identification of $E_j$ with $\C$, and hence of
$\Omega_x$ with $\C^ m$, 
with the symplectic form defined by (\ref{sigmaW}).
The element $c \in\mathbb{T}^m$ acts 
on $\C^m$ by sending $z \in \C^m$ to the element 
$c\cdot z$ such that $(c\cdot z)^j=c^j\, z^j$ for every 
$1\leq j\leq m$. There is a unique monomorphism
of Lie groups $\iota :T_x\to \mathbb{T}^m$ such that 
$h\in T_x$ acts on $\Omega_x=\C^m$ by sending 
$z\in \C^m$ to $\iota (h)\cdot z$, 
hence $T_x$ acts on $T \times \C^m$ by 
\begin{eqnarray} \label{act}
h \star (t,\,z)=(h^{-1}t,\,\iota(h)\,z).
\end{eqnarray}

Consider the disk bundle $T \times_{T_x} \C^m$ where $T_x$ acts by (\ref{act}). The translational
action of $T$ on $T \times \C^m$ descends to an action of $T$ on $T \times_{T_x} \C^m$. 
By Lemma \ref{sigmat}, the antisymmetric bilinear form
$\sigma^{\mathfrak{t}} \colon \mathfrak{t} \times \mathfrak{t} \to \R$
is non--degenerate and hence it determines a unique symplectic
form $\sigma^T$ on $T$. In view of this, $\sigma|_{\op{T}_x(T \cdot x)}=\sigma^{\got{t}}$
in the proof of Proposition \ref{easynow} does not depend on $x \in M$.
The product symplectic form $\sigma^T \oplus \sigma^{\C^m}$ descends to a symplectic form on $T \times_{T_x}D$. With this terminology the proof of Proposition \ref{easynow} implies 
the following.

\begin{cor}[Tube theorem for symplectic orbits] \label{tubetheorem2}
Let $(M, \,\sigma)$ be a compact, connected, symplectic 
manifold equipped with an effective symplectic action of a torus $T$ for which at least
one, and hence every  $T$\--orbit is a $\op{dim}T$\--dimensional symplectic submanifold of $(M, \, \sigma)$. 
Then there is an open $\mathbb{T}^m$\--invariant 
neighborhood $D$ of $0 \in \C ^m$,  
an open $T$\--invariant neighborhood $U$ of $x$ in $M$,
and a $T$\--equivariant symplectomorphism 
$
\Phi \colon T \times_{T_x}D \to U
$
with $\Phi ([1,\, 0]_{T_x})=x$.
\end{cor}

In Corollary \ref{tubetheorem2} the product symplectic form is defined pointwise as
$$
(\sigma^T \oplus \sigma^{\C^m})_{(t, \, z)}((X,\,u),\, (X',\,u'))=\sigma ^{\mathfrak{t}}(X,\, X')+\sigma ^{\C ^m}(u,\,u').
\nonumber
$$
Here we identify each tangent space of the torus $T$ 
with the Lie algebra $\mathfrak{t}$ of $T$ and each tangent space of a vector space with the 
vector space itself.

\subsection{From Local to Global} 

We prove the main theorem 
of Section 2, Theorem \ref{tt}, which gives a model description
of $(M,\, \sigma)$ up to $T$\--equivariant
symplectomorphisms.

\subsubsection{Orbifold coverings of $M/T$}

We start by recalling the notion of orbifold covering.

\begin{definition}\cite[Sec.\,2.2]{BMP} \label{orbifoldcovering} 
A \emph{covering of a connected orbifold $\mathcal{O}$} is a connected orbifold
$\widehat{\mathcal{O}}$ together with a continuous mapping 
$p \colon \widehat{\mathcal{O}} \to \mathcal{O}$,
called an \emph{orbifold covering map}, such that every point
$x \in \mathcal{O}$ has an open neighborhood $U$ with the property that
for each component $V$ of $p^{-1}(U)$ there is a chart
$\phi \colon \widetilde{V} \to V$ of $\widehat{\mathcal{O}}$ such that
$p \circ \phi$ is a chart of $\mathcal{O}$.
Two coverings $p_1 \colon \mathcal{O}_1 \to
\mathcal{O}, \, \, p_2 \colon \mathcal{O}_2 \to \mathcal{O}$
are \emph{equivalent} if there exists a diffeomorphism $f \colon
\mathcal{O}_1 \to \mathcal{O}_2$ such that $p_2\circ f=p_1$.
A \emph{universal cover of $\mathcal{O}$} is a covering 
$p \colon \widehat{\mathcal{O}} \to \mathcal{O}$ such that for every
covering $q \colon \widetilde{\mathcal{O}} \to \mathcal{O}$, there
exists a unique covering $r \colon \widehat{\mathcal{O}} \to \widetilde{\mathcal{O}}$
such that $q \circ r=p$. 
The \emph{deck transformation group of a covering $p \colon
\mathcal{O}' \to \mathcal{O}$} is the group of all self--diffeomorphisms
$f \colon \mathcal{O}' \to \mathcal{O}'$ such that $p \circ f=p$,
and it is denoted by $\textup{Aut}(\mathcal{O}',\,p)$. 
 \end{definition}

\begin{remark} \label{localdiffeo}
Assuming the terminology introduced in Definition \ref{orbifoldcovering}, the
open sets $U,\,V$ equipped with the restrictions of the charts for $\widehat{\mathcal{O}},
\,\mathcal{O}$, are smooth orbifolds, and the restriction $p \colon V \to U$ is
an orbifold diffeomorphism.
\end{remark}

\begin{lemma} \label{obvious}
Let $(M, \,\sigma)$ be a compact connected symplectic manifold equipped
with an effective symplectic action of a torus $T$ for which at least one,
and hence every $T$\--orbit is a $\op{dim}T$\--dimensional
symplectic submanifold of $M$, and let $x \in M$. 
Let $\pi \colon M \to M/T$ be the canonical projection, and let
$\mathcal{I}_x$ be the maximal integral manifold of the
distribution of symplectic orthogonal complements to the tangent spaces to the
$T$\--orbits  which goes through $x$ (c.f. Proposition \ref{easynow}).
Then the inclusion $i_x \colon \mathcal{I}_x \to M$ is an injective immersion 
between smooth manifolds and 
the composite $\pi \circ i_x \colon \mathcal{I}_x \to M/T$ is an
orbifold covering map. 
\end{lemma}

\begin{proof}
Since $\Omega$ is a smooth distribution, the maximal integral manifold $\mathcal{I}_x$ in 
injectively immersed in $M$, c.f. \cite[Sec.\,1]{warner}.
Let $p \in M/T$, let $x \in \pi^{-1}(\{p\})$. 
Then $(\pi \circ i_x)^{-1}(\{p\})=T \cdot x \cap \mathcal{I}_x$.
By Corollary \ref{tubetheorem2}
there is an open $\mathbb{T}^m$\--invariant 
neighborhood $D_x$ of $0 \in \C ^m$,  
an open $T$\--invariant neighborhood $U_x$ of $x$ in $M$,
and a $T$\--equivariant symplectomorphism $\Phi_x \colon T \times_{T_x}D_x \to U_x$
with $\Phi_x([1,\, 0]_{T_x})=x$. 
For each $t \in T$, let the mapping $\rho_x(t) \colon D_x \to T \times_{T_x}D_x$ be given by
\begin{eqnarray} \label{rhotx}
\rho_x(t)(z)=[t,\,z]_{T_x}.
\end{eqnarray}
Since $(\pi \circ i_x)^{-1}(\{p\}) \subset T \cdot x$, there exists a collection
$\mathcal{C}:=\{t_k\} \subset T$ such that $(\pi \circ i_x)^{-1}(\{p\})=\{t_k \cdot x \, | \,t_k \in \mathcal{C}\}$.
Since $D_x$ is an open neighborhood of $0 \in D_x$ and $\Phi_x$ is a $T$\--equivariant
symplectomorphism, the image $\Phi_x(\op{Im}(\rho_x(t_k)))$ is an open neighborhood of
$t_k \cdot x$ in $\mathcal{I}_x$, and the image
$
V_p:=\pi(U_x \cap \mathcal{I}_x)
$
is an open neighborhood of $p \in M/T$. 
Let
$
V(t_k \cdot x):=\Phi_x(\op{Im}(\rho_x(t_k))$.
For each $t_k \in \mathcal{C}$, the set $V(t_k \cdot x)$ is a connected subset of
$\mathcal{I}_x$, since $\Phi_x$ and $\rho_x$ are continuous. Each connected 
component $K$ of $U_x \cap \mathcal{I}_x$ is of the form
$K=V(t_k \cdot x)$  for a unique $t_k \in \mathcal{C}$.
By the commutativity of diagram (\ref{keydiagram}), the composite mapping
\begin{eqnarray} \label{phirhox}
\Phi_x \circ \rho_x(t_k) \colon D_x \to V(t_k \cdot x) \subset \mathcal{I}_x
\end{eqnarray}
is a homeomorphism and hence a chart for $\mathcal{I}_x$ around $t_k \cdot x$, 
for each $t_k \in \mathcal{C}$. Since for 
each $t_k \in \mathcal{C}$ we have that $\pi(V(t_k \cdot x))=V_p$, the mapping
$$
\psi_x:=(\pi \circ i_x) \circ (\Phi_x \circ \rho_x(1)) \colon D_x \to V_p
$$
is onto. Moreover, $\psi_x$ is
smooth and factors through the homeomorphism $\widehat{\Phi}_x \colon D_x/T_x \to V_p$, and hence it is an orbifold chart for $M/T$ around $p$.
The definitions of the maps in (\ref{keydiagram}) 
imply that $\pi(V(t_k \cdot x))=V_p$ for all $t_k \in \mathcal{C}$.
The $T$\--equivariance of $\Phi_x$ and (\ref{rhotx}) imply that
$
(\pi \circ i_x) \circ (\Phi_x \circ \rho_x(t_k))=\psi_x.
$
for all $t_k \in \mathcal{C}$.
By taking the neighborhood $V_p$ around $p$ the map $\pi \circ i_x$ is
an orbifold covering as in Definition \ref{orbifoldcovering}, where therein we take
$U:=V_p$, and the charts in (\ref{phirhox}) as charts for $\mathcal{I}_x$.
\end{proof}

We will see in Section 2.4.2. that the map $\pi \circ i_x$ in Lemma \ref{obvious}
respects the symplectic forms, where the form on $\mathcal{I}_x$ is the restriction 
of $\sigma$, and the form on $M/T$ is the natural symplectic form that we define therein.

\begin{remark} \label{re1}
Let $M$ be an arbitrary smooth manifold.
It is a basic fact of foliation theory \cite[Sec.\,1] {warner} that in general, the integral 
manifolds of a smooth distribution $D$ on $M$ are injectively immersed manifolds 
in $M$, but they are not necessarily embedded or compact. For example 
the one--parameter 
subgroup of tori $\{ (t, \, \lambda \, t)+\Z^2 \in \R^2/\Z^2 \, | \, 
 t \in \R \}$, 
in which the constant $\lambda$ is an irrational real number. This is 
the maximal integral manifold through $(0, \, 0)$ of the distribution 
which is spanned by the constant vector field $(1,\,  \lambda)$, and it
is non--compact. 

It is an exercise to verify that if $I_x$ is a maximal integral manifold of a
smooth distribution $D$ which passes 
through $x$, then $I_x$ must contain every end point of a smooth curve 
$\gamma$ which starts at $x$ and satisfies the condition that for each 
$t$ its velocity vector $\op{d}\gamma(t)/\op{d}t$ belongs to $D_{\gamma(t)}$, 
and conversely each such end point is 
contained in an integral manifold through $x$. Therefore 
$I_x$ is the set of all such endpoints, which is the unique maximal 
integral manifold through $x$. It remains is to show that 
this set $I_x$ is an injectively immersed manifold in $M$.
This is one of the basic results of foliation theory, c.f. \cite[Sec.\,1]{warner}. 
In Remark \ref{e2} we give a self--contained proof of this fact for the
particular case of our orbibundle $\pi \colon M \to M/T$ in Proposition \ref{easynow}.
\end{remark}

Next we define a notion of orbifold fundamental group 
that extends the classical definition for when the orbifold is a manifold. The
first difficulty is to define a ``loop in an orbifold'', which
we do mostly but not entirely following
Boileau--Maillot--Porti \cite[Sec.\,2.2.1]{BMP}. 
The definitions for orbifold loop and homotopy of loops in an orbifold which we give
are not the most general ones, but give a convenient definition of
orbifold fundamental group, which we use in the definition of
the model for $(M, \, \sigma)$ in Definition \ref{themodel}.

\begin{definition} \label{equalpaths}
An \emph{orbifold loop} $\alpha \colon [0,\,1] \to \mathcal{O}$ 
in a smooth orbifold $\mathcal{O}$ is represented by:
\begin{itemize}
\item
a continuous map $\alpha  \colon [0, \,1] \to |\mathcal{O}|$
such that $\alpha(0)=\alpha(1)$ and there are at most finitely many $t$ such that
$\alpha(t)$ is a singular point of $\mathcal{O}$,
\item
for each $t$ such that $\alpha(t)$ is singular, a chart
$\phi \colon \widetilde{U} \to U$, $\alpha(t) \in U$,
a neighborhood $V(t)$ of $t$ in $[0,\,1]$ such that for all
$u \in V(t) \setminus \{t\}$,
$\alpha(u)$ is regular and lies in $U$, and a lift
$\widetilde{\alpha}|_{V(t)}$
of $\alpha|_{V(t)}$ to $\widetilde{U}$. We say that $\widetilde{\alpha}|_{V(t)}$ is
\emph{a local lift of $\alpha$ around $t$}.
\end{itemize}
We say that two orbifold loops $\alpha, \, \alpha' \colon [0,\,1] \to \mathcal{O}$ 
respectively equipped with lifts
of charts $\widetilde{\alpha}_i, \, \widetilde{\alpha}_j'$ 
\emph{represent the same loop} if the underlying maps are equal and
the collections of charts satisfy: for each $t$ such that 
$\alpha'(t)=\alpha(t) \in U_i \cap U'_j$ is singular, 
where $U_i, \, U'_j$ are the corresponding charts 
associated to $t$,
there is a diffeomorphism
$\psi$ between a neighborhood of 
$\widetilde{\alpha}_i(t)$ and a neighborhood of $\widetilde{\alpha}'_j(t)$
such that $\psi(\widetilde{\alpha}_i(s))=\widetilde{\alpha}'_j(s)$ for all $s$
in a neighborhood of $t$ and $\phi'_j \circ \psi=\phi_i$.
\end{definition}

In Definition \ref{equalpaths}, if $\mathcal{O}$ is a smooth manifold, two orbifold loops in $\mathcal{O}$
are \emph{equal} if their underlying maps are equal.

\begin{definition} \label{homotopic}
Let $\gamma, \, \lambda \colon [0, \, 1] \to \mathcal{O}$
be orbifold loops with common initial and end point in a smooth orbifold $\mathcal{O}$ of which
the set of singular points has codimension at least $2$ in $\mathcal{O}$.
We say that \emph{$\gamma$ is homotopic to $\lambda$ with fixed end points} if
there exists a continuous map $H  \colon [0, \,1] \to |\mathcal{O}|$
such that:
\begin{itemize}
\item
for each $(t,\,s)$ such that $H(t, \,s)$ is singular, there is a chart
$\phi \colon \widetilde{U} \to U$, $H(t, \,s) \in U$,
a neighborhood $V(t, \,s)$ of $(t, \,s)$ in $[0,\,1]^2$ such that for all
$(u,\,v) \in V(t,\,s) \setminus \{(t,\,s)\}$,
$H(u,\,v)$ is regular and lies in $U$, and a lift
$\widetilde{H}|_{V(t,\,s)}$ of
$H|_{V(t,\,s)}$ to $\widetilde{U}$. We call  $\widetilde{H}|_{V(t,\,s)}$
\emph{a local homotopy lift of $H$ around $(t, \,s)$}.
\item
the orbifold loops $t \mapsto H(t,\,0)$ and $t \mapsto H(t,\,1)$ from
$[0,\,1]$ into $\mathcal{O}$ respectively
endowed with the local homotopy
lifts $\widetilde{H}|_{V(t,\,0)}$, $\widetilde{H}|_{V(t,\,1)}$ of $H$ for each
$t \in [0,\,1]$, are respectively equal to $\gamma$ and $\lambda$ as orbifold loops, c.f. Definition \ref{equalpaths}.
\item
$H$ fixes the initial and end point: $H(0,\,s)=H(1,\,s)=\gamma(0)$ for all $s \in [0,\,1]$.
\end{itemize}
\end{definition}

The assumption in Definition \ref{homotopic} on the codimension of the singular
locus of the orbifold always holds for symplectic orbifolds, so such requirement is
natural in our context, since the orbifold which we will be working with, the orbit space $M/T$,
comes endowed with a symplectic structure. From now on we assume this 
requirement for all orbifolds.

\begin{definition}\cite[Sec.\,2.2.1]{BMP} \label{orbifoldpi1}
Let $\mathcal{O}$ be a connected orbifold. The \emph{orbifold
fundamental group $\pi^{\textup{orb}}_1(\mathcal{O},
\, x_0)$ of $\mathcal{O}$ based at the point $x_0$}
as the set of homotopy classes of orbifold loops with initial and end point point $x_0$
with the usual composition law by concatenation of loops,
as in the classical sense. The set $\pi^{\textup{orb}}_1(\mathcal{O},
\, x_0)$ is a group, and a change of base point
results in an isomorphic group. 
 \end{definition}

It is a theorem of Thurston, c.f. \cite[Th.\,2.5]{BMP}, that any connected smooth orbifold 
$\mathcal{O}$ has a unique up to equivalence orbifold covering $\widetilde{\mathcal{O}}$ which is universal and
whose orbifold fundamental group based at any regular point is trivial, so this definition of
universal covering extends to smooth orbifolds the classical definition for smooth manifolds.
Next we exhibit a construction, which is analogous 
to the construction of the universal cover of a generic orbifold,
of such covering for the case $\mathcal{O}=M/T$.

Let $x_0 \in M$, write $p_0=\pi(x_0) \in M/T$. 
By Proposition \ref{atlaslemma} the orbit space $M/T$ is a compact connected and 
smooth $(\op{dim}M-\op{dim}T)$\--dimensional orbifold. 
For each $p \in M/T$, let $\widetilde{M/T}_p$ denote the space of  homotopy 
classes of (orbifold) paths $\gamma \colon 
[0,\,1] \to M/T$ which start at $p_0$ and end at $p$. Let $\widetilde{M/T}$ denote the 
set--theoretic disjoint 
union of the spaces $\widetilde{M/T}_p$, where $p$ ranges over the orbit--space $M/T$.
Let $\psi$
be the set--theoretic mapping $\psi \colon \widetilde{M/T} \to M/T$ which sends $\widetilde{M/T}_p$ to $p$. As every smooth orbifold, the orbit space $M/T$ has a smooth orbifold atlas
$
\{(U_i, \, \widetilde{U}_i, \, \phi_i, \, \Gamma_i)\}_{i \in I}
$
in which the sets $\tilde{U}_i$ are simply connected. 
Pick $p_i \in U_i$, let $\mathcal{F}_i$ be the fiber of $\psi$ over $p_i$, and let $x_i$ be 
a point in $\mathcal{F}_i$, for each $i \in I$. 
Let $x \in \widetilde{U}_i, \, y \in \mathcal{F}_i$, 
choose $q_i \in \psi_i^{-1}(p_i)$ and choose
a path $\widetilde{\lambda}$ in $\widetilde{U}_i$ from $q_i$ to $x$, and let 
$\lambda:=\phi_i (\widetilde{\lambda})$ equipped with the
lift $\widetilde{\lambda}$ at each point, which is an (orbifold) path in $U_i$ such that $\lambda(0)=p_i$.
By definition, the element $y$ is a  homotopy class of (orbifold)
paths from $p_0$ to $p_i$. Let
$\gamma$ be a representative of $y$, and let $\alpha$ be the (orbifold) path obtained by concatenating
$\gamma$ with $\lambda$, where $\alpha$ travels along the points in $\gamma$ first. Since $\widetilde{U}_i$ is simply connected, the  homotopy class
of $\alpha$ does not depend on the choice of $\gamma$ and $\lambda$, and we define the
surjective map 
$
\psi_i(x, \, y):=[\alpha]
$
from the Cartesian product $\widetilde{U}_i \times \mathcal{F}_i$ to $\psi^{-1}(U_i)$.
Equip each fiber $\mathcal{F}_i$ with the discrete topology and the Cartesian product
$\widetilde{U}_i \times \mathcal{F}_i$ with the
product topology. Then $\psi$ induces a topology on $\psi^{-1}(U_i)$ whose connected
components are the images of sets of the form $\widetilde{U}_i \times \{y\}$ and we set
$(\widetilde{U}_i \times \{y\}, \, \psi_i|_{\widetilde{U_i} \times \{y\}})$, 
$y \in \mathcal{F}_i$, $i \in I$, as orbifold charts for an orbifold atlas for $\widetilde{M/T}$.
The pair $(\widetilde{M/T}, \, \psi)$ endowed with the
equivalence class of atlases of the orbifold atlas for $\widetilde{M/T}$ whose orbifold
 charts are $(\widetilde{U}_i \times \{y\}, \,\psi_i|_{\widetilde{U}_i \times \{y\}})$, $y \in
 \mathcal{F}_i$, $i \in I$, is a regular and universal orbifold covering of 
the orbifold $M/T$ called the \emph{universal cover of $M/T$ based at $p_0$}, and denoted simply
by $\widetilde{M/T}$, whose orbifold fundamental group
$\pi^{\textup{orb}}_1(\widetilde{M/T},\,p_0)$ is trivial. 

\subsubsection{Symplectic structure on $M/T$}

We prove that $M/T$ comes endowed with a symplectic structure.
A \emph{(smooth) differential form $\omega$} (resp. \emph{symplectic form}) \emph{on the smooth orbifold $\mathcal{O}$} 
is given by a collection $\{\widetilde{\omega}_i\}$ where $\widetilde{\omega}_i$ is a 
$\Gamma_i$\--invariant differential form (resp. symplectic form) on each
$\widetilde{U}_i$ and such that any two of them agree on overlaps. 
Here the notion of ``agreeing on overlaps'' can be made precise as follows:
for each $x \in \widetilde{U}_i$
and $y \in \widetilde{U}_j$ with $\phi_i(x)=\phi_j(y)$, there is a diffeomorphism
$\psi$ between a neighborhood of $x$ and a neighborhood of $y$
such that $\phi_j(\psi(z))=\phi_i(z)$ for all $z$ in such neighborhood
and $\psi^*\omega_j=\omega_i$.
A \emph{symplectic orbifold} is a a smooth orbifold equipped with a symplectic 
form.  

\begin{remark} \label{rmk}
Let $\omega$ be a differential form on an orbifold $\mathcal{O}$, and suppose that
$\omega$  is given by a collection $\{\widetilde{\omega}_i\}$ where $\widetilde{\omega}_i$ is a 
$\Gamma_i$\--invariant differential form (resp. symplectic form) on each
$\widetilde{U}_i$. Because the $\widetilde{\omega}_i$'s which define it
are $\Gamma_i$\--invariant and agree on overlaps,
$\omega$ is uniquely determined 
by its values on any orbifold atlas for $\mathcal{O}$ even if the atlas is not maximal.
\end{remark}

If $f \colon \mathcal{O}' \to \mathcal{O}$ is an orbifold diffeomorphism, the 
\emph{pull--back of a differential form $\omega$ on $\mathcal{O}$} is the 
unique differential form $\omega'$ on $\mathcal{O}'$ given by $\widetilde{f}^*\widetilde{\omega}_i$
on each chart $\widetilde{f}^{-1}(\widetilde{U}_i)$; we write $\omega':=f^*\omega$.
Analogously we define the pullback under a principal $T$\--orbibundle
$p \colon Y \to X/T$ as in Definition \ref{ic}, of a form $\omega$ on $X/T$,
where the maps $\widetilde{p}$ are given in Remark \ref{pullback}.

\begin{lemma} \label{3.7}
Let $(M, \,\sigma)$ be a compact, connected, symplectic  manifold equipped
with an effective (resp. free) symplectic action of a torus $T$ for which 
at least one, and hence every $T$\--orbit is 
a symplectic $\op{dim}T$\--dimensional submanifold of $(M, \, \sigma)$. 
Then there exists a unique $2$\--form $\nu$ on the orbit space $M/T$ such that
$\pi^* \nu|_{\Omega_x}=\sigma|_{\Omega_x}$ for every $x \in M$,
where $\{\Omega_x\}_{x \in M}$ is the distribution on $M$ of
symplectic orthogonal complements to the tangent spaces to the $T$\--orbits,
and $\pi \colon M \to M/T$ is the projection map, c.f. Proposition \ref{easynow}.
Moreover, the form $\nu$ is symplectic, and so the pair $(M/T,\, \nu)$
is a compact, connected symplectic orbifold (resp. manifold).
\end{lemma}

\begin{proof}
By Corollary \ref{tubetheorem2}, for each $x \in M$ there is a $T$\--invariant neighborhood
$U_x$ of the $T$\--orbit $T\cdot x$ and a $T$\--equivariant symplectomorphism
$\Phi_x$ from $U_x$ onto $T\times_{T_x}D_x$ with the symplectic form
$\widetilde{\sigma^T \oplus \sigma^{\C^m}}$. 
As in Proposition \ref{atlaslemma}, the collection of charts 
$\{(\pi(U_x), \, D_x,\,\phi_x,\,T_x)\}$ given by expression (in the general definition of 
smooth orbifold $\widetilde{U}_i=D_x, \, U_i=\pi(U_x)$)
 (\ref{mtatlas}) is an orbifold atlas
for the orbit space $M/T$, which is equivalent to the one given in Proposition \ref{atlaslemma} and hence
defines the orbifold structure of $M/T$ given in Definition \ref{maximal}. Because by Remark \ref{rmk}
it suffices to define a smooth differential form on the charts of any atlas of our choice,
the collection $\{\nu_x\}_{x \in M}$ given by $\nu_x:=\sigma^{\C^m}$ defines
a unique smooth differential $2$\--form on the orbit space $M/T$. Because $\sigma^{\C^m}$
is moreover symplectic, each $\nu_x$ is a symplectic form, and hence so is $\nu$ on $M/T$.
\end{proof}

We say that two symplectic orbifolds
$(\mathcal{O}_1,\,\nu_1), \, (\mathcal{O}_2,\,\nu_2)$ are \emph{symplectomorphic}
if there is an orbifold diffeomorphism $f \colon \mathcal{O}_1 \to \mathcal{O}_2$
with $f^*\nu_2=\nu_1$. $f$ is called
an \emph{orbifold symplectomorphism}.

\begin{lemma} \label{orbifoldsymplectomorphic}
Let $(M, \,\sigma), \, (M', \, \sigma')$ be compact connected symplectic
manifolds equipped with an effective (resp. free) symplectic action of a torus $T$ for which at least
one, and hence every $T$\--orbit is a $\op{dim}T$\--dimensional 
symplectic submanifold of $(M, \, \sigma)$ and $(M', \, \sigma')$, respectively. 
Suppose additionally that $(M, \, \sigma)$ is $T$\--equivariantly symplectomorphic to $(M',\, \sigma')$.
Then the symplectic 
orbit spaces $(M/T, \, \nu)$ and $(M'/T, \, \nu')$ are symplectomorphic. 
\end{lemma}

\begin{proof}
By Lemma \ref{finitestabilizers00}, every $T$\--orbit is a $\op{dim}T$\--dimensional
symplectic submanifold of $(M', \, \sigma')$.
Indeed, assume that $\Phi$ is a $T$\--equivariant symplectomorphism from the symplectic manifold
$(M, \, \sigma)$ to the symplectic manifold $(M', \, \sigma')$. 
Because the mapping $\Phi$ is $T$\--equivariant, it descends to an orbifold
diffeomorphism $\widehat{\Phi}$ from the orbit space $M/T$ onto the orbit space
$M'/T$. Because of the definition of the connections $\Omega$ and of $\Omega'$ for
the orbibundle projections $\pi \colon M \to M'/T$ and $\pi' \colon M' \to M'/T$
respectively, as the connections given by the symplectic orthogonal complements to the
tangent spaces to the $T$\--orbits, we have that $\Phi^*\Omega'=\Omega$. 

Since by Lemma \ref{3.7} the symplectic forms $\nu$ on $M/T$ and $\nu'$ 
on $M'/T$ are respectively induced by the projections $\pi$ and $\pi'$, and both
$\Phi^*\Omega'=\Omega$ and $\Phi^* \sigma'=\sigma$, we have that $\widehat{\Phi}^* \nu'=\nu$
and therefore $\widehat{\Phi}$ is an orbifold symplectomorphism between
$(M/T, \, \nu)$ and $(M'/T, \, \nu')$.
\end{proof}

We use the notation $\int_{\mathcal{O}} \omega$ 
for the integral of the differential form $\omega$ on the orbifold
$\mathcal{O}$. If $(\mathcal{O},\, \omega)$
is a symplectic orbifold, the integral $\int_{\mathcal{O}} \omega$
is known as the \emph{total symplectic area of $(\mathcal{O}, \,
\omega)$}. 

\begin{remark} \label{moser}
Under the assumptions of Lemma \ref{orbifoldsymplectomorphic},
if $\op{dim}T=\op{dim}M-2$, then
the total symplectic area of the symplectic orbit space $(M/T, \, \nu)$
equals the total symplectic area of $(M'/T, \, \nu')$. 
Moreover, recall that by the orbifold Moser's theorem \cite[Th.\,3.3]{MW}, 
if $(M, \, \sigma)$, $(M', \, \sigma')$ are any symplectic manifolds equipped
with an action of a torus $T$ of dimension $\op{dim}T=\op{dim}M-2$, for which
at least one and hence every $T$\--orbit is a $\op{dim}T$\--dimensional symplectic
submanifold of $(M, \, \sigma)$, and such that
the symplectic area of the orbit space
$(M/T, \, \nu)$ equals the symplectic area of the orbit space
$(M'/T, \, \nu')$, and the Fuchsian signature of the orbisurfaces $M/T$ and $M'/T$ are equal,
then $(M/T, \, \nu)$ is $T$\--equivariantly (orbifold) symplectomorphic to $(M'/T,\, \nu')$; 
we emphasize that this is the case because the orbit spaces $M/T$ and $M'/T$
are $2$\--dimensional, so their genus determines their diffeomorphism type. 
This result is used in the proof of Theorem \ref{FTT}
and Theorem \ref{TT}.
 \end{remark}

\subsubsection{Model of $(M, \, \sigma)$: definition}

In the following definition we define a $T$\--equivariant symplectic model for
$(M, \, \sigma)$: for each regular point $p_0 \in M/T$ we construct a smooth
manifold $M_{\textup{model},\, p_0}$, a $T$\--invariant symplectic form
$\sigma_{\textup{model}}$ on  $M_{\textup{model},\, p_0}$, and an effective
symplectic action of the torus $T$ on  $M_{\textup{model},\, p_0}$. 
The ingredients for the construction of
such model are Proposition \ref{easynow} and Lemma \ref{3.7}.

If $p \colon \mathcal{O}' \to \mathcal{O}$ is an orbifold covering map, 
$p$ is a local diffeomorphism in the sense of Remark \ref{localdiffeo},
and the 
\emph{pull--back of a differential form $\omega$ on $\mathcal{O}$} is the 
unique differential form $\omega'$ on $\mathcal{O}'$ given by $\widetilde{p}^*\widetilde{\omega}_i$
on each chart $\widetilde{p}^{-1}(\widetilde{U}_i)$; we write $\omega':=p^*\omega$.

\begin{remark} \label{modelremark}
This remark justifies that item ii) in Definition \ref{themodel} below is correctly
defined. The reader may proceed with Definition \ref{themodel}
and return to this remark as therein indicated.

The pull--back of the symplectic form $\nu$ on $M/T$, given by Lemma \ref{3.7}, to 
the universal cover $\widetilde{M/T}$, by means of the smooth covering map 
$\psi \colon \widetilde{M/T} \to M/T$, is 
a ${\pi^{\textup{orb}}_1(M/T, \,p_0)}$\--invariant symplectic form on the orbifold
universal cover $\widetilde{M/T}$. 
The symplectic form on the torus $T$ is translation invariant and therefore ${\pi^{\textup{orb}}_1(M/T, \,p_0)}$\--invariant. 
The direct sum of the 
symplectic form on $\widetilde{M/T}$ and the symplectic form on $T$ is a ${\pi^{\textup{orb}}_1(M/T, \,p_0)}$\--invariant 
(and $T$\--invariant) symplectic form on the Cartesian product
$\widetilde{M/T}\times T$, and therefore there exists a 
unique symplectic form on the associated bundle 
$\widetilde{M/T}\times_{{\pi^{\textup{orb}}_1(M/T, \,p_0)}} T$ of which the pull--back by the 
covering map $\widetilde{M/T} \times T \to \widetilde{M/T} \times_{{\pi^{\textup{orb}}_1(M/T, \,p_0)}} T$ is equal to the 
given symplectic form on $\widetilde{M/T} \times T$. 
\end{remark}

\begin{definition} \label{themodel}
Let $(M, \,\sigma)$ be a compact, connected, symplectic manifold equipped
with an effective symplectic action of a torus $T$ for which at least
one, and hence every $T$\--orbit is  
a $\op{dim}T$\--dimensional
symplectic submanifold of $(M, \, \sigma)$. We define the \emph{$T$\--equivariant symplectic model 
$(M_{\textup{model},\, p_0},\, \sigma_{\textup{model}})$
of $(M, \, \sigma)$ based at a regular point $p_0 \in M/T$} as follows.
\begin{itemize}
\item[i)]
The space $M_{\textup{model},\, p_0}$ is the associated bundle
$$
M_{\textup{model}, \, p_0}:=\widetilde{M/T} \times_{\pi_1^{\textup{orb}}(M/T, \, p_0)} T,
$$
where the space $\widetilde{M/T}$ denotes the orbifold universal
cover of the orbifold $M/T$ 
based at a regular point $p_0 \in M/T$, and the orbifold fundamental group 
$\pi^{\textup{orb}}_1(M/T,\,p_0)$ acts on the product $\widetilde{M/T} 
\times T$ by the diagonal action $x \, (y,\,t)=(x \star y^{-1},\,\mu(x) \cdot t)$, where 
$\star \colon \pi_1^{\textup{orb}}(M/T, \, p_0) \times \widetilde{M/T} \to \widetilde{M/T}$ denotes the natural
action of $\pi_1^{\textup{orb}}(M/T, \, p_0)$ on $\widetilde{M/T}$, 
and $\mu \colon \pi_1^{\textup{orb}}(M/T, \, p_0)
\to T$ denotes the monodromy homomorphism of the flat connection
$\Omega:=\{\Omega_x\}_{x \in M}$ given by the symplectic orthogonal complements to the 
tangent spaces to the $T$\--orbits (c.f. Proposition \ref{easynow}). 
\item[ii)]
The symplectic form  $\sigma_{\textup{model}}$ is induced on the quotient by the product
symplectic form on the Cartesian product $\widetilde{M/T} \times T$. 
The symplectic form on $\widetilde{M/T}$ is defined as the pullback by the
orbifold universal covering map $\widetilde{M/T} \to M/T$ of the unique 
$2$\--form $\nu$ on $M/T$ such that
$\pi^* \nu|_{\Omega_x}=\sigma|_{\Omega_x}$ for every $x \in M$ (c.f. Lemma \ref{3.7}).
The symplectic form on the torus
$T$ is the unique $T$\--invariant symplectic form determined by the non--degenerate
antisymmetric bilinear form $\sigma^{\mathfrak{t}}$ such that
$\sigma_x(X_M(x), \,Y_M(x))=\sigma^{\mathfrak{t}}(X,Y)$, 
for every $X, \,Y \in \mathfrak{t}$, and every $x \in M$ (c.f. Lemma \ref{sigmat}).
See Remark \ref{modelremark}.
\item[iii)]
The action of $T$ on the space $M_{\textup{model},\, p_0}$  is the action of $T$ by translations
which descends from the action of $T$ by translations on the right factor of the product
$\widetilde{M/T} \times T$.
\end{itemize}
 \end{definition}

\normalfont

\begin{example} \label{exx}
\normalfont
When $\op{dim}M-\op{dim}T=2$ the orbit space $M/T$ is $2$\--dimensional and
it is an exercise to describe locally the monodromy homomorphism
$\mu$ which appears in part i) of Definition \ref{themodel}. A small $T$\--invariant open subset
of our symplectic manifold looks like $T \times_{T_x}D$, where $T$ is the standard 
$2$\--dimensional torus $(\R/\Z)^2$, and $D$ is a standard
$2$\--dimensional disk centered at the origin in the complex plane $\C$. Here the
quotient $M/T$ is the orbisurface $D/T_x$. (Recall that we know that
$T_x$ is a finite cyclic group, c.f. Lemma \ref{finitestabilizers}). Suppose that $T_x$ has order $n$. The
monodromy homomorphism $\mu$ in Definition \ref{themodel} part i)
is a map from 
$\pi^{\textup{orb}}_1:=\pi^{\textup{orb}}_1(D/T_x,\,p_0)$ into $T$. Note that 
$\pi^{\textup{orb}}_1=\langle \gamma \rangle$ with $\gamma$ of order $n$. 
If $t \in T$,
there exists a homomorphism $f \colon \pi^{\textup{orb}}_1 \to T$ such that $f(\gamma)=t$
of and only if $t^n=1$. 

If we identify $T$ with $(\R/\Z)^2$, we can write $t=(t_1,\, t_2)$,
there exists a homomorphism $f \colon \pi^{\textup{orb}}_1 \to T$ such that $f(\gamma)=(t_1,\, t_2)$
of and only if $n$ divides the order of $(t_1, \, t_2)$, which means that $t_1,\,t_2$
must be rational numbers such that $n$ divides the smallest integer $m$ such
that $m\,t_i \in \Z$, $i=1,\,2$. If for example $n=2$, this condition says
that the smallest integer $m$ such that $m\, t_i \in \Z$ must
be an even number. In other words, not every element in $T$ can be
achieved by the monodromy homomorphism. In fact, all the achievable elements
are of finite order, but as we see from this example, more restrictions must take place.
 \end{example}

\subsubsection{Model of $(M,\,\sigma)$: proof}

We prove that the associated bundle in
Definition \ref{themodel}, which we called ``the model of $M$'', is 
$T$\--equivariantly symplectomorphic to $(M, \, \sigma)$.
The main ingredient of the proof is
the existence of the flat connection for $\pi \colon M \to M/T$
in Proposition \ref{easynow}.

We start with the observation that the universal
cover $\widetilde{M/T}$ is a smooth manifold and the orbit space $M/T$ is a good orbifold.
Recall that an orbifold $\mathcal{O}$ is said to be \emph{good} (resp. \emph{very good}) if it is obtained as the quotient of manifold by a discrete (resp. finite) group of diffeomorphisms of such manifold.

\begin{lemma} \label{goodorbifoldlemma}
Let $(M, \,\sigma)$ be a compact, connected, symplectic manifold equipped
with an effective symplectic action of a torus $T$ for which at least one, and hence
every $T$\--orbit is a $\op{dim}T$\--dimensional
symplectic submanifold of $(M, \, \sigma)$. Then the orbifold universal
cover $\widetilde{M/T}$ is a smooth manifold and the orbit space $M/T$ is a 
good orbifold. Moreover, if $\op{dim}T=\op{dim}M-2$, 
the orbit space $M/T$ is a very good orbifold.
\end{lemma}

\begin{proof}
Let $\psi \colon \widetilde{M/T} \to M/T$ be the universal cover of orbit space
$M/T$ based at the regular point $p_0=\pi(x_0) \in M/T$.
Recall from Proposition \ref{easynow} the connection $\Omega$ for the principal $T$\--orbibundle
$\pi \colon M \to M/T$ projection mapping, whose elements are the
symplectic orthogonals to the tangent spaces to the $T$\--orbits of $(M, \, \sigma)$.
By Corollary \ref{obvious}, if the mapping $i_x \colon \mathcal{I}_x \to M$ is the inclusion
mapping of the integral manifold $\mathcal{I}_x$ through $x$ to the distribution
$\Omega$, then the composite
map $\pi \circ i_x \colon \mathcal{I}_x \to M/T$ is an orbifold covering mapping for each $x \in M$,
in which the total space $\mathcal{I}_x$ is a smooth manifold. Because the
covering map $\psi \colon \widetilde{M/T} \to M/T$ is universal, there
exists an orbifold covering $r \colon \widetilde{M/T} \to \mathcal{I}_x$
such that $\pi \circ i_x \circ r=\psi$, and in particular the orbifold universal cover
$\widetilde{M/T}$ is a smooth manifold, since an orbifold covering
of a smooth manifold must be a smooth manifold itself. 
Since the orbit $M/T$ is obtained as a quotient of $\widetilde{M/T}$ by the discrete group
$\pi^{\textup{orb}}_1(M/T, \, p_0)$, by definition $(M/T, \, \mathcal{A})$ is a good orbifold.
Now, it is well--known \cite[Sec.\,2.1.2]{BMP} that in dimension $2$, an orbifold is good if and only if it is very good.
\end{proof}

Let $\alpha  \colon [0, \,1] \to \mathcal{\widehat{O}}$
be an orbifold path in $\widehat{\mathcal{O}}$.
The \emph{projection
of the path $\alpha$ under a covering mapping
$p \colon \widehat{\mathcal{O}} \to \mathcal{O}$}, 
which we write as $p(\alpha)$, is a path in the orbifold 
$\mathcal{O}$ whose underlying map is
$p \circ \alpha$, and such that for each $t$ for which
$(p \circ \alpha)(t)$ is a singular point, if $\alpha(t)$
is regular then there exists a neighborhood
$V$ of $\alpha(t)$ such that $p|_V$ is a chart at $(p \circ \alpha)(t)$,
and this can be used to define the local lift of
$p \circ \alpha$ around $t$; if otherwise $\alpha(t)$
is singular, then by definition of $\alpha$ there is a chart
$\phi \colon \widetilde{V} \to V$ at
$\alpha(t)$ and a local lift
to $\widetilde{V}$ of $\alpha$ around $t$, and 
by choosing $\widetilde{V}$ to be small
enough, $\phi \circ p$ is a chart of $(p \circ \alpha)(t)$
giving the local lift of $p \circ \alpha$ around $t$. We say
that \emph{$\alpha$ is a lift of a path $\beta$ in $\mathcal{O}$}
if $p(\alpha)=\beta$. 
The following generalizes the classical manifold property. Recall Definition \ref{ic}.

\begin{lemma} \label{liftlemma}
Let $T$ be a torus, let $(X, \, \omega)$ be a compact, connected symplectic manifold 
endowed with an effective symplectic action of $T$ for which at least one, and hence
every $T$\--orbit is a $\op{dim}T$\--dimensional symplectic submanifold of $(X, \, \omega)$.
Equip $X/T$ with the orbifold structure in 
Definition \ref{maximal}.
Let $p_0 \in X/T$ be a regular point.
Let $H$ be a connection for a smooth principal $T$\--orbibundle $p \colon X \to X/T$.  
Then for any loop $\gamma \colon [0,\,1]\to X/T$ in the orbifold $X/T$
such that $\gamma(0)=p_0$ there exists
a unique horizontal lift  $\lambda_{\gamma} \colon [0,\,1]\to X$ with respect to the 
connection $H$ for $p \colon X \to X/T$, such that 
$\lambda_{\gamma}(0)=x_0$,
where by horizontal we mean that $\op{d}\lambda_{\gamma}(t)/\op{d}t\in 
H_{\lambda_{\gamma(t)}}$ for every $t \in [0,\, 1]$.
\end{lemma}

The following is the main result of Section 2.

\begin{theorem} \label{tt}
Let $(M, \,\sigma)$ be a compact, connected, symplectic manifold equipped
with an effective symplectic action of a torus $T$, for which at least one, and hence every $T$\--orbit is 
a $\op{dim}T$\--dimensional symplectic submanifold of $(M,\, \sigma)$. Then
$(M,\,\sigma)$ is $T$\--equivariantly
symplectomorphic to its $T$\--equivariant symplectic model based 
at a regular point $p_0 \in M/T$, c.f. Definition \ref{themodel}, for every regular point
$p_0 \in M/T$.
\end{theorem}

\begin{proof} 
Let $\psi \colon \widetilde{M/T} \to M/T$ be the universal cover of orbit space
$M/T$ based at the point $p_0=\pi(x_0) \in M/T$.	
Recall from Proposition \ref{easynow} the connection $\Omega$ for the $T$\--orbibundle
$\pi \colon M \to M/T$ projection mapping given by the
symplectic orthogonals to the tangent spaces to the $T$\--orbits of $(M, \, \sigma)$.
By Lemma \ref{goodorbifoldlemma}, $\widetilde{M/T}$ is a smooth manifold.
The mapping $\pi^{\textup{orb}}_1(M/T, \, p_0) \times \widetilde{M/T} \to 
\widetilde{M/T}, \, \, \, ([\lambda], \,[\gamma]) \mapsto [\gamma \, \lambda]$,
is a smooth action of the orbifold fundamental group
$\pi^{\textup{orb}}_1(M/T, \,p_0)$ on the orbifold universal cover
$\widetilde{M/T}$, 
which is transitive on each fiber $\widetilde{M/T}_p$ of $\psi \colon \widetilde{M/T} \to M/T$. 
By Lemma \ref{liftlemma}, for any loop $\gamma \colon [0,\,1]\to M/T$ in the orbifold $M/T$
such that $\gamma(0)=p_0$, denote 
by $\lambda_{\gamma} \colon [0,\,1]\to M$ its unique horizontal lift with respect to the 
connection $\Omega$ for $\pi \colon M \to M/T$ such that 
$\lambda_{\gamma}(0)=x_0$,
where by horizontal we mean that $\op{d}\lambda_{\gamma}(t)/\op{d}t\in 
\Omega_{\lambda_{\gamma(t)}}$ for every $t \in [0,\, 1]$.
Proposition \ref{easynow} says that $\Omega$ is an orbifold flat connection, 
which means that $\lambda_{\gamma}(1)=\lambda_{\delta}(1)$ if $\delta$ is 
homotopy equivalent 
to $\gamma$ in the space of all orbifold paths in the orbit space
$M/T$ which start at $p_0$ and end at the given end point 
$p=\gamma(1)$. The gives existence of a unique group homomorphism $\mu \colon \pi^{\textup{orb}}_1(M/T, \, p_0) \to 
T$ such that $\lambda_{\gamma}(1)=\mu([\gamma]) \cdot x_0$. The homomorphism
$\mu$ does not depend
on the choice of the base point $x_0 \in M$.
For any homotopy class 
$[\gamma]\in \widetilde{M/T}$ and $t\in T$, define the element \begin{eqnarray} 
\label{kk}
\Phi([\gamma],\, t) := t\cdot\lambda_{\gamma}(1) \in M.
\end{eqnarray} 
The assignment $([\gamma],\,t) \mapsto \Phi([\gamma],\,t)$ defines a 
smooth covering mapping $\Phi \colon \widetilde{M/T} \times T\to M$ between smooth manifolds. Let
$[\delta]\in\pi^{\textup{orb}}_1(M/T, \, p_0)$ act on the product  $\widetilde{M/T} \times T$ by sending 
the pair $([\gamma],\, t)$ to the pair
$([\gamma \, \delta^{-1}],\,  \mu([\delta])\, t)$.
The mapping $\Phi$ induces a diffeomorphism $\phi$ from 
the associated bundle $\widetilde{M/T}\times_{\pi^{\textup{orb}}_1(M/T, \, p_0)} T$ onto $M$.
Indeed, $\phi$ is onto because $\Phi$ is onto, and
$\Phi([\gamma],\, t) = \Phi([\gamma'],\, t')$ if and only if 
$t\cdot\lambda_{\gamma}(1) = t'\cdot\lambda_{\gamma'}(1)$, if and only if 
$\lambda_{\gamma}(1) = (t^{-1} \,t')\cdot\lambda_{\gamma'}(1)$, if and only if $t^{-1} 
\,t' = \mu([\delta])$, in which $\delta$ is equal to the loop starting and ending at 
$p_0$, which is obtained by first doing the path $\gamma$ and then going back by means of 
the path $\gamma'^{-1}$. Since $\Phi$ is a smooth covering map, the mapping $\phi$ is a local
diffeomorphism, and we have just proved that it is bijective, so $\phi$ must be
a diffeomorphism. By definition, $\phi$ intertwines the action of $T$ by translations on 
the right factor of the associated bundle
$\widetilde{M/T} \times_{\pi^{\textup{orb}}_1(M/T, \, p_0)} T$ with the action of $T$ on $M$.
Recall that the symplectic form on
$\widetilde{M/T} \times_{\pi^{\textup{orb}}_1(M/T, \, p_0)} T$ is
the unique symplectic form of which the pullback
$\widetilde{M/T} \times T \to \widetilde{M/T} \times_{\pi^{\textup{orb}}_1(M/T, \, p_0)} T$ is equal to the
product form on $\widetilde{M/T} \times T$, where the symplectic form
on $\widetilde{M/T}$ is given by Definition \ref{themodel} part ii). It follows
from the definition of the symplectic form 
on $\widetilde{M/T} \times_{\pi^{\textup{orb}}_1(M/T, \, p_0)} T$
and Corollary \ref{tubetheorem2} that the $T$\--equivariant 
diffeomorphism $\phi$ from $\widetilde{M/T} \times_{{\pi^{\textup{orb}}_1(M/T, \,p_0)}} T$ onto $M$ pulls back the 
symplectic form $\sigma$ on $M$ to the just obtained symplectic form on $\widetilde{M/T} 
\times_{{\pi^{\textup{orb}}_1(M/T, \,p_0)}} T$.
\end{proof}

\normalfont
\begin{remark} \label{e2}
\normalfont
We assume the terminology of Definition \ref{themodel} and Theorem \ref{tt}.
In this remark we describe a covering isomorphic to the orbifold covering
$\pi \circ i_x$ of Lemma \ref{obvious}, which will be of use in Theorem \ref{altermodel}.
Additonally, the construction of this new covering makes transparent the relation between the
universal cover $\widetilde{M/T} \to M/T$ and the covering $\mathcal{I}_x \to M/T$.
Indeed, recall the distribution $\Omega$ of the symplectic orthogonal complements 
of the tangent spaces of the $T$\--orbits given in Proposition \ref{easynow},
and the principal $T$\--orbibundle $\pi : M \to M/T$, for which $\Omega$ 
is a $T$\---invariant flat connection. As in the proof of Theorem \ref{tt},
if $x\in M$ then each smooth curve $\delta$ in $M/T$ which starts at $\pi(x)$ has a unique 
horizontal lift $\gamma$ which starts at $x$. The endpoints of such lifts form 
the injectively immersed manifold $\mathcal{I}_x$, c.f. Remark \ref{re1}. 
As in the proof of Theorem \ref{tt} ,
if we keep the endpoints of $\delta$ fixed, then the endpoint of 
$\gamma$ only depends on the homotopy class of $\delta$. 
As explained prior to Definition \ref{themodel},
the homotopy classes of $\delta$'s in $M/T$ with fixed endpoints are 
by definition the elements of the universal covering space 
$\widetilde{M/T}$ of $M/T$, which is a principal $\pi^{\textup{orb}}_1(M/T, \, \pi(x))$\--bundle 
over $M/T$. The corresponding endpoints of the $\gamma$'s exhibit the
integral manifold
$\mathcal{I}_x$ as the image of an immersion from $\widetilde{M/T}$ into $M$,
but this immersion is not necessarily injective. 
Recall that the monodromy homomorphism $\mu$ of $\Omega$
tells what the endpoint of $\gamma$ is when $\delta$ 
is a loop.  By replacing $\mu$ by the induced injective homomorphism 
$\mu'$ from $\pi^{\textup{orb}}_1(M/T, \, \pi(x))/\textup{ker}\, \mu$ to $T$, we get an injective immersion. 
This procedure is equivalent to replacing the universal covering $\widetilde{M/T}$ 
by the covering $\widetilde{M/T}/\textup{ker}\, \mu$ with fiber 
$\pi^{\textup{orb}}_1(M/T, \, \pi(x))/\textup{ker}\, \mu$.
So $\widetilde{M/T}/\textup{ker}\, \mu$ is injectively immersed in $M$ by
a map whose image is $\mathcal{I}_x$. It follows that the coverings
$\pi \circ i_x \colon \mathcal{I}_x \to M/T$ and $\widetilde{M/T}/\textup{ker}\, \mu \to M/T$
are isomorphic coverings of $M/T$, of which the total space is a smooth manifold. 
We use this new covering  $\widetilde{M/T}/\textup{ker}\, \mu \to M/T$ in
order to construct an alternative model of $(M, \, \sigma)$ to the one given
in Theorem \ref{tt}, c.f. Theorem \ref{altermodel}. Note: in the statement of
Theorem \ref{altermodel}, $N:=\widetilde{M/T}/\textup{ker}\, \mu$. 
\end{remark}

\begin{example}
\normalfont
Assume the terminology of Definition \ref{themodel}. Assume moreover that $\op{dim}M-
\op{dim}T=2$ and that
the action of $T$ on $M$ is free, so the quotient space $M/T$ is a compact, connected,
smooth surface and it is classified by its genus.
If the genus is zero, then the orbit space $M/T$ is a sphere, 
which is simply connected, and the orbibundle $\pi :\mathcal{I}_x\to M/T$ is a diffeomorphism, c.f.
Lemma \ref{obvious}. In such case $M$ is the Cartesian product of a sphere with a torus.
If $M/T$ has genus $1$, then the orbit space $M/T$ is a two--dimensional torus, with
fundamental group at any point isomorphic to the free abelian group on two generators $t_1, \, t_2$,
and the monodromy homomorphism is determined by the 
images $t_1$ and $t_2$ of the generators $(1,\,  0)$ and $(0,\, 1)$ of 
$\Z^2$. $\mathcal{I}_x$ is compact if and only if $\mathcal{I}_x$ is a closed subset of $M$ 
if and only if $t_1$ and $t_2$ generate a closed subgroup of $T$ 
if and only if $t_1$ and $t_2$ generate a finite subgroup of $T$. 
$\mathcal{I}_x$ is dense in $M$ if and only if the subgroup 
of $T$ generated by $t_1$ and $t_2$ is dense in $T$. 
If the genus of $M/T$ is strictly positive, then
$\mathcal{I}_x$ is compact if and only if the monodromy elements form a finite subgroup 
of $T$, which is a very 
particular situation (since the $2g$ generators of the monodromy subgroup of $T$
can be chosen arbitrarily, this is very rare: even one 
element of $T$ usually generates a dense subgroup of $T$). 
\end{example}

\subsection{Model up to $T$\--equivariant diffeomorphisms}

\subsubsection{Generalization of Kahn's theorem}

In \cite[Cor. 1.4]{Kahn} P. Kahn states that if a compact connected $4$\--manifold 
$M$ admits a free action of a $2$\--torus $T$ such that the $T$\--orbits are 
$2$\--dimensional symplectic
submanifolds, then $M$ splits as a product $M/T \times T$, the following being the 
statement in \cite{Kahn}.

\begin{theorem}[P. Kahn, \cite{Kahn}, Cor. 1.4] \label{Kahntheorem}
Let $(M, \, \sigma)$ be a compact, connected symplectic $4$\--dimensional manifold. Suppose that
$M$ admits a free action of a $2$\--dimensional torus $T$ for which the $T$\--orbits
are $2$\--dimensional symplectic submanifolds of $(M, \, \sigma)$. Then there exists a $T$\--equivariant
diffeomorphism between $M$ and the cartesian product $M/T \times T$, where
$T$ is acting by translations on the right factor of $M/T \times T$.
\end{theorem} 

Next we generalize this result of Kahn's to the case when the torus and the manifold are of arbitrary dimension, c.f Theorem \ref{mainlemma}. We start with Corollary \ref{KK}, the proof 
of which relies on the forthcoming Theorem \ref{mainlemma}, so we
suggest to return to its proof at the end of the section.

\begin{cor} \label{KK}
Let $(M, \,\sigma)$ be a compact connected symplectic $2n$\--dimensional manifold equipped
with a free symplectic action of a $(2n-2)$\--dimensional torus $T$ such that at least
one, and hence every $T$\--orbit is a $(2n-2)$\--dimensional 
symplectic submanifold of $(M, \, \sigma)$. Then  $M$ is $T$\--equivariantly diffeomorphic to
the Cartesian product $M/T \times T$, where $M/T \times T$ 
is equipped with the action of $T$ on the right factor of $M/T \times T$ 
by translations.
\end{cor}

\begin{proof}
If $\textup{dim}(M/T)=2$, by the first part of Theorem \ref{mainlemma}, the
quotient space $M/T$ is a compact, connected, smooth, orientable surface.
Therefore, by the classification theorem for surfaces, $\op{H}_1(M/T, \,\Z)$
is isomorphic to $\Z^{2g}$, where $g$
is the topological genus of $M/T$. Therefore 
$\op{H}_1(M/T, \,\Z)$ does not have torsion elements, and
by Theorem \ref{mainlemma}, $M$ is $T$\--equivariantly diffeomorphic to the 
cartesian product $M/T \times T$. 
\end{proof}

\subsubsection{Smooth equivariant splittings}

We give a characterization of
the existence of $T$\--equivariant
splittings of $M$ as a Cartesian product $M/T \times T$ , up to $T$\--equivariant diffeomorphisms c.f. Theorem \ref{mainlemma}.

Recall that if $X$ is a smooth manifold,
$\op{H}_1(X,\,\Z)_{\op{T}}$ denotes the subgroup of torsion elements of 
$\op{H}_1(X,\,\Z)$
i.e. the set of $[\gamma] \in \op{H}_1(X,\,\Z)$ such that there exists
a strictly positive integer
$k$ satisfying $k \,[\gamma]=0$. 
$\op{H}_1(X,\,\Z)_{\op{T}}$ is a subgroup.
Let $\mu^{\op{T}}_{\textup{h}}$ denote
the homomorphism given by restricting the homomorphism
$\mu_{\textup{h}}$ to $\op{H}_1(M/T,\,\Z)_{\op{T}}$. 
Recall that by Proposition \ref{atlaslemma}, $M$ is a compact, connected
$(\op{dim}M-\op{dim}T)$\--dimensional orbifold,
and that if the action of torus $T$ on $M$ 
is free, the local groups are trivial, and hence 
the orbit space $M/T$ is a $(\op{dim}M-\op{dim}T)$\--dimensional 
smooth manifold, c.f. Remark \ref{manifoldcase}.

\begin{theorem} \label{mainlemma}
Let $(M, \,\sigma)$ be a compact, connected, symplectic manifold equipped
with a free symplectic action of a torus $T$ for which at least one, and hence
every of its $T$\--orbits is 
a $\op{dim}T$\--dimensional symplectic submanifold of $(M, \, \sigma)$, 
and let $\mu_{\textup{h}} \colon \op{H}_1(M/T,\,\Z) \to T$ 
be the homomorphism induced on homology via the Hurewicz map 
by the monodromy homomorphism $\mu$ from 
$\pi_1(M/T, \,p_0)$ into $T$  with
respect to the connection on $M$ given by the symplectic orthogonal complements
to the tangent spaces to the $T$\--orbits, for the $T$\--bundle projection map
$\pi \colon M \to M/T$ onto the orbit space $M/T$ (c.f. expression (\ref{muh})). 
Then $M$ is $T$\--equivariantly diffeomorphic to the 
Cartesian product $M/T \times T$ equipped with the action of $T$ by translations
on the right factor of $M/T \times T$,
if and only if the torsion part $\mu_{\textup{h}}^{\op{T}}$
of the homomorphism $\mu_{\textup{h}}$ is trivial, i.e. $\mu_{\textup{h}}$
satisfies that $\mu_{\textup{h}}([\gamma])=1$ for every $[\gamma]\in 
\op{H}_1(M/T,\,\Z)$ of finite order.  
\end{theorem}

\begin{proof} By means of the same argument that
we used in the proof of Theorem \ref{tt}, using any $T$\--invariant flat 
connection $\mathcal{D}$ for $\pi$
instead of $\Omega$, we may define a mapping 
\begin{eqnarray} \label{map}
\phi^{\mathcal{D}} \colon \widetilde{M/T} \times_{\pi_1(M/T, \, p_0)} T \to M.
\end{eqnarray}
The mapping $\phi^{\mathcal{D}}$, which is induced by (\ref{kk}), is a $T$\--equivariant diffeomorphism 
between the smooth manifolds $\widetilde{M/T} \times_{\pi_1(M/T, \, p_0)} T$ and $M$, where 
$\widetilde{M/T}$ is the 
universal cover of the manifold $M/T$ based at $p_0$. 
Since $\widetilde{M/T}$ is a regular covering c.f. Remark \ref{natural},
$\widetilde{M/T}/\pi_1(M/T, \, p_0)=M/T$ 
and the monodromy homomorphism mapping
$\mu^{\mathcal{D}} \colon \pi_1(M/T, \, p_0) 
\to T$ of $\mathcal{D}$ is trivial if and only if $\widetilde{M/T} \times_{\pi_1(M/T, \, p_0)} T=M/T \times 
T$. Hence the smooth manifold
$M$ is $T$\--equivariantly diffeomorphic to the Cartesian product $M/T \times T$ if and only if 
there exists a $T$\--invariant flat connection $\mathcal{D}$ for the principal 
$T$\--bundle $\pi \colon M \to M/T$ for which its monodromy homomorphism $\mu^{\mathcal{D}} 
\colon \pi_1(M/T, \, p_0) \to T$ is trivial.

If $\mu^{\op{T}}_{\textup{h}}=0$, 
there exists a homomorphism $\widehat{\mu_{\textup{h}}}
\colon \op{H}_1(M/T,\,\Z) \to \mathfrak{t}$ such that
$
\textup{exp} \circ \widehat{\mu_{\textup{h}}}=
\mu_{\textup{h}}$, which by
viewing $\widehat{\mu_{\textup{h}}}$
as an element $[\beta] \in \op{H}^1_{\textup{de Rham}}(M/T) \otimes \mathfrak{t}$,
can be rewritten as 
\begin{eqnarray} \label{easy?}
\mu_{\textup{h}}([\gamma])=\textup{exp}\, \int_{\gamma}\,\beta.
\end{eqnarray}
The $T$\--invariant connection $\Omega$ of symplectic orthogonal complements
to the tangent spaces to the $T$\--orbits defines a $T$\--invariant 
connection one--form $\theta$ for the canonical projection mapping
$\pi \colon M \to M/T$, which, by $T$\--invariance of 
$\theta$, is identified with a connection one--form 
$\widehat{\theta}\in \Omega^1(M/T) \otimes
\mathfrak{t}$. By Proposition \ref{easynow}, the connection $\Omega$ is flat,
and since the torus $T$ is an abelian group, $\widehat{\theta}$ is a closed form.
The one--form $\Theta:=\widehat{\theta} - \beta \in \Omega_{\textup{de 
Rham}}^1(M/T) \otimes \mathfrak{t}$ is a $\mathfrak{t}$\--valued connection one--form on $M/T$, in the
sense that its pull--back to $M$ is a connection one--form on $M$, and therefore there 
exists
a $T$\--invariant connection $\mathcal{D}$ on $M$, such that 
$\Theta=\widehat{\theta}^{\mathcal{D}}$.
Since $\beta$ is closed, and $\Omega$ is flat, the connection one-form 
$\widehat{\theta}^{\mathcal{D}}$ is flat, and hence the connection $\mathcal{D}$ on $M$ 
is flat. This means that $\mathcal{D}$ has a monodromy homomorphism $\mu^{\mathcal{D}} 
\colon \pi_1(M/T, \, p_0) \to T$, and a corresponding monodromy homomorphism 
on homology $\mu^{\mathcal{D}}_{\textup{h}} \colon \op{H}_1(M/T,\,\Z)
\to \mathfrak{t}$, and for any closed curve $\gamma$ in the orbit space $M/T$, 
$$
\mu^{\mathcal{D}}_{\textup{h}}([\gamma])=\mu_{\textup{h}}
([\gamma]) \, \, \, \textup{exp}\, \int_{\gamma}\,-\beta=
\mu_{\textup{h}}
([\gamma]) \, \, \, (\textup{exp}\, \int_{\gamma}\,\beta)^{-1}=1,
$$
where the second equality follows from the fact that 
$$
\textup{exp}\, \int_{\gamma}\,-\beta=
(\textup{exp} \, \int_{\gamma}\,\beta)^{-1}, 
$$
and the third equality follows from (\ref{easy?}), and we have proven that $\mathcal{D}$ is a 
flat connection for the orbibundle $\pi \colon M \to M/T$ whose monodromy homomorphism
on homology $\mu^{\mathcal{D}}_{\textup{h}}$ is trivial, and hence so it is 
$\mu^{\mathcal{D}}$. 
If conversely, $M$ is $T$\--equivariantly diffeomorphic to the Cartesian product $M/T \times T$, then 
there exists a $T$\--invariant flat connection $\mathcal{D}$ for the
orbibundle $\pi \colon M \to M/T$ such 
that the monodromy homomorphism for $\pi$ with respect to
$\mathcal{D}$, $\mu^{\mathcal{D}} \colon \pi_1(M/T, \, p_0) \to T$ is trivial. On the other hand,
if $[\gamma] \in \op{H}_1(M/T,\,\Z)_{\textup{T}}$, there exists a strictly positive
integer $k$ with
$k \,[\gamma]=0$ and
$
0=\langle k \,[\gamma],\, [\alpha]\rangle=k \, \langle [\gamma],\,[\alpha]\rangle,
$
which means that $\langle [\gamma], \, [\alpha] \rangle=0$ for all $[\gamma] \in 
\op{H}_1(M/T,\,
\Z)_{\op{T}}$, and hence that 
$(\mu^{\mathcal{D}}_{\textup{h}})^{\op{T}}=\mu_{\textup{h}}^{\op{T}}$.
Since the homomorphism $\mu^{\mathcal{D}}$ is trivial, 
$\mu_{\textup{h}}^{\op{T}}$ is trivial.
\end{proof}

\begin{remark}
\normalfont
Under the assumptions of Theorem \ref{mainlemma} on our symplectic manifold
$(M, \, \sigma)$, if the dimension of the orbit space $M/T$ is strictly greater than 2, 
then it may happen 
that the first integral orbifold homology group $\op{H}_1(M/T,\,  \Z)$ has 
no torsion and therefore $M$ is 
$T$\--equivariantly diffeomorphic to the Cartesian product $M/T \times T$. 
However, already in the case that $M/T$ is $4$\--dimensional there are examples 
of symplectic manifolds $M/T$ of which the the integral homology
group $\op{H}_1(M/T,\, \Z)$ has 
nontrivial torsion. For instance, in \cite[Sec.\,8]{DuPe}
we computed the fundamental group and the first integral homology group
of the whole manifold, and this computation shows that 
even in dimension $4$ there are examples $X$ where $\op{H}_1(X, \, \Z)$ 
is isomorphic to $\Z^3 \times (\Z/k\, \Z)$, where 
$k$ can be any positive integer. Taking such manifolds 
as the base space $M/T$, and a monodromy homomorphism 
which is nontrivial on the torsion subgroup $\Z/k\, \Z$, we 
arrive at an example where $M$ is not $T$\--equivariantly 
diffeomorphic to the Cartesian product $M/T\times T$.
\end{remark}

\begin{remark} \label{natural}
\normalfont
Under the assumptions of Theorem \ref{mainlemma}, 
since $\widetilde{M/T}$ is a regular covering of $M/T$, we have a diffeomorphism
$\widetilde{M/T}/\pi_1(M/T, \, p_0) \simeq M/T$ 
naturally, or in other words, the symbol 
$\simeq$ may be taken to be an equality. Hence the monodromy homomorphism 
$\mu^{\mathcal{D}} \colon \pi_1(M/T, \, p_0) 
\to T$ of a $T$\--invariant flat connection $\mathcal{D}$ 
is trivial if and only if $\widetilde{M/T} \times_{\pi_1(M/T, \, p_0)} T=M/T \times 
T$. Strictly speaking, this is not an equality, it is a $T$\--equivariant
diffeomorphism 
$$
\widetilde{M/T} \times_{\pi_1(M/T, \, p_0)} T \to \widetilde{M/T}/\pi_1(M/T,\,p_0) \times T \to M/T 
\times T, 
$$
where the
first arrow is the identity map, and the second arrow is the identity on the $T$ component 
of the Cartesian product
$(\widetilde{M/T}/\pi_1(M/T, \, p_0)) \times T$, while on the first component it is given by the map 
$[\,[\gamma]\,]_{\pi_1(M/T,\,p_0)} \mapsto \gamma(1)$. If the action of $T$ on $M$ is not free,
the same argument works by replacing the fundamental group of $M/T$ at $p_0$,
by the corresponding orbifold fundamental group.
 \end{remark}

\subsection{Alternative model}

Next we present a model for $(M, \, \sigma)$, up to
$T$\--equivariant symplectomorphisms, which does not involve the universal cover of $M/T$, but
rather a smaller cover of $M/T$.

Following Borzellino's article \cite[Def.\,6]{Borzellino} we make the following definition.

\begin{definition}
Let $\mathcal{O}$ be a smooth orbifold. 
The \emph{first integral orbifold homology $\op{H}_1^{\textup{orb}}(\mathcal{O},\, \Z)$
of a smooth orbifold  $\mathcal{O}$} is defined
as the abelianization of the orbifold fundamental group $\pi_1^{\textup{orb}}(\mathcal{O}, \,x)$
of $\mathcal{O}$ at $x$, where $x$ is any point in $\mathcal{O}$. 
\end{definition}

\begin{remark}
The definition of $\op{H}_1^{\textup{orb}}(\mathcal{O},\, \Z)$ does not depend on
the choice of the point $x \in \mathcal{O}$ in the sense that all 
abelianizations of $\pi_1^{\textup{orb}}(\mathcal{O}, \,x)$, $x \in \mathcal{O}$,
are naturally identified with each other.
\end{remark}

The Hurewicz map at the level of smooth orbifolds may be defined analogously
to the usual Hurewicz map at the level of smooth manifolds. Indeed, if $\mathcal{O}$
is a smooth orbifold, the
\emph{orbifold Hurewicz map} $h_1$, from the orbifold
fundamental group $\pi_1^{\textup{orb}}(\mathcal{O}, \, x_0)$ to 
the first orbifold homology group
$\op{H}_1^{\textup{orb}}(\mathcal{O},\, \Z)$, assigns to a 
homotopy class of a loop based at $x_0$, c.f. Hatcher's \cite[Sec.\,4.2]{Ha} or 
Spanier's \cite[Sec.\,7.4]{Spanier}, which is also a one--dimensional cycle, the 
homology class of that cycle. The mapping
$h_1$ is a homomorphism from $\pi_1^{\textup{orb}}(\mathcal{O}, \, x_0)$ to 
$\op{H}^{\textup{orb}}_1(\mathcal{O},\,\Z)$.

\begin{definition} \label{below}
Let $(M, \,\sigma)$ be a compact, connected, symplectic manifold equipped
with an effective symplectic action of a torus $T$ for which at least one, and hence
every $T$\--orbit is 
a $\op{dim}T$\--dimensional symplectic submanifold of $(M, \, \sigma)$, 
and let $h_1$ from $\pi_1^{\textup{orb}}(M/T, \, p_0)$ to 
$\op{H}_1^{\textup{orb}}(M/T,\, \Z)$ be the orbifold Hurewicz mapping.
There exists a unique homomorphism $\op{H}_1^{\textup{orb}}(M/T,\,\Z)\to T$, which
we call $\mu_{\textup{h}}$, such that 
\begin{eqnarray} \label{muh}
\mu = \mu_{\textup{h}} \circ h_1,
\end{eqnarray}
where 
$\mu \colon \pi_1^{\textup{orb}}(M/T, \, p_0) \to T$ is the monodromy homomorphism of the
connection $\Omega=\{\Omega_x=(\op{T}_x(T \cdot x))^{\sigma_x}\}$ of symplectic orthogonal
complements to the tangent spaces to the $T$\--orbits in $M$ (c.f. Proposition \ref{easynow}).
The homomorphism $\mu_{\textup{h}}$ is independent of the choice of base point $p_0$.
 \end{definition}

See Remark \ref{e2} for a description
of the ingredients involved in the following theorem. 

\begin{theorem}  \label{altermodel}
Let $(M, \,\sigma)$ be a compact, connected, symplectic manifold equipped
with an effective symplectic action of a torus $T$ for which at least one, and hence
every $T$\--orbit, is a $\op{dim}T$\--dimensional 
symplectic submanifold of $(M,\,\sigma)$. Then there exist
a $(\op{dim}M-\op{dim}T)$\--dimensional 
symplectic manifold $(N,\,\sigma^N)$,
a commutative group $\Delta$ which acts properly on $N$ with finite stabilizers
and such that $N/\Delta$ is compact, and 
a group monomorphism $\mu'_{\textup{h}} \colon \Delta \to T$,
such that the symplectic $T$\--manifold $(M,\,\sigma)$ is $T$\--equivariantly
symplectomorphic to $N \times_{\Delta} T$, where $\Delta$ acts on $N 
\times T$ by the diagonal action $x \, (y,\,t)=(x \star y^{-1},\,\mu'_{\textup{h}}(x) \cdot t)$, where 
$\star \colon \Delta \times N \to N$ denotes the action of $\Delta$ 
on $N$.
$N \times_{\Delta} T$ 
is equipped with the action of $T$ by translations
which descends from the action of $T$ by translations on the right factor of the product
$N \times T$, and the symplectic form 
induced on the quotient by the product
symplectic form $\sigma^{N} \oplus \sigma^T$ on the
Cartesian product $N \times T$. Here $\sigma^T$ is the unique
translation invariant symplectic form on $T$ induced by the antisymmetric bilinear form
$\sigma^{\mathfrak{t}}$, where $\sigma_x(X_M(x), \,Y_M(x))=\sigma^{\mathfrak{t}}(X,Y)$
for every $X, \,Y \in \mathfrak{t}$, and every $x \in M$, and $\sigma^{N}$
is inherited from the symplectic form 
$\nu$ on the orbit  space $M/T$, c.f. Lemma \ref{3.7} by means of the
covering map $N \to M/T$, c.f. Remark \ref{model2remark}.
The structure of smooth manifold for $N$
is inherited from the smooth manifold structure of the orbifold universal cover
$\widetilde{M/T}$, since $N$ is defined as the quotient
$\widetilde{M/T}/K$ where $K$ is the kernel
of the monodromy homomorphism $\mu \colon \pi_1^{\textup{orb}}(M/T,\,p_0) \to T$.
Moreover, $N$ is a regular
covering of the orbifold $M/T$ with covering group $\Delta$.
\end{theorem}

\begin{proof} We proved in Theorem \ref{tt} that for any element 
$[\gamma]\in \widetilde{M/T}$ and $t\in T$, 
the mapping defined by
$
\Phi([\gamma], \, t) := t\cdot\lambda_{\gamma}(1) \in M,
$ 
induces a induces a $T$\--equivariant symplectomorphism $\phi$ from 
the associated bundle 
$\widetilde{M/T}\times_{\pi^{\textup{orb}}_1(M/T, \, p_0)} T$ onto the symplectic
$T$\--manifold $(M, \, \sigma)$.
See Definition \ref{themodel} 
for the construction of the symplectic form and $T$\--action on
this associated bundle.
Let $K$ be the kernel subgroup of 
the monodromy homomorphism $\mu$ from the orbifold fundamental
group $\pi^{\textup{orb}}_1(M/T, \, p_0)$ into the torus $T$. The kernel
subgroup $K$ is a normal subgroup of $\pi^{\textup{orb}}_1(M/T, \, p_0)$ which contains the commutator subgroup $C$ of  $\pi^{\textup{orb}}_1(M/T, \, p_0)$. 
There is a unique homomorphism 
$\mu_c \colon \pi^{\textup{orb}}_1(M/T, \, p_0)/C\to T$ such that $\mu = 
\mu_c\circ\chi$, where $\chi$ is the quotient homomorphism from 
the orbifold fundamental group
$\pi^{\textup{orb}}_1(M/T, \, p_0)$ onto $\pi^{\textup{orb}}_1(M/T, \, p_0)/C$. 
The orbifold Hurewicz map 
$h_1$, from the orbifold fundamental group $\pi^{\textup{orb}}_1(M/T, \, p_0)$
 to the first orbifold homology group 
 $\op{H}^{\textup{orb}}_1(M/T,\, \Z)$, assigns to a 
homotopy class of a loop based at $p_0$, which is also a one--dimensional cycle, the 
homology class of that cycle, is a homomorphism from 
$\pi^{\textup{orb}}_1(M/T, \, p_0)$ onto 
$\op{H}^{\textup{orb}}_1(M/T,\,\Z)$ with kernel the commutator subgroup $C$.
Therefore the quotient group $K/C \le \pi^{\textup{orb}}_1(M/T, \, p_0)/C$ can be viewed as a subgroup of the first orbifold homology group
$$
\op{H}^{\textup{orb}}_1(M/T,\, \Z)\simeq \pi^{\textup{orb}}_1(M/T, \, p_0)/C,
$$
where the symbol $\simeq$ stands for the projection induced by the Hurewicz map $h_1$
from the orbifold fundamental group $\pi_1(M/T,\,p_0)$ into the first orbifold
homology group $\op{H}_1^{\textup{orb}}(M/T,\, \Z)$.
Let $N := \widetilde{M/T}/K$. By Lemma \ref{obvious} and Remark
\ref{e2}, $N$ is a smooth manifold (diffeomorphic to the integral manifold
of the distribution $\Omega$, c.f. Proposition \ref{easynow}). Because the universal cover
of $M/T$ is a regular orbifold covering of which the orbifold fundamental group is the
covering group, the quotient
$N$ is a regular orbifold covering of the orbit space $M/T$ with 
covering group 
$$
\Delta := \op{H}^{\textup{orb}}_1(M/T,\, \Z)/(K/C)
\simeq \pi^{\textup{orb}}_1(M/T, \, p_0)/C. 
$$
The group $\Delta$ is a commutative group. 
Let $\Delta$ acts on the torus $T$ by the mapping 
\begin{eqnarray} \label{mu'h}
(x, \,t) \mapsto \mu'_{\textup{h}}(x) \, t,
\end{eqnarray}
where the mapping $\mu'_{\textup{h}} \colon \Delta \to T$ is the quotient homomorphism
induced by $\mu_{\textup{h}}$ and $\pi$, where recall that 
$\mu_{\textup{h}}$ denotes the homomorphism induced on homology from 
the monodromy $\mu$ associated to the connection $\Omega$, c.f. Proposition \ref{easynow}.
The mapping $\mu'_{\textup{h}}$ is injective because the subgroup
$K$ equals the kernel of $\mu$. Our construction produces an identification of 
the bundle $\widetilde{M/T} 
\times_{\pi^{\textup{orb}}_1(M/T, \, p_0)} T$ and
$N \times_{\Delta} T$, where $\Delta$ is acting by the diagonal
action on $N \times T$ to give rise to $N \times_{\Delta} T$ (c.f
expression (\ref{mu'h})),
which intertwines the actions of $T$ by translation
on the right $T$\--factor of both spaces. In this way the mapping
$\phi$ induces a diffeomorphism from 
the associated bundle $N \times_{\Delta}T$ to $M$,
which intertwines the action of $T$ by translations on the right factor of
$N \times_{\Delta} T$ with the action of $T$ on $M$. 
By the same proof as in Theorem \ref{mainlemma}, the $T$\--equivariant 
diffeomorphism $\phi$ from $N \times_{\Delta} T$ onto $M$ pulls back the 
symplectic form on $M$ to the symplectic form on $N 
\times_{\Delta} T$ given in Remark \ref{model2remark}.
\end{proof}

\begin{remark} \label{model2remark}
This remark justifies why the symplectic form on the model space defined in Theorem
\ref{altermodel} is correctly defined. Let us assume the terminology used in the statement
of Theorem \ref{altermodel}.
The pull--back of the $2$\--form $\nu$ on $M/T$ 
to the smooth manifold $N$ by means of the covering map 
$\phi$ (c.f. Lemma \ref{3.7})
such that $\pi^* \nu|_{\Omega_x}=\sigma|_{\Omega_x}$ for every $x \in M$,
where $\Omega:=\{\Omega_x\}_{x \in M}$ is the distribution on $M$ of
symplectic orthogonal complements to the tangent spaces to the $T$\--orbits,
defined by Proposition \ref{easynow},  is a $\Delta$\--invariant symplectic form on $N$. 
The symplectic form on 
$T$ determined by the antisymmetric bilinear form $\sigma^{\mathfrak{t}}$
given by Lemma \ref{sigmat} is translation invariant,
and therefore $\Delta$\--invariant. The direct sum of the 
symplectic form on $N$ and the symplectic form on $T$ is a $\Delta$\--invariant 
and $T$\--invariant symplectic form on $N\times T$, and therefore there is a 
unique symplectic form on $N\times_{\Delta} T$ of which the pull-back by the 
covering map $N \times T \to N \times_{\Delta} T$ is equal to the 
given symplectic form on $N \times T$.
\end{remark}

\normalfont

\section{Classification: free case}

Throughout this 
section $(M, \, \sigma)$  is a compact and connected symplectic manifold and 
$T$ is a torus which acts freely on $(M, \, \sigma)$ by means of
symplectomorphisms, and such that at least one $T$\--orbit is a
$\op{dim}T$\--dimensional symplectic submanifold of $(M, \, \sigma)$. 
Our goal is to use the model for $(M, \, \sigma)$ which we constructed in Definition \ref{themodel}
to provide a classification of $(M, \, \sigma)$ when $\op{dim}T=\op{dim}M-2$,
in terms of a collection of invariants.

\subsection{Monodromy Invariant}
We define what we call the free monodromy invariant of $(M, \, \sigma)$,
ingredient 4) in Definition \ref{Fdef2}.

\subsubsection{Intersection Forms and Geometric Maps}

Let $\Sigma$ be a compact, connected, orientable, smooth surface of
genus $g$, where $g$ is a positive integer. Recall the \emph{algebraic intersection number}
\begin{eqnarray} \label{if0}
\cap \colon
\op{H}_1(\Sigma,\, \Z) \otimes \op{H}_1(\Sigma, \, \Z)
\to \Z,
\end{eqnarray}
which extends uniquely to the \emph{intersection form}
\begin{eqnarray} \label{if}
\cap \colon
\op{H}_1(\Sigma,\, \R) \otimes \op{H}_1(\Sigma, \, \R) \to \R,
\end{eqnarray}
which turns $\op{H}_1(\Sigma,\, \R)$ into a symplectic vector space.
It is always possible to find, c.f. \cite[Ex. 2A.2]{Ha}, a so called ``symplectic'' basis, in the sense of
\cite[Th.\,2.3]{MS}.

\begin{definition}
Let $\Sigma$ be a compact, connected, orientable, smooth surface of
genus $g$, where $g$ is a positive integer.
A collection of elements
$\alpha_i,\, \beta_i$, $1 \le i \le g$, of $\op{H}_1(\Sigma, \, \Z) \subset \op{H}_1(\Sigma, \, \R)$
such that
\begin{eqnarray} \label{Fcap}
\alpha_i \cap \alpha_j=\beta_i \cap \beta_j=0, \, \, \, \alpha_i \cap \beta_j=\delta_{ij}
\end{eqnarray}
for all $i, \, j$ with $1 \le i, \, j \le g$ is called a \emph{symplectic basis of the group}
$\op{H}_1(\Sigma, \, \Z)$ or a \emph{symplectic basis of the symplectic
vector space} $(\op{H}_1(\Sigma,\, \R), \, \cap)$ .
\end{definition}

Hence the matrix associated to the antisymmetric bilinear
form $\cap$ on the basis $\alpha_i, \, \beta_i$ is the block diagonal matrix
\[ \textup{J}_0=\left( \begin{array}{ccccccc}
 0 & 1 &               &  & \\
-1 & 0 &   &   & \\
 & &    0 &1  &  &&\\
&  &   -1& 0  & &&\\
&   &  &  & \ddots &&\\
& &      &&                & 0 &1 \\
& &      &&              &-1 & 0
\end{array} \right).\]
Notice that the symplectic linear group $\textup{Sp}(2g, \, \R) \subset
\textup{GL}(2g, \, \R)$ is the group of matrices
$A$ such that
\begin{eqnarray} \label{Fsp2gz}
A \cdot \textup{J}_0 \cdot A^{\textup{t}}=\textup{J}_0,
\end{eqnarray}
and therefore it is natural to denote by $\textup{Sp}(2g, \, \Z)$ the group
of matrices $A \in \textup{GL}(2g, \, \Z)$ which satisfy expression (\ref{Fsp2gz}).

\begin{remark}
The group $\textup{Sp}(2g, \, \Z)$ is also called Siegel's modular group, and denoted
by $\Gamma_g$, see for example J. Birman's article \cite{birman}. Generators
for $\textup{Sp}(2g, \, \Z)$ were first determined by L. K. Hua and and I. Reiner \cite{rk}. Later
H. Klingen found a characterization of \cite{klingen} for $g \ge 2$ by a finite system of
relations. Birman's article \cite{birman} reduces Klingen's article to a more usable form,
in which she explicitly describes the calculations in Klingen's paper, among other results.
I thank J. McCarthy for making me aware of J. Birman's article and for pointing me
to his article \cite{mc} which contains a generalization of it.
\end{remark}

\begin{remark}
Let $g$ be a non--negative integer, and let $\Sigma$ be a compact, connected,
oriented, smooth surface of genus $g$. 
If the first homology group $\op{H}_1(\Sigma, \, \Z)$ is identified with
$\textup{GL}(2g, \, \Z)$ by means of a choice of symplectic basis, the
group of  automorphisms 
of $\op{H}_1(\Sigma, \, \Z)$
which preserve the intersection form gets identified wih
the group $\textup{Sp}(2g, \, \Z)$ of matrices in $\textup{GL}(2g, \, \Z)$
which satisfy expression (\ref{Fsp2gz}).
Indeed, 
let $\alpha_i, \, \beta_i$  be a symplectic
basis of $\op{H}_1(\Sigma, \, \Z)$ (i.e. whose elements satisfy
expression (\ref{Fcap})).
The $(2g \times 2g)$\--matrix $\op{M}(\alpha_i, \, \beta_i, \, f)$ of $f$ 
with respect to the basis $\alpha_i, \, \beta_i$
is an element of the linear group $\textup{GL}(2g, \, \Z)$ of $(2g \times 2g)$\--invertible
matrices with integer coefficients and it  is an exercise to check that $f$ preserves the intersection
form on $\op{H}_1(\Sigma, \, \R)$ if and only if
$
\op{M}(\alpha_i, \, \beta_i, \, f) \cdot \textup{J}_0 \cdot \op{M}(\alpha_i, \, \beta_i, \, f))^{\textup{t}}=\textup{J}_0.
$
\end{remark}

An orientation preserving diffeomorphism induces an isomorphism in
homology which preserves the intersection 
form. Moreover,
the converse also holds, 
c.f. \cite[pp. 355-356]{MKS}. An algebraic proof of this result is given by J. Birman in \cite[pp. 66--67]{birman} as a consequence
of \cite[Thm. 1]{birman} (see also the references therein). To be self contained we present a sketch of proof here, following a preliminary draft by B. Farb and D. Margalit \cite{farbmargalit}. 
 Denote by
$\op{MCG}(\Sigma)$ the mapping class group of orientation preserving diffeomorphisms
of $\Sigma$ modulo isotopies.  After a choice of basis, there is a natural
homomorphism
\begin{eqnarray} \label{mu}
\mu \colon \op{MCG}(\Sigma) \to \op{Sp}(2g,\,\Z).
\end{eqnarray}

\begin{lemma} \label{triviallemma}
Let $g$ be a non--negative integer, and let $\Sigma, \, \Sigma'_{}$ be compact, connected,
oriented, smooth surfaces of the same genus $g$. Then a group isomorphism
$f$ from the first homology group $\op{H}_1(\Sigma, \, \Z)$
onto the first homology group $\op{H}_1(\Sigma'_{}, \, \Z)$
preserves the intersection form if and only if there exists an orientation preserving 
diffeomorphism $i \colon \Sigma \to \Sigma'_{}$
such that $f=i_*$. 
Moreover, if $\{U_k\}_{k=1}^m \subset \Sigma$, $\{U'_k\}_{k=1}^m \subset \Sigma'$ are
finite disjoint collections of embedded disks, $i$ can
be chosen to map $U_k$ to $U'_k$ for all $k$, $1 \le k \le m$.
\end{lemma}

\begin{proof}
W.l.o.g. assume that $\Sigma=\Sigma'$.
It is immediate that an orientation preserving diffeomorphism of $\Sigma$ induces an intersection
form preserving automorphism of $\op{H}_1(\Sigma, \, \Z)$.
 Let $\gamma$ be the permutation of $\{1, \ldots, 2g\}$
which transposes $2i$ fand $2i-1$ for $1 \le i \le 2g$. We define the $ij^{\textup{th}}$
elementary symplectic matrix by
$$
\Sigma^{ij} = \left\{ \begin{array}{ll}
        I_{2g} +E_{ij}& \mbox{if $i=\gamma(j)$};\\
        I_{2g}+E_{ij}-(-1)^{i+j}E_{\gamma(i)\gamma(j)}& \mbox{otherwise}.\end{array} \right.
$$
where $I_{2g}$ stands for the identity matrix of dimension $2g$ and
$E_{ij}$ is the matrix with a $1$ in the $ij^{\textup{th}}$ position and
$0's$ elsewhere. It is a classical fact 
that $\op{Sp}(2g, \, \mathbb{Z})$
is generated by the matrices $\Sigma^{ij}$. To prove the lemma is equivalent
to showing that $\mu$ in (\ref{mu}) is surjective onto $\op{Sp}(2g,\,\Z)$.
Let $\tau_b$ be the Dehn twist
about a simple closed curve $b$. 
Then for integer values of $k$, the image
$
\mu(\tau^k_b)
$
is given by
\begin{eqnarray} \label{hatform2}
a \mapsto a + k \cdot \cap(a, \,b) \, b.
\end{eqnarray}
We may restrict our attention to a subsurface $\Sigma_0$ of $\Sigma$
of genus $1$ in the case of $i=\gamma(j)$, or of genus 2 otherwise, 
such that the $i^{\op{th}}$ and $j^{\op{th}}$ basis
elements are supported on it, as well as assume that $i$ is odd,
and that $\op{H}_1(\Sigma_0,\,\Z)$ is spanned by a symplectic basis
$\alpha_1,\, \beta_1, \, \alpha_2, \, \beta_2$.
Then using that $\mu$ is a homomorphism and (\ref{hatform2}) one shows 
$\mu(\tau^{-1}_{\alpha_1})=\Sigma^{1,2}$,
$\mu(\tau_{\alpha_2}^{-1}\tau_{\alpha_1}^{-1}\tau_{\alpha_1+\beta_2})=\Sigma^{1,3}$
and 
$\mu(\tau_{\alpha_2} \tau_{\alpha_1} \tau^{-1}_{\alpha_1+\beta_2})=\Sigma^{3,2}$. 
\end{proof}

\subsubsection{Construction}

In this section we construct the monodromy invariant of our symplectic manifold
$(M, \, \sigma)$, c.f. Definition \ref{Ffreeinvariant}. The fact that such invariant
is well defined uses only linear algebra.
To any group homomorphism 
$$
f \colon \op{H}_1(\Sigma, \, \Z) \to T, 
$$
we can assign the $2g$\--tuple 
$(f(\alpha_i),\, f(\beta_i))_{i=1}^g$,
where $\alpha_i, \,\beta_i$, $1 \le i \le g$ is a basis
of the homology group $\op{H}_1(\Sigma, \, \Z)$ satisfying formulas (\ref{Fcap}). Conversely, given a
$2g$\--tuple
$(a_1, \, b_1, \ldots, a_g, \, b_g) \in T^{2g}$, the commutativity of $T$ implies,
by the universal property of free abelian groups  \cite[Ch. I.3]{df}
that there exists a unique group homomorphism from the homology group
$\op{H}_1(\Sigma, \, \Z)$ into the torus
$T$  which sends $\alpha_i$ to $a_i$ and $\beta_j$ to $b_j$, for all values of $i, \, j$
with $1 \le i, \, j \, \le g$. Notice that this tuple depends on the choice of
basis. For each such basis we have an isomorphism of groups
\begin{eqnarray}
f \mapsto (f(\alpha_i),\, f(\beta_i))_{i=1}^g, \label{Fbasismap}
\end{eqnarray}
between the homomorphism
group $\textup{Hom}(\op{H}_1(\Sigma, \, \Z), \, T)$ and the Cartesian
product $T^{2g}$.

\begin{definition} \label{Fkey}
Let $T$ be a  torus and $m$ an even positive integer. Let
$\mathcal{H}$ be a subgroup of $\textup{GL}(m, \, \Z)$.
A matrix $A$ in $\mathcal{H}$ acts
on the Cartesian product $T^m$ by sending an $m$\--tuple $x$ to $x \cdot A^{-1}$
by identifying $x$ with a homomorphism from $\Z^{2g}$ to $T$.
We say that two $m$\--tuples
$x, \, y \in T^{m}$ are $\mathcal{H}$\--equivalent if they lie in the
same $\mathcal{H}$\--orbit, i.e. if there exists a matrix
$A \in \mathcal{H}$ such that $y=x \cdot A$.
We write $\mathcal{H} \cdot x$ for the 
$\mathcal{H}$\--orbit of $x$, and 
$T^{m}/\mathcal{H}$ for the set of all 
$\mathcal{H}$\--orbits.
 \end{definition}

\begin{lemma} \label{Fprevious0}
Let $T$ be an abelian group, let $\Sigma$ be a compact, connected, smooth orientable surface,
let $\alpha_i, \, \beta_i$ and $\alpha'_i, \, \beta'_i$ be symplectic
bases of the integral homology group
$\op{H}_1(\Sigma, \, \Z)$, and
let $f \colon \op{H}_1(\Sigma, \, \Z) \to T$
be a group homomorphism. Then there exists a matrix
in the symplectic group $\textup{Sp}(2g, \, \Z)$ which takes
the image tuple $(f(\alpha_i), \, f(\beta_i))$ to the image
tuple $(f(\alpha'_i), \, f(\beta'_i))$.
\end{lemma}

\begin{proof}
Because any two symplectic bases $\alpha_i, \, \beta_i$ and $\alpha'_i, \, \beta'_i$
of  the group
$\op{H}_1(\Sigma, \, \Z)$ are taken onto each 
other by an element of the symplectic group
$\textup{Sp}(2g, \, \Z)$, the change of basis matrix 
from the basis $\alpha_i, \, \beta_i$, to the basis $\alpha'_i, \, \beta'_i$ 
is in the group $\textup{Sp}(2g, \, \Z)$. Notice that
the $(2j-1)^{\textup{th}}$\--column of this matrix consists of the coordinates of $\alpha'_j$, 
with respect to the basis $\alpha'_i, \, \beta'_i$, and its $(2j)^{\textup{th}}$\--column 
consists of the coordinates of $\beta'_j$, and hence we have that the tuple
$(f(\alpha_i),
f(\beta_i))_{i=1}^g$ is obtained by applying such matrix to the tuple
$(f(\alpha'_i),f(\beta'_i))_{i=1}^g$.
\end{proof}

Lemma \ref{Fprevious0} shows that the assignment
\begin{eqnarray} \label{Falphabetaassignment}
f \mapsto   \textup{Sp}(2g, \Z) \cdot (f(\alpha_i), \, f(\beta_i))_{i=1}^g   
\end{eqnarray}
induced by expression (\ref{Fbasismap})
is well defined independently of the choice of basis 
$\alpha_i, \, \beta_i$, $1 \le i \le g$ of $\op{H}_1(\Sigma, \, \Z)$, as long as it is a symplectic basis.  
The following is a consequence of Lemma \ref{Fprevious0} and Definition \ref{Fkey}.

\begin{lemma} \label{Fprevious}
Let $(M,\,\sigma)$ be a compact, connected $2n$\--dimensional symplectic manifold equipped 
with a free symplectic
action of a $(2n-2)$\--dimensional torus $T$, for which at least one, and hence
every $T$\--orbit is a $(2n-2)$\--dimensional symplectic submanifold of $(M, \, \sigma)$. 
If $\alpha_i,\,\beta_i$ and 
$\alpha'_i,\,\beta'_i$, $1 \le i \le g$, are symplectic
bases of $\op{H}^{}_1(M/T, \, \Z)$, 
then for every homomorphism $f \colon
\op{H}_1(M/T, \, \Z) \to T$,  we have that the $\textup{Sp}(2g, \, \Z)$\--orbits
of the tuples $(f(\alpha_i),\, f(\beta_i))_{i=1}^g$ and
$(f(\alpha'_i),
f(\beta'_i))_{i=1}^g$ are equal.
\end{lemma}

\begin{definition} \label{Ffreeinvariant}
Let $(M, \,\sigma)$ be a compact, connected, $2n$\--dimensional symplectic manifold equipped
with a free symplectic action of a torus $T$ for which at least one, and hence every 
$T$\--orbit is a $(2n-2)$\--dimensional symplectic submanifold of $(M, \, \sigma)$. 
The \emph{free monodromy invariant of $(M, \, \sigma, \,T)$} is the $\textup{Sp}(2g, \, \Z)$\--orbit
$$
 \textup{Sp}(2g, \, \Z)  \cdot (\mu_{\textup{h}}(\alpha_i), \, 
\mu_{\textup{h}}(\beta_i))_{i=1}^g  \in T^{2g}/\textup{Sp}(2g,\,\Z),
$$
where $\alpha_i, \, \beta_i$, $1 \le i \le g$, is a basis of the 
homology group $\op{H}^{}_1(M/T, \, \Z)$ satisfying expression 
(\ref{Fcap}) and $\mu_{\textup{h}}$ is the homomorphism induced on homology by the monodromy
homomorphism 
$\mu$ of the connection of symplectic orthogonal complements to the tangent
spaces to the $T$\--orbits (c.f. Proposition \ref{easynow})
by means of the Hurewicz map (c.f. Definition \ref{Fkey}, formula (\ref{muh})).
 \end{definition}

\subsection{Uniqueness}

\subsubsection{List of ingredients of $(M,\, \sigma, \, T)$}

We start with the following definition.

\begin{definition} \label{Fdef2}
Let $(M, \,\sigma)$ be a compact connected $2n$\--dimensional symplectic
manifold equipped
with a free symplectic action of a $(2n-2)$\--dimensional torus $T$ for which 
at least one, and hence every
$T$\--orbit is a $(2n-2)$\--dimensional
symplectic submanifold of $(M, \, \sigma)$. The \emph{list of ingredients of
$(M,\,\sigma,\,T)$}
consists of the following items. 
\begin{itemize}
\item[1)]
The genus $g$ of the surface $M/T$ (c.f. Remark \ref{manifoldcase}).
\item[2)]
The total symplectic area of the symplectic surface $(M/T,\, \nu)$, where 
the symplectic form $\nu$ is defined
by the condition 
$\pi^* \nu|_{\Omega_x}=\sigma|_{\Omega_x}$ for every $x \in M$,
$\pi \colon M \to M/T$ is the projection map, and 
where for each $x \in M$, $\Omega_x=(\op{T}_x(T \cdot x))^{\sigma_x}$ (c.f. Lemma \ref{3.7}).
\item[3)]
The non--degenerate antisymmetric bilinear form 
$\sigma^{\mathfrak{t}} \colon \mathfrak{t} \times \mathfrak{t} \to \R$ on the Lie
algebra $\mathfrak{t}$ of $T$ such that for all $X, \, Y \in \mathfrak{t}$
and all $x \in M$  
$\sigma_x(X_M(x),\, Y_M(x))=\sigma^{\mathfrak{t}}(X, \, Y)$ (c.f. Lemma \ref{sigmat}).
\item[4)]
The free monodromy invariant of $(M, \, \sigma, \, T)$, i.e. the
$\textup{Sp}(2g, \, \Z)$\--orbit 
$$
\textup{Sp}(2g, \, \Z) \cdot (\mu_{\textup{h}}(\alpha_i), \, 
\mu_{\textup{h}}(\beta_i))_{i=1}^g
$$
of the $2g$\--tuple $(\mu_{\textup{h}}(\alpha_i), \, 
\mu_{\textup{h}}(\beta_i))_{i=1}^g \in T^{2g}$, c.f. Definition \ref{Ffreeinvariant}.
\end{itemize} 
\end{definition}

\subsubsection{Uniqueness Statement}

The next two results say that 
the list of ingredients of $(M,\,\sigma,\,T)$ as in Definition
\ref{Fdef2} is a complete set of invariants of $(M,\,\sigma,\,T)$. We start with
a preliminary remark.

\begin{remark} \label{Ftrivialremark}
\normalfont
Let $T$ be a torus. Let $x_0 \in X$. Let $p \colon X \to B$, $p' \colon X' \to B'$ be smooth principal
$T$\--bundles, equipped with flat connections $\Omega, \, \Omega'$, and
let $\Phi \colon X \to X'$ be a $T$\--bundle isomorphism such that $\Phi^*\Omega'=\Omega$
(i.e. $\Omega'_{\Phi(x)}=\op{T}_x\Phi(\Omega_x)$ for each $x \in X$). Let 
$x'_0:=\Phi(x_0)$, and
let $\widehat{\Phi} \colon B \to B'$ be induced by $\Phi$ and such that $\widehat{\Phi}(b_0)=b'_0$,
where $p(x_0)=b_0$, $p'(x'_0)=b'_0$. 
The monodromy homomorphism $\mu \colon \pi^{}_1(B, \, b_0) \to T$ associated 
to the
connection $\Omega$ is the unique homomorphism such 
that 
\begin{eqnarray} \label{Fexpression}
\lambda_{\gamma}(1)=\mu([\gamma]) \cdot x_0,
\end{eqnarray}
for every path $\gamma \colon [0, \, 1] \to B$ such that $\gamma(0)=\gamma(1)=b_0$, 
where
$\lambda_{\gamma}$ is the unique horizontal lift of $\gamma$ with respect to 
the connection $\Omega$, such 
that $\lambda_{\gamma}(0)=x_0$.

Applying $\Phi$ to both sides of expression (\ref{Fexpression}),
and using that $\Phi$ preserves the horizontal subspaces, denoting by 
$\lambda'_{\gamma'}$
the unique horizontal lift with respect to $\Omega'$,
of any loop $\gamma' \colon [0, \,1] \to B'$ with 
$\gamma'(0)=\gamma'(1)=b'_0$,
we obtain that
$
\lambda'_{\widehat{\Phi}\, \gamma}(1)=(\mu \circ (\widehat{\Phi}_*)^{-1})([\widehat{\Phi}\, 
\gamma]) \cdot x'_0,
$
where
$\widehat{\Phi}_* \colon \pi_1(B, \, b_0) \to
\pi_1(B', \,b'_0)$ is the isomorphism induced by $\widehat{\Phi}$, which by the uniqueness
of the monodromy homomorphism, implies that $\mu'=\mu \circ (\widehat{\Phi}_*)^{-1} \colon
\pi_1(B', \, b_0) \to T$, and since $\mu=\mu_{\textup{h}} \circ h_1$ and 
$\mu'=\mu_{\textup{h}} \circ h_1$, see
expression (\ref{muh}), that $\mu'_{\textup{h}}=\mu_{\textup{h}} \circ (\widehat{\Phi}_*)^{-1} \colon
\pi^{}_1(B', \, b_0) \to T$, where in this case $\widehat{\Phi}_* \colon \op{H}^{}_1(B, \, \Z)
\to \op{H}^{}_1(B', \, \Z)$ is the homomorphism induced by $\widehat{\Phi}$ in homology.
 \end{remark}

\begin{lemma} \label{Fmt1part1}
Let $(M, \,\sigma)$ be a compact connected $2n$\--dimensional symplectic
manifold equipped
with a free symplectic action of a $(2n-2)$\--dimensional torus $T$ for which
at least one, and hence every
$T$\--orbit is a $(2n-2)$\--dimensional symplectic submanifold of $(M, \, \sigma)$.
If $(M',\,\sigma')$ is 
a compact connected $2n$\--dimensional symplectic
manifold equipped
with a free symplectic action of $T$ for which at least one, and hence every
$T$\--orbit is a $(2n-2)$\--dimensional symplectic submanifold of $(M', \, \sigma')$, 
and $(M', \, \sigma')$ is
$T$\--equivariantly symplectomorphic to $(M, \, \sigma)$, then 
the list of ingredients of $(M',\,\sigma', \, T)$ is equal to the list of ingredients of
$(M,\,\sigma, \, T)$.
\end{lemma}

\begin{proof}
Let $\Phi$ be a $T$\--equivariant
symplectomorphism from $(M, \, \sigma)$ onto $(M', \, \sigma')$. Like in the proof of
Lemma \ref{orbifoldsymplectomorphic}, the mapping $\Phi$ descends to a
symplectomorphism $\widehat{\Phi}$ from  the orbit space $(M/T, \, \nu)$ onto the
orbit space $(M'/T,\, \nu')$. By Remark \ref{manifoldcase}, 
the orbit spaces $M/T$ and $M'/T$ are 
compact, connected, orientable, smooth surfaces, and because they are diffeomorphic,
$M/T$ and $M'/T$ must have the same genus $g$.

Lemma \ref{orbifoldsymplectomorphic} and Remark \ref{moser} imply
that ingredient 2) of $(M,\, \sigma)$ equals ingredient 2) of $(M', \, \sigma')$.

If $X, \, Y \in \mathfrak{t}$, then the $T$\--equivariance of $\Phi$ implies that $\Phi^*(X_{M'})=X_M$,
$\Phi^*(Y_{M'})=Y_M$. In combination with $\sigma=\widehat{\Phi}^*\sigma'$, this implies,
in view of Lemma \ref{sigmat}, that
$$
\sigma^{\mathfrak{t}}(X, \, Y)=\Phi^*(\sigma'(X_{M'}, \, Y_{M'}))=
\Phi^*((\sigma')^{\mathfrak{t}}(X, \, Y))=(\sigma')^{\mathfrak{t}}(X, \, Y),
$$
where we have used in the last equation that $(\sigma')^{\mathfrak{t}}(X, \,Y)$ is a constant on $M'$.
This proves that $\sigma^{\mathfrak{t}}=(\sigma')^{\mathfrak{t}}$.
Since $\Phi^*\Omega'=\Omega$, we have that
$$
\mu'_{\textup{h}}= \mu_{\textup{h}} \circ (\widehat{\Phi}_*)^{-1},
$$ 
as mappings from the orbifold homology group $\op{H}^{}_1(M'/T, \, \Z)$ into the torus $T$, where the mapping $\widehat{\Phi}_*$ from the homology group $\op{H}^{}_1(M/T, \, \Z)$
to $\op{H}^{}_1(M'/T, \, \Z)$ is the group isomorphism induced by 
the orbifold diffeomorphism $\widehat{\Phi}$ from the orbit space $M/T$ onto $M'/T$ 
(see Remark \ref{Ftrivialremark}). By Lemma \ref{triviallemma},
$\widehat{\Phi}_*$ preserves the intersection--form, and therefore
the images under $\widehat{\Phi}_*$ of $\alpha_i,\,\beta_i$, $1 \le i \le g$,
which we will call $\alpha'_i, \, \beta'_i$,
form a symplectic basis of $\op{H}^{}_1(M'/T, \, \Z)$.
Hence $\mu_{\textup{h}}(\alpha_i)=\mu'_{\textup{h}}(\alpha'_i)$ and
$\mu_{\textup{h}}(\beta_i)=\mu'_{\textup{h}}(\beta'_i)$, which in turn
implies that ingredient 4) of $(M, \, \sigma)$ equals ingredient 4)
of $(M', \, \sigma')$.
\end{proof}

\begin{prop} \label{Fmt1part2}
Let $(M, \,\sigma)$ be a compact connected $2n$\--dimensional symplectic
manifold equipped
with a free symplectic action of a $(2n-2)$\--dimensional torus $T$ for 
which at least one, and hence every
$T$\--orbit is a $(2n-2)$\--dimensional symplectic submanifold of $(M, \sigma)$.
Then if $(M',\,\sigma')$ is 
a compact connected $2n$\--dimensional symplectic
manifold equipped
with a free symplectic action of $T$ for which at least one, and hence every
$T$\--orbit is a $(2n-2)$\--dimensional 
symplectic submanifold of $(M', \, \sigma')$, and the list of ingredients of $(M,\,\sigma,\,T)$ is 
equal to the list of ingredients of $(M',\,\sigma',\,T)$, then
$(M, \, \sigma)$ is $T$\--equivariantly symplectomorphic to $(M', \, \sigma')$. 
\end{prop}

\begin{proof}
Suppose that the list of ingredients of $(M, \, \sigma)$ equals the list of
ingredients of $(M', \, \sigma')$.
Let $\alpha_i, \, \beta_i$, and $\alpha'_i, \, \beta'_i$ 
be, respectively, symplectic bases 
of the first integral homology groups 
$\op{H}^{}_1(M/T, \, \Z)$ and
$\op{H}^{}_1(M'/T, \, \Z)$, and suppose 
that ingredient 4)  of $(M, \, \sigma)$ is equal to
ingredient 4)  of $(M', \, \sigma')$. 
Let $\mu, \, \mu',\, \mu_{\textup{h}}, \, \mu_{\textup{h}}'$ be the corresponding
homomorphisms respectively associated to $(M, \, \sigma)$, $(M', \, \sigma')$ as
in Definition \ref{below}.
Then, by  Definition \ref{Fkey}, there exists a matrix in the
integer symplectic group $\textup{Sp}(2g, \, \Z)$  which takes
the tuple of images of the $\alpha_i, \, \beta_i$ under $\mu_{\textup{h}}$
to the tuple of images of $\alpha'_i, \, \beta'_i$ under $\mu'_{\textup{h}}$,
which means that 
\begin{eqnarray} \label{FfreepartequationF}
\mu_{\textup{h}}=\mu'_{\textup{h}} \circ G,
\end{eqnarray}
as maps from the homology group
$\op{H}_1(M/T, \, \Z)$ to the torus $T$,
where the mapping $G$ is the intersection form preserving automorphism from 
the homology group
$\op{H}^{}_1(M/T, \, \Z)$ 
to $\op{H}^{}_1(M'/T, \, \Z)$
whose matrix with respect to the bases $\alpha_i, \, \beta_i$
and $\alpha'_i, \, \beta'_i$ is precisely the aforementioned matrix. Since $M/T$
and $M'/T$ have the same genus, by Lemma \ref{triviallemma}, there exists a surface diffeomorphism $F \colon M/T \to M'/T$
such that $F_*=G$, and hence by (\ref{FfreepartequationF}), 
\begin{eqnarray} \label{monorelation}
\mu=\mu' \circ F_*.
\end{eqnarray}

 Let $\nu$ and $\nu'$ be the symplectic 
forms given by Lemma \ref{3.7}. If $\nu_0:=F^*\nu'$, the symplectic manifold
$(M/T, \, \nu_0)$ is symplectomorphic to $(M'/T, \, \nu')$,
by means of $F$. Let $\widetilde{\nu_0}$
be the pullback of the $2$\--form $\nu_0$ by the universal cover $\psi\colon \widetilde{M/T} \to M/T$
of $M/T$ based at $p_0$, and similarly we define $\widetilde{\nu'}$ by means of the
universal cover $\psi' \colon \widetilde{M'/T} \to M'/T$ based at $F(p_0)$. 
Choose $x_0, \,x'_0$ such that $p_0=\psi(x_0),\,p'_0=\psi'(x'_0)$.
By standard covering space theory, the 
symplectomorphism $F$ between $(M/T, \, \nu_0)$ and $(M'/T, \, \nu')$ lifts to a unique
symplectomorphism $\widetilde{F}$ between $(\widetilde{M/T}, \, \widetilde{\nu_0})$ 
and $(\widetilde{M'/T}, \, \widetilde{\nu'})$ such that $\widetilde{F}(x_0)=x'_0$.
Now, let $\sigma^T$ be the unique translation invariant symplectic form on the torus 
$T$ which is uniquely determined by the antisymetric bilinear form $\sigma^{\mathfrak{t}}$,
which since ingredient 3) of $(M, \, \sigma)$ is equal to ingredient 3) of $(M', \, \sigma')$,
equals the antisymetric bilinear form $(\sigma')^{\mathfrak{t}}$. Then the product
$\widetilde{M/T} \times T$ is $T$\--equivariantly symplectomorphic to the product
$\widetilde{M'/T} \times T$ by means of the map 
\begin{eqnarray} \label{Fproductmap}
([\gamma],\,t) \mapsto (\widetilde{F}([\gamma]), \, t),
\end{eqnarray}
where $T$ is acting by translations on the right factor of both
spaces, $\widetilde{M/T} \times T$ is equipped with the symplectic form $\widetilde{\nu_0} \oplus
\sigma^T$, and  $\widetilde{M'/T} \times T$ is equipped with the symplectic form
$\widetilde{\nu'} \oplus \sigma^T$. 
As in Definition \ref{themodel}, let $[\delta]\in\pi^{}_1(M/T, \, p_0)$ act diagonally on the product $\widetilde{M/T} \times T$
by sending the pair $([\gamma],\, t)$ to 
$([\gamma \, \delta^{-1}], \, \mu([\delta])\, t)$, and similarly
let $[\delta']\in\pi^{}_1(M'/T, \, p_0)$ act on $\widetilde{M'/T} \times T$
by sending $([\gamma'],\, t)$ to 
$([\gamma' \, (\delta')^{-1}], \mu'([\delta'])\, t)$, hence giving rise to the
well--defined quotient spaces $\widetilde{M/T} \times_{\pi^{}_1(M/T, \,p_0)}T$
and  $\widetilde{M'/T} \times_{\pi^{}_1(M'/T, \,p'_0)}T$. 
(Here
$\mu$ and $\mu'$ are the monodromy homomorphisms respectively associated to
the connections $\Omega$ and $\Omega'$). Because the based fundamental
groups $\pi^{}_1(M/T, \, p_0)$ and
$\pi^{}_1(M'/T,\,p'_0)$ act properly and discontinuously, both the action of $T$ on the products, as well as the symplectic
form, induce well--defined actions and symplectic forms on these quotients. 
Therefore, because of expression (\ref{monorelation}), the assignment induced by the product mapping (\ref{Fproductmap})
$$
\Psi \colon [[\gamma], \, t]_{\pi^{}_1(M/T, \, p_0)} \mapsto 
[\widetilde{F}([\gamma]), \, t]_{\pi^{}_1(M'/T, \, p'_0)},
$$
is a $T$\--equivariant symplectomorphism between the associated bundle
$\widetilde{M/T} \times_{\pi^{}_1(M/T, \,p_0)} T$
and the associated bundle $\widetilde{M'/T} \times_{\pi^{}_1(M'/T, \, p'_0)} T$. 
Because $M/T$ and $M'/T$
have the same genus  and symplectic area, by Moser's theorem, 
c.f. \cite[Th.\,3.17]{MS}, the (compact, connected, smooth, orientable)
surfaces $(M/T, \, \nu)$ and $(M'/T, \, \nu')$ are symplectomorphic,
and hence $(M/T,\,\nu)$ is symplectomorphic to $(M/T,\,\nu_0)$.
By Theorem \ref{tt}, 
$(M, \, \sigma)$ is $T$\--equivariantly symplectomorphic to the associated bundle $\widetilde{M/T} 
\times_{\pi^{}_1(M/T, \, p_0)} T$, with the symplectic form $\widetilde{\widetilde{\nu} \oplus
\sigma^T}$,
by means of a $T$\--equivariant symplectomorphism $\phi$, and hence
$T$\--equivariantly symplectomorphic to the associated bundle $\widetilde{M/T} 
\times_{\pi^{}_1(M/T, \, p_0)} T$ with the symplectic form $\widetilde{\widetilde{\nu_0} \oplus
\sigma^T}$, say by means of $\widetilde{\phi}$.
Similarly, $(M', \, \sigma')$ is
$T$\--equivariantly symplectomorphic to $\widetilde{M'/T} 
\times_{\pi^{}_1(M'/T, \, p'_0)} T$, by means of a $T$\--equivariant
symplectomorphism $\varphi$, and therefore the composite
$\varphi \circ \Psi \circ \widetilde{\phi}^{-1} \colon M \to M'$
is a $T$\--equivariant symplectomorphism between $(M, \, \sigma)$ and $(M', \, \sigma')$.
\end{proof}

\subsection{Existence}

\subsubsection{List of ingredients for $T$}

We start by making an abstract
list of ingredients which we associate to a torus $T$.

\begin{definition} \label{Fdef1}
Let $T$ be a torus. The \emph{list of ingredients for $T$} consists of the following items.
\begin{itemize}
\item[i)]
A non--negative integer  $g$.
\item[ii)]
A positive real number $\lambda$.
\item[iii)]
An non--degenerate antisymmetric bilinear form $\sigma^{\mathfrak{t}}$ on the Lie
algebra $\mathfrak{t}$ of $T$.
\item[iv)]
A $\textup{Sp}(2g, \, \Z)$\--orbit $\gamma \in T^{2g}/
\textup{Sp}(2g, \, \Z)$, where $\textup{Sp}(2g,\, \Z)$
denotes the group of $2g$\--dimensional square symplectic matrices with integer entries, 
c.f. Definition \ref{Fkey}.
\end{itemize}
\end{definition}

\subsubsection{Existence Statement}

Any list of ingredients as in Definition \ref{Fdef1} gives rise to one of our manifolds with symplectic
$T$\--action. 

\begin{prop} \label{Fexistence}
Let $T$ be a $(2n-2)$\--dimensional torus. Then given a list of ingredients for $T$, as in Definition \ref{Fdef1}, there exists a symplectic $2n$\--dimensional manifold $(M,\, \sigma)$ with a
free symplectic action of $T$ for which at least one, and hence every 
$T$\--orbit is $(2n-2)$\--dimensional symplectic submanifold of $(M, \, \sigma)$, and such that the list of ingredients of $(M,\,\sigma,\,T)$ in Definition \ref{Fdef2} is equal to the
list of ingredients for $T$ in Definition \ref{Fdef1}.
\end{prop}

\begin{proof} 
Let $\mathcal{I}$ be a list of ingredients for the torus $T$, as in Definition \ref{Fdef1}.
Let the pair $(\Sigma, \, \sigma^{\Sigma})$ be a compact, connected
symplectic surface of genus $g$
given by ingredient i) of $\mathcal{I}$ in Definition \ref{Fdef1}, and
with total symplectic area equal to the positive real number $\lambda$ given
given by ingredient ii) of $\mathcal{I}$. Let the space $\widetilde{\Sigma}$
be the universal cover of the surface $\Sigma$ 
based at an arbitrary regular point $p_0 \in \Sigma$, which we fix for the
rest of the proof.
Let $T$ be the $(2n-2)$\--dimensional torus that we started with,  
equipped with the unique $T$\--invariant symplectic form $\sigma^T$ on $T$
whose associated non--degenerate antisymmetric bilinear form is $\sigma^{\mathfrak{t}} \colon
\mathfrak{t} \times \mathfrak{t} \to \R$, given by ingredient iii) of $\mathcal{I}$. 
Write ingredient
iv) of Definition \ref{Fdef1} as $\gamma= \textup{Sp}(2g, \, \Z) \cdot  (a_i, \, b_i)_{i=1}^g \in T^{2g}/
\textup{Sp}(2g,\, \Z)$,
and recall the existence \cite[Ex. 2A.1]{Ha}
of a symplectic $\Z$\--basis $\alpha_i, \, \beta_i$, $1 \le i \le g$,
of the integral 
homology group $\op{H}^{}_1(\Sigma, \, \Z)$, satisfying
expression (\ref{Fcap}).
Let $\mu_{\textup{h}}$ be the unique homomorphism 
from $\op{H}^{}_1(\Sigma, \, \Z)$ into $T$ such that
$$
\mu_{\textup{h}}(\alpha_i):=a_i, \,\,\,\, \mu_{\textup{h}}(\beta_i):=b_i, 
$$
for all $1 \le i \le 2g$. Define $\mu:=\mu_{\textup{h}} \circ h_1$,
where $h_1$ denotes the Hurewicz homomorphism from
$\pi^{}_1(\Sigma, \, p_0)$ to $\op{H}^{}_1(\Sigma, \, \Z)$.
$\mu$ is a homomorphism from $\pi_1^{}(\Sigma,\,p_0)$ into $T$.
Let the fundamental group 
$\pi^{}_1(\Sigma, \, p_0)$ act on the Cartesian product 
$\widetilde{\Sigma} \times T$  by 
the diagonal action  $[\delta] \cdot ([\gamma],\, t)=([\delta \, \gamma^{-1}],\, \mu([\delta]) \, t)$.
We equip the universal cover
$\widetilde{\Sigma}$ with the symplectic form
$\sigma^{\widetilde{\Sigma}}$ obtained as the pullback of
$\sigma^{\Sigma}$ by the universal covering mapping 
$\widetilde{\Sigma} \to
\Sigma$, and we
equip the product space $\widetilde{\Sigma} \times T$ with the product symplectic form
$\sigma^{\widetilde{\Sigma}} \oplus \sigma^T$, and let the torus
$T$ act by translations on the right factor of $\widetilde{\Sigma} \times T$. 
Because $\pi^{}_1(\Sigma, \, p_0)$ 
is acting properly and discontinuously, the symplectic form on $\widetilde{\Sigma} \times T$
passes to a unique symplectic 
$\sigma^{\Sigma}_{\textup{model}}:=\widetilde{\sigma^{\widetilde{\Sigma}} \oplus \sigma^T}$ 
form on $M^{\Sigma}_{\textup{model}}$ (as in the proof of Theorem \ref{tt}). 
Similarly, the action of the torus $T$ by translations
on the right factor of $\widetilde{\Sigma} \times T$ passes to an action of $T$ on 
$M^{\Sigma}_{\textup{model}}$, 
which is free. Define the associated bundle
\begin{eqnarray} \label{Fmodell}
M^{\Sigma}_{\textup{model}\,p_0}:=
\widetilde{\Sigma} \times_{\pi_1(\Sigma, \, p_0)}T.
\end{eqnarray}
It follows from the construction that 
$(M^{\Sigma}_{\textup{model}}, \, \sigma^{\Sigma}_{\textup{model}})$
is a compact, connected symplectic manifold, with a free $T$\--action for which every 
$T$\--orbit is a $\op{dim}T$\--dimensional symplectic submanifold 
of $(M^{\Sigma}_{\textup{model}}, \, \sigma^{\Sigma}_{\textup{model}})$.

It is left to show that the list $\mathcal{I}$ equals the list of ingredients of 
$(M^{\Sigma}_{\textup{model}}, \, \sigma^{\Sigma}_{\textup{model}})$. 
Since the action of the torus 
$T$ on $M^{\Sigma}_{\textup{model}}$ is induced by the action of $T$ on the right factor of 
$\widetilde{\Sigma} \times T$, $M_{\textup{model}}/T$ is symplectomorphic to 
$\widetilde{\Sigma}/\pi^{}_1(\Sigma, \, p_0)$ with the 
symplectic form induced by the mapping $\widetilde{\Sigma}
\to \widetilde{\Sigma}/\pi_1(\Sigma, \, p_0)$, which by construction
(i.e. by regularity of the orbifold universal cover)
is symplectomorphic to the quotient 
$(\Sigma, \, \sigma^{\Sigma})$. 
Therefore $M^{\Sigma}_{\textup{model}}/T$
with the symplectic form $\nu_{\Sigma}$ given by Lemma \ref{3.7}, is
symplectomorphic to $(\Sigma, \, \sigma^{\Sigma})$, and in particular 
the total symplectic area of $(M^{\Sigma}_{\textup{model}}/T,  \, \nu_{M^{\Sigma}_{\textup{model}}}/T)$
equals the total symplectic area of $(\Sigma,  \, \sigma^{\Sigma})$, which is equal to the positive
real number $\lambda$.
Let $p \colon \widetilde{\Sigma}
\times T \to M_{\textup{model}}$ be the projection map. 
It follows from the definition of $\sigma^{\Sigma}_{\textup{model}}$ that for
all $X, \, Y \in \mathfrak{t}$, the real number
$$
(\sigma^{\Sigma}_{\textup{model}})_{[[\gamma], \, t]_{\pi_1(\Sigma, \, p_0)}}(X_{M^{\Sigma}_{\textup{model}}}([[\gamma], \, t]_{\pi_1(\Sigma, \, p_0)}),
Y_{M^{\Sigma}_{\textup{model}}}([[\gamma], \, t]_{\pi_1(\Sigma, \, p_0)}))
$$ 
is equal to
$$
(\widetilde{\sigma^{\widetilde{\Sigma}} \oplus \sigma^T})_{[[\gamma], \, t]_{\pi_1(\Sigma, \, p_0)}}
(\op{T}_{([\gamma], \, t)}p(X_{\widetilde{\Sigma} \times T}([\gamma], \, t)), \, \op{T}_{([\gamma], \, t)}p(Y_{\widetilde{\Sigma} \times T}([\gamma], \, t))) \nonumber \\
$$
which is equal to
\begin{eqnarray}
(\sigma^{\widetilde{\Sigma}} \oplus \sigma^T)_{([\gamma], \, t)}
(X_{\widetilde{\Sigma} \times T}([\gamma], \, t)), \, Y_{\widetilde{\Sigma} \times T}([\gamma], \, t)) 
&=& \sigma_t^T(X_T(t), \, Y_T(t)) \nonumber \\
&=& \sigma^{\mathfrak{t}}(X, \, Y). \label{Feq}
\end{eqnarray}
where in the first equality of (\ref{Feq}) we have used that the vector fields 
$X_{\widetilde{\Sigma} \times T},  
\, Y_{\widetilde{\Sigma} \times T}$ are tangent to the $T$\--orbits
$\{y\} \times T$ of $\widetilde{\Sigma} \times T$, 
and the symplectic form vanishes on the orthogonal complements
$\op{T}_{(u, \,t)} (\widetilde{\Sigma} \times \{t\})$ to the $T$\--orbits. The last equality follows from the definition of $\sigma^{\mathfrak{t}}$.

Finally let 
$\Omega^{\Sigma}_{\textup{model}}$
stand for the flat connection on $M^{\Sigma}_{\textup{model}}$ given by the symplectic orthogonal
complements to the tangent spaces to the orbits of the $T$\--actions, see 
Proposition \ref{easynow}, and let
$\mu^{\Omega^{\Sigma}_{\textup{model}}}_{\textup{h}}$ stand for the induced homomorphism
$\mu^{\Omega^{\Sigma}_{\textup{model}}}_{\textup{h}} \colon \op{H}_1(M^{\Sigma}_{\textup{model}}/T, \, \Z)
\to T$ in homology. 
If $f \colon \op{H}_1(M^{\Sigma}_{\textup{model}}/T,
\, \Z) \to \op{H}_1(\Sigma, \, \Z)$ is
the group isomorphism induced by the symplectomorphism
\begin{eqnarray} \label{compositemap}
M^{\Sigma}_{\textup{model}}/T \to \widetilde{\Sigma}
/\pi^{}_1(\Sigma, \, p_0) \to \Sigma,
\end{eqnarray}
where each arrow in (\ref{compositemap}) represents the natural map,
\begin{eqnarray} \label{moneqtn1}
\mu_{\textup{h}}^{\Omega^{\Sigma}_{\textup{model}}}=\mu_{\textup{h}} \circ f.
\end{eqnarray}
Because $f$ is induced by a diffeomorphism, by Lemma \ref{triviallemma} $f$
preserves the intersection form and hence 
the unique collection of elements 
$\alpha'_i, \, \beta'_i$, $1 \le i \le g$
such that
$f(\alpha'_i)=\alpha_i$ and $f(\beta'_i)=\beta_i$, for all $1 \le i \le g$,
is a symplectic basis
of the homology group 
$\op{H}^{}_1(M^{\Sigma}_{\textup{model}}/T, \, \Z)$.
Let $\widehat{\gamma}$ be the $2g$\--tuple of elements 
$\mu_{\textup{h}}^{\Omega^{\Sigma}_{\textup{model}}}(\alpha'_i)$,
$\mu_{\textup{h}}^{\Omega^{\Sigma}_{\textup{model}}}(\beta'_i)$, $1 \le i \le g$.
Therefore by (\ref{moneqtn1})
$$
\widehat{\gamma}=(\mu_{\textup{h}}(\alpha_i),\, \mu_{\textup{h}}(\beta_i))_{i=1}^{g}.
$$
The result follows because the $\op{Sp}(2g, \, \Z)$\--orbit of $(\mu_{\textup{h}}(\alpha_i),\, \mu_{\textup{h}}(\beta_i))_{i=1}^{g}$ is equal to item 4) in Definition \ref{Fdef2}.
\end{proof}

\subsection{Classification Theorem}

We state and prove the two main results of Section 3, by putting together
previous results.

\begin{theorem} \label{FTT}
Let $T$ be a $(2n-2)$\--dimensional torus. Let $(M,\,\sigma)$ be 
a compact connected $2n$\--dimensional symplectic manifold 
on which $T$ acts freely and symplectically and such that at least one, and hence
every $T$\--orbit is a $(2n-2)$\--dimensional symplectic submanifold of $(M, \, \sigma)$.

Then the list of ingredients of $(M,\,\sigma,\,T)$ as in Definition
\ref{Fdef2} is a complete set of invariants of $(M,\,\sigma,\,T)$, 
in the sense that, if $(M',\,\sigma')$ is 
a compact connected $2n$\--dimensional symplectic
manifold equipped
with a free symplectic action of $T$ for which at least one, and hence every
$T$\--orbit is $(2n-2)$\--dimensional symplectic submanifold of $(M', \, \sigma')$, $(M', \, \sigma')$ is
$T$\--equivariantly symplectomorphic to $(M, \, \sigma)$ if and only if 
the list of ingredients of $(M',\,\sigma',\,T)$ is equal to the list of ingredients of
$(M,\,\sigma, \,T)$.

And given a list of ingredients for $T$, as in Definition \ref{Fdef1}, there 
exists a symplectic $2n$\--dimensional 
manifold $(M,\, \sigma)$ with a free torus action of $T$ for which one 
at least one, and hence every
$T$\--orbit is a $(2n-2)$\--dimensional
symplectic submanifold of $(M, \, \sigma)$, 
such that the list of ingredients of $(M,\,\sigma,\,T)$ is equal to the
list of ingredients for $T$.
\end{theorem}
\begin{proof}
It follows by putting together Lemma \ref{Fmt1part1}, Proposition \ref{Fmt1part2} and Proposition 
\ref{Fexistence}. Observe that the combination of Lemma \ref{Fmt1part1}, Proposition \ref{Fmt1part2}
gives the existence part of the theorem, while Proposition \ref{Fexistence} gives the existence
part.
\end{proof}

\begin{cor} \label{FTTT}
Let $(M, \,\sigma)$ be a compact connected symplectic $2n$\--dimensional manifold equipped
with a free symplectic action of $(2n-2)$\--dimensional torus $T$ for which at least 
one, and hence every $T$\--orbit is 
a $(2n-2)$\--dimensional symplectic submanifold of $(M,\,\sigma)$.  

Then the genus $g$ of the surface $M/T$ is a complete invariant of the $T$\--equivariant
diffeomorphism type of $M$, in the following sense. 
If $(M', \, \sigma')$ is a compact, connected, $2n$\--dimensional symplectic
manifold equipped with a free action of a $(2n-2)$\--dimensional
torus such that at least one, and hence every $T$\--orbit
is a $(2n-2)$\--dimensional symplectic submanifold of $(M', \, \sigma')$, 
$M'$ is $T$\--equivariantly diffeomorphic to $M$ if and only if the genus of the orbit space $M'/T$ is equal to the genus of $M/T$. 

Moreover, given any non--negative integer,
there exists a symplectic $2n$\--dimensional manifold $(M, \, \sigma)$ with
at least one, and hence every $T$\--orbit being a $(2n-2)$\--dimensional
symplectic submanifold of $(M, \, \sigma)$, such that the genus of $M/T$ is precisely the aforementioned integer.
\end{cor}

The proof of Corollary \ref{FTTT} is immediate. Indeed, by 
Corollary \ref{KK}, $M$ is $T$\--equivariantly diffeomorphic to the Cartesian product
$M/T \times T$ equipped with the action of $T$ by translations on the right factor
of the product.
Suppose that $\Phi$ is a $T$\--equivariant
symplectomorphism from $(M, \, \sigma)$ to $(M', \, \sigma')$. Because $\Phi$ is $T$\--equivariant, it descends to a 
diffeomorphism $\widehat{\Phi}$ from the orbit space $M/T $ onto the orbit space 
$M'/T$, and hence the genus
of $M/T$ equals the genus of $M'/T$. Conversely, suppose that $(M, \, \sigma)$
and $(M', \, \sigma')$ are such that the genus of $M/T$ equals the genus
of $M'/T$. By the classification theorem of compact, connected, orientable
(boundaryless of course) surfaces, there exists a diffeomorphism
$F \colon M/T \to M'/T$. Hence the map $M/T \times T \to M'/T \times T$
given by $(x, \, t) \to (F(x), \, t)$ is a $T$\--equivariant
diffeomorphism, and by Corollary \ref{KK} we are done. 
Now let $g$ be a non--negative integer, and let
$\Sigma$ be a (compact, connected, orientable) surface of genus $g$, and
let $T$ be a $(2n-2)$\--dimensional torus. Then $M_g:=\Sigma \times T$
is a $2n$\--dimensional manifold. Equip it with any product symplectic form, and 
with the action of $T$ by translations on the right factor. Then
the $T$\--orbits, which are of the form $\{u\} \times T$, $u \in \Sigma$,
are symplectic submanifolds of $M_g$. Clearly $M_g /T$ is diffeomorphic to $\Sigma$
by means of $[u,\,t] \mapsto u$.

\begin{remark}
We summarize Corollary \ref{FTTT}
in the language of categories, c.f. MacLane's book \cite{maclane}.
Let $T$ be a torus and let $\mathcal{M}$ denote the category of which 
the objects are the compact connected symplectic manifolds 
$(M,\,\sigma )$ together with a free symplectic 
$T$\--action on $(M,\,\sigma )$, such that at least one, and hence every $T$\--orbit
is a $(2n-2)$\--dimensional symplectic submanifold of $M$, and of which the 
morphisms are the $T$\--equivariant symplectomorphisms
of $(M, \, \sigma)$.
 Let $\mathcal{Z}^+$ denote 
category which consists of the set of non--negative
integers, and of which the identity is
the only endomorphism of categories. 
Then the assignment $M \mapsto g$, where
$g$ is the genus of the surface $M/T$,
is a full functor of categories from the
category $\mathcal{M}$ onto the category
$\mathcal{Z}^+$.
In particular
 the proper class $\mathcal{M}/ \sim$ of isomorphism classes in $\mathcal{M}$ 
is a set, and the functor $\iota \colon M \mapsto g$, where
$g$ is the genus of the surface $M/T$, 
induces a bijective mapping $\iota/\sim$ 
from the category $\mathcal{M}/\sim$ onto the category
$\mathcal{Z}^+$. 

Let $\mathcal{I}$ denote the set of all lists of 
ingredients as in Definition \ref{Fdef1}, 
viewed as a category, and of which the identities are
the only endomorphisms of categories. 
Then the assignment $\iota$ in Definition \ref{Fdef2} 
is a full functor of categories from the
category $\mathcal{M}$ onto the category
$\mathcal{I}$. 
In particular the proper class $\mathcal{M}/ \sim$ of isomorphism classes in $\mathcal{M}$ 
is a set, and the functor $\iota  \colon \mathcal{M} \to \mathcal{I}$ 
induces a bijective mapping $\iota/\sim$ 
from $\mathcal{M}/\sim$ onto $\mathcal{I}$. 
The fact that the mapping $\iota  \colon \mathcal{M} \to \mathcal{I}$ is a functor 
and the mapping $\iota/ \sim$ is injective 
follows from the uniqueness part of Theorem \ref{FTT}, while the surjectivity of $\iota$,
follows from the existence part.
\end{remark}

%%%%%%%%%%%%

\section{Classification}

This section extends the results of Section 3 to non--free actions. 

A tool needed to obtain the classification is Theorem \ref{veryhard},
which is a characterization of geometric isomorphisms of orbifold homology, c.f. 
Definition \ref{geomiso}. Such classification appears to be of independent interest
as it generalizes a classical result about smooth surfaces to smooth orbisurfaces.
Some of the statements and proofs in the present section are analogous to those of
Section 3, and we do not repeat them. 

\subsection{Geometric Torsion in Homology of Orbifolds}

Like for compact, connected, orientable smooth surfaces, there exists a classification
for compact, connected, orientable smooth orbisurfaces. 

\begin{remark}
It follows from Definition \ref{orb} that there are only finitely many points in compact, connected, smooth orientable orbifold $\mathcal{O}$ which are singular, so the singular locus $\Sigma_{\mathcal{O}}$ is finite. 
\end{remark}

The compact, connected, boundaryless, orientable smooth
orbisurfaces are classified
by the genus of the underlying surface and the 
$n$\--tuple of cone point orders $(o_k)_{k=1}^n$, where $o_i \le o_{i+1}$, for all $1 \le 
i \le n-1$, in the sense that the following 
two statements hold.

\begin{theorem} \label{2dorbifoldtheorem}
 First, given a positive integer $g$ and an $n$\--tuple $(o_k)_{k=1}^n$, $o_i \le o_{i+1}$
 of positive integers, 
there exists a compact, connected, boundaryless, orientable smooth orbisurface
orbifold $\mathcal{O}$ with underlying topological
space a compact, connected
surface of genus $g$ and $n$ cone points of respective orders $o_1,...,o_n$.
Secondly, let $\mathcal{O}, \, \mathcal{O}'$ be compact, connected, boundaryless, orientable smooth
orbisurfaces. Then $\mathcal{O}$ is diffeomorphic to $\mathcal{O}'$ if and only if the genus of their underlying surface is the same, and their 
associated increasingly ordered $n$\--tuple of orders of cone points are equal.
\end{theorem}

\begin{proof}
If $\mathcal{O}, \, \mathcal{O}'$ are compact, connected, boundaryless, orientable smooth $2$\--dimensional orbifolds which are moreover diffeomorphic, it follows from the definition of diffeomorphism
of orbifolds that the genus of their underlying surface is the same, and their associated increasingly ordered $n$\--tuple of orders of cone points are equal. 

Conversely suppose that $\mathcal{O}, \, \mathcal{O}'$ are
boundaryless, orientable smooth $2$\--dimensional orbifolds with the property that
the genus of their underlying surface is the same, and their 
associated increasingly ordered $n$\--tuples of orders of cone points are equal. 
Write $p_k$ for the cone points of
$\mathcal{O}$ ordered so that $p_k$ has order $o_k$ for all $k$, where $1 \le k \le n$.
Because  $\mathcal{O}, \, \mathcal{O}'$ are compact, connected, boundaryless, orientable 
and $2$\--dimensional, for each cone point in 
$\mathcal{O}$ and each cone point in $\mathcal{O}'$
there exists a neighborhood that is orbifold diffeomorphic to the standard 
model of the plane modulo a rotation. Therefore, since the tuples of orders of $\mathcal{O}$ and
$\mathcal{O}'$ are the same, for each $k$ there is a neighborhood $D_k$ of $p_k$,
a neighborhood $D'_k$ of $p'_k$ and a diffeomorphism $f_k \colon D_k \to D'_k$ such that
$f(p_k)=p'_k$. By shrinking each $D_k$ or $D'_k$ if necessary, we may assume that
the topological boundaries $\partial_{|\mathcal{O}|}(D_k)$, $\partial_{|\mathcal{O}'|}(D'_k)$
are simple closed curves. The
map $f \colon D:=\bigcup D_k \to D':=\bigcup D'_k$ defined by $f|_{D_k}:=f_k$ is an orbifold diffeomorphism such that $f(p_k)=p'_k$ and $f(D_k)=D'_k$ for all $k$, where $1 \le k \le n$.
Since
the topological boundaries $C:=\partial_{|\mathcal{O}|}(D_k)$, $C:=\partial_{|\mathcal{O}'|}(D'_k)$
are simple closed curves,
$|\mathcal{O}| \setminus \cup D_k$ and
$|\mathcal{O}'| \setminus \cup D'_k$ are surfaces with boundary, and their 
corresponding boundaries
consist of precisely $k$ boundary components,
$\partial_{|\mathcal{O}|}(D_1), \ldots, \partial_{|\mathcal{O}|}(D_n)$ and
$\partial_{|\mathcal{O}'|}(D'_1), \ldots, \partial_{|\mathcal{O}'|}(D'_n)$,
each of which is a circle. 
Then by the classification of surfaces with boundary, there exists a
diffeomorphism $g \colon C \to C'$. By definition of diffeomorphism, $g(\partial C)=\partial C'$,
and hence there exists a permutation $\rho$ of $\{1,\ldots,n\}$
such that $g(\partial_{|\mathcal{O}|}(D_k))=\partial_{|\mathcal{O}'|}(D'_{\rho(k)})$ for all $k$,
where $1 \le k \le n$. There exist diffeomorphisms of a surface with boundary that 
permute the boundary circles in any way one wants,  and hence 
by precomposing with
an appropiate such diffeomorphism we may assume that $\tau$ is the identity.
Because of this together with the fact that a diffeomorphism of a circle which preserves orientation is isotopic to the identity
map, we can smoothly deform $g$ near the boundary in $\mathcal{O}$ of each $D_k$ 
so that $g$ agrees with
$f$ on $\partial_{|\mathcal{O}|}(C)$. Hence the map $F \colon \mathcal{O}
\to \mathcal{O}$ given as $F|_{C}:=g$ and $F|_{\mathcal{O} \setminus C}:=f$ is
a well defined orbifold diffeomorphism between $\mathcal{O}$ and $\mathcal{O}'$.
\end{proof}

\begin{remark}
A geometric classification of orbisurfaces which considers hyperbolic, elliptic
and parabollic structures, is given by Thurston in \cite[Th.\,13.3.6]{thurston3} -- while
such statement is very interesting and complete, the most convenient classification
statement for the purpose of this paper is only in differential topological terms, c.f Theorem 
\ref{2dorbifoldtheorem}. The author is grateful to A. Hatcher
for providing him with the precise classification statement, and indicating to him
how to prove it.
\end{remark}

\begin{definition} \label{geometrictorsionbasis}
Let $\Sigma$ be a smooth, compact, connected, orientable smooth orbisurface
with $n$ singular points. Fix an order in the singular points, say $p_1,\ldots,p_n$, such
that the order of $p_k$ is less than or equal to the order of $p_{k+1}$.
We say that a collection of $n$ elements
$\{\gamma_k\}_{k=1}^n \subset \op{H}^{\textup{orb}}_1(\Sigma,\, \Z)$
is a \emph{geometric torsion basis\footnote{We use the word ``torsion'' because a geometric torsion basis
will generate the torsion subgroup of the orbifold homology group. Similarly, we use ``geometric'' because the homology classes
come from geometric elements, loops around singular points.} 
of $\op{H}^{\textup{orb}}_1(\Sigma,\, \Z)$} 
if $\gamma_k$ is the homology class of a loop
$\widetilde{\gamma}_k$ obtained as an oriented boundary 
of a closed small disk containing the $k^{\textup{th}}$ singular point 
of the orbisurface $\Sigma$ with respect to the aforementioned ordering, and where no two such disks intersect. 
 \end{definition}

\begin{definition} \label{key2}
Let $\vec{o}=(o_k)_{k=1}^n$ be an $n$\--tuple of integers.
We call $\textup{S}_n^{\vec{o}}$ the subgroup of the permutation group $\textup{S}_n$ 
which preserve the $n$\--tuple $\vec{o}$, i.e.
$
\textup{S}_n^{\vec{o}}:=\{\tau \in \textup{S}_n \, | \, (o_{\tau(k)})_{k=1}^n=
\vec{o} \}. 
$
 \end{definition}

\begin{lemma} \label{scott2}
Let $\Sigma$ be a compact, connected, orientable, boundaryless smooth orbisurface.
Assume moreover that $\Sigma$ is a good orbisurface. Then the order of any 
cone point of $\Sigma$ is equal to the order of the homotopy class of a small loop around that 
point. 
\end{lemma}
\begin{proof}
Take a cone point with cone angle $2\, \pi/n$, and let 
$\gamma$ denote the associated natural generator of $\pi_1^{\textup{orb}}(\Sigma, \, x_0)$. 
We already know that $\gamma^n=1$, 
because we have that relation in the presentation of the orbifold fundamental group. Since the 
orbifold has a manifold cover, the projection around the pre--image of the cone point is a 
$n$\--fold branched cover, which implies that for any $k<n$, $\gamma^k$ does not lift to the cover 
and so it must be nontrivial.
\end{proof}

Let $g, \, n, \, o_k $, $1 \le k \le n$, be non--negative integers, 
and let $\Sigma$ be a compact, connected, orientable smooth
orbisurface with
underlying topological space a surface of genus $g$, and with $n$ singular points
$p_k$ of respective orders $o_k$.
The orbifold fundamental group $\pi_1^{\textup{orb}}(\Sigma, \, p_0)$ has group presentation 
\begin{eqnarray} \label{fg}
\langle 
\{\alpha_i,\, \beta_i\}_{i=1}^g,  \, \{\gamma_k\}_{k=1}^n \, \, | \, \,
\prod_{k=1}^n \gamma_k=\prod_{i=1}^g [\alpha_i, \, \beta_i], \, \, \gamma_k^{o_k}=1, \, \, 1 \le k \le n
\rangle,
\end{eqnarray}
where the elements $\alpha_i, \, \beta_i$, $1 \le i \le g$ form a symplectic basis of the free homology
$\op{H}^{\textup{orb}}_1(\Sigma, \, \Z)_{\textup{F}}$,  
and the $\gamma_k$ form a geometric torsion basis of 
the orbifold homology group $\op{H}_1^{\textup{orb}}(\Sigma, \, \Z)$.
By abelianizing expression (\ref{fg}) we obtain 
the first integral orbifold homology group $\op{H}_1^{\textup{orb}}(\Sigma,
\, \Z)$
\begin{eqnarray}
\langle 
\{\alpha_i,\, \beta_i\}_{i=1}^g,  \, \{\gamma_k\}_{k=1}^n \, \, | \, \,
\sum_{k=1}^n \gamma_k=0, \, \, o_k \, \gamma_k=0, \, \, 1 \le k \le n
\rangle.
\end{eqnarray}

The torsion subgroup $\op{H}_1^{\textup{orb}}(\Sigma, \, \Z)_{\op{T}}$ 
of the first orbifold integral homology group of the orbisurface $\Sigma$ is generated
by the geometric torsion $\{\gamma_k\}_{k=1}^n$ with the relations
$o_k \, \gamma_k=0$ and $\sum_{k=1}^n \gamma_n=0$.
There are many free subgroups $F$ of the first orbifold homology group
for which $
\op{H}_1^{\textup{orb}}(\Sigma, \, \Z)= F \oplus  \op{H}_1^{\textup{orb}}(\Sigma, \, \Z)_{\op{T}}$. In what follows we will use the definition
\begin{eqnarray} \label{orbifoldhomologyquotient}
\op{H}^{\textup{orb}}_1(\Sigma, \, \Z)_{\textup{F}}:=\op{H}_1^{\textup{orb}}(\Sigma, \, \Z) /  \op{H}_1^{\textup{orb}}(\Sigma, \, \Z)_{\op{T}},
\end{eqnarray}
and we call the left hand side of expression (\ref{orbifoldhomologyquotient}), the \emph{free first orbifold homology group of 
the orbisurface $\Sigma$}; such quotient group is isomorphic to 
the free group on $2g$ generators $\Z^{2g}$,
and there is an isomorphism of groups from
$\op{H}_1^{\textup{orb}}(\Sigma, \, \Z)$ onto
$\op{H}^{\textup{orb}}_1(\Sigma, \, \Z)_{\textup{F}}
\oplus \op{H}^{\textup{orb}}_1(\Sigma, \, \Z)_{\op{T}} .$
As in (\ref{if}), there is a natural \emph{intersection form}
\begin{eqnarray} \label{if2}
\cap_{\scriptop{F}} \colon
\op{H}^{\scriptop{orb}}_1(\Sigma,\, \R)_{\scriptop{F}} 
\otimes \op{H}^{\scriptop{orb}}_1(\Sigma, \, \R)_{\scriptop{F}} \to \R,
\end{eqnarray} 
which for simplicity we write $\cap_{\scriptop{F}}=\cap$, and 
a natural isomorphism
\begin{eqnarray} \label{if3}
\op{H}^{\scriptop{orb}}_1(\Sigma,\, \R)_{\scriptop{F}} 
\to \op{H}_1(\widehat{\Sigma}, \, \R)
\end{eqnarray} 
which pullbacks $\cap$ to $\cap=\cap_{\scriptop{F}}$, where $\widehat{\Sigma}$
denotes the underlying surface to $\Sigma$.

\begin{example} 
\label{assumptioniii)example0}
Let $\Sigma$ be a compact, connected, orientable smooth orbisurface
with underlying topological space equal to a $2$\--dimensional torus $(\R/\Z)^2$. Suppose that $\Sigma$ has precisely one cone point $p_1$ of order $o_1=2$.
Let $\widetilde{\gamma}$ be obtained as boundary loop
of closed small disk containing the singular point $p_1$, and let $\gamma=[\widetilde{\gamma}]$.
Let $\alpha, \, \beta$ be representative of the standard basis of loops
of the surface underlying $\Sigma$ (recall that $\alpha, \, \beta$ are a
basis of the quotient free first orbifold homology group of $\Sigma$).
Then for any point $x_0 \in \Sigma$
$$
\pi_1^{\textup{orb}}(\Sigma,\,x_0)=\langle \alpha, \, \beta, \, \gamma \, | \,
[\alpha, \, \beta]=\gamma, \, \gamma^2=1 \rangle,
$$
and
$$
\op{H}_1^{\textup{orb}}(\Sigma,\, \Z)=\langle \alpha, \, \beta, \, \gamma \, | \,
\gamma=1, \, 2 \, \gamma=0 \rangle \simeq \langle \alpha, \, \beta \rangle.
$$
 \end{example}

\begin{remark} \label{effective}
If $(M, \, \sigma)$ is a compact connected symplectic manifold,  
$T$ is a torus which acts effectively but non--freely on it by 
symplectomorphisms and with $\op{dim}T$\--dimensional symplectic $T$\--orbits, 
there exists at least one
non--trivial element $[\gamma] \in \op{H}_1^{\textup{orb}}(M/T, \, \Z)$ for which there
is an integer $k \in \Z^+$ such that $k \, [\gamma]=0$ and $\mu_{\textup{h}}([\gamma])
\neq 0$. (c.f. Theorem \ref{mainlemma} for a global result).
\end{remark}

\subsection{Monodromy invariant}

We define the Fuchsian signature  monodromy space, whose elements give one
of the ingredients of the classification theorems. 

\subsubsection{Geometric Isomorphisms}

An orbifold diffeomorphism between orbisurfaces induces an isomorphism at the level of orbifold fundamental groups, and at the level of first orbifold homology groups.

\begin{definition} \label{geomiso}
Let $\mathcal{O}, \, \mathcal{O}'$ be compact,
connected, orientable smooth orbisurfaces. We say that an isomorphism
$\op{H}_1^{\textup{orb}}(\mathcal{O},\,\Z) \to \op{H}_1^{\textup{orb}}(\mathcal{O}',\,\Z)$
is \emph{orbisurface geometric} if 
there exists an orbifold diffeomorphism $\mathcal{O} \to \mathcal{O'}$ which induces it. 
\end{definition}

\begin{example}
\normalfont
If $(M, \, \sigma)$ is a compact and connected symplectic manifold of dimension
$2n$, and if $T$ is a $(2n-2)$\--dimensional torus which acts freely on $(M, \, \sigma)$ by means of symplectomorphisms and whose
$T$\--orbits are $(2n-2)$\--dimensional symplectic
submanifold of $(M, \, \sigma)$, then the
torsion part of the first integral orbifold homology group of the surface $M/T$
is trivial. But if the action of the torus
is not free then $M/T$ is an orbisurface, and the torsion subgroup is frequently non--trivial, although in a few cases it is trivial.  For example:
\begin{itemize}
\item[a)]
if there is only one singular point in $M/T$, or
\item[b)]
if there are precisely two cone points of orders $2$ and $3$ in $M/T$, or in general of orders
$k$ and $k+1$, for a positive integer $k$. In this case it is easy to see 
that the torsion subgroup of the first integral orbifold homology group
is trivial because it is generated by $\gamma_1, \, \gamma_2$ with the relations $k \, \gamma_1=(k+1)
\, \gamma_2=0$, and $\gamma_1+\gamma_2=0$. 
\item[c)]
if there are precisely three cone points of orders $3,\,4,\,5$ in $M/T$, or more generally any three points
whose orders are coprime. 
\end{itemize}
However,
in cases a) and b) the orbifolds are not good, so they do not arise as $M/T$, c.f. Lemma
\ref{goodorbifoldlemma}.
On the other hand, in 
case c) the orbifold is good, and as we will prove later
it does arise as the orbit space of many
symplectic manifolds $(M, \, \sigma)$. In this case it is possible to give  a 
description of the monodromy invariant
which is analogous to the free case done in Section 3. 
\end{example}

Recall that, by Lemma \ref{triviallemma}, 
if $\Sigma, \, \Sigma'$ are
compact, connected, oriented, smooth
orbisurfaces of the same  Fuchsian signature, and if
$f \colon \op{H}_1^{\textup{orb}}(\Sigma, \, \Z) \to \op{H}_1^{\textup{orb}}(\Sigma', \, \Z)$ is a group
isomorphism for which there exists an orientation preserving orbifold diffeomorphism 
$i \colon \Sigma \to \Sigma'$
such that $f=i_*$, then $f$ preserves the intersection form. 

\begin{example}
Let $\mathcal{O}$ be
any orbifold with two cone points of orders $10$ and $15$. The 
torsion part of the orbifold homology is isomorphic to
the quotient of the additive group $\Z_{10} \oplus \Z_{15}$ 
by the sum of the two obvious generators. This group is isomorphic
to $\Z_5$. Take an isomorphism of $\Z_5$, which squares each element. 
This isomorphism cannot be realized by an orbifold diffeomorphism. This is the case
because an orbifold diffeomorphism has to be multiplication by $1$ or $-1$ 
on each of $\Z_{10}$ and $\Z_{15}$, so it has to be multiplication
by $1$ or $-1$ on the quotient $\Z_5$. 
\end{example}

Not every
automorphism of the first orbifold 
homology group preserves the order
of the orbifold singularities.

\begin{example} \label{assumptioniii)example1}
Let $\mathcal{O}$ be a compact, connected, orientable smooth orbisurface
with underlying topological space equal to a $2$\--dimensional torus. Suppose
that $\mathcal{O}$ has precisely two cone points $p_1, \, p_2$ of respective orders $o_1=5$
and $o_2=10$.
Then for any point $x_0 \in \mathcal{O}$
$$
\pi_1^{\textup{orb}}(\mathcal{O},\,x_0)=\langle \alpha, \, \beta, \, \gamma_1, \, \gamma_2 \, | \,
\gamma_1 \, \gamma_2=[\alpha, \, \beta], \, \, \gamma_1^5=\gamma_2^{10}=1 \rangle,
$$
and
$$
\op{H}_1^{\textup{orb}}(\mathcal{O},\, \Z)=
\langle \alpha, \, \beta, \, \gamma_1, \, \gamma_2 \, | \,
\gamma_1 + \gamma_2=0, \, 5 \, \gamma_1= 10 \, \, \gamma_2=0 \rangle.
$$
The assignment
\begin{eqnarray} \label{F}
F \colon \alpha \to \alpha, \, \beta \to \beta, \, \gamma_1 \mapsto \gamma_2,  \, \gamma_2 \mapsto \gamma_1
\end{eqnarray}
defines a group automorphism of  
$\op{H}_1^{\textup{orb}}(\mathcal{O},\, \Z)$. If this automorphism
is induced by a diffeomorphism, then the same map on 
$\pi_1^{\textup{orb}}(\mathcal{O},\, x_0)$ should be an isomorphism,
which is false since in $\pi_1^{\textup{orb}}(\mathcal{O},\, x_0)$ the classes
$\gamma_1$ and $\gamma_2$ have different orders.
Thus the assignment (\ref{F}) is not geometric. 
 \end{example}

\subsubsection{Symplectic and torsion geometric maps}

Next we introduce the notion of symplectic isomorphism as well as that
of singularity--order preserving isomorphism.

Let $K, \, L$ be arbitrary groups, and let $h \colon K \to L$ be a group isomorphism. 
We denote by $K_{\op{T}}, \, L_{\op{T}}$ the corresponding torsion subgroups,
by $K_{\textup{F}}, \, L_{\textup{F}}$ the quotients $K/K_{\op{T}}$, $L/L_{\op{T}}$,
by $h_{\op{T}}$ the restriction of $h$ to $K_{\op{T}} \to L_{\op{T}}$, and
by $h_{\textup{F}}$ the restriction of $h$ to $K_{\textup{F}} \to L_{\textup{F}}$.

\begin{definition} \label{preiso}
Let $\mathcal{O}, \, \mathcal{O}'$ be compact,
connected, orientable smooth orbisurfaces. We say that an isomorphism
$Z=\op{H}_1^{\textup{orb}}(\mathcal{O},\,\Z) \to Z'=\op{H}_1^{\textup{orb}}(\mathcal{O}',\,\Z)$
is \emph{torsion geometric} if and only if
the isomorphism $Z _{\op{T}} \to Z'_{\op{T}}$ sends a geometric torsion basis to a geometric torsion basis preserving the 
order of the orbifold singularities.
\end{definition}

\begin{definition} \label{symiso}
Let $\mathcal{O}, \, \mathcal{O}'$ be compact,
connected, orientable smooth orbisurfaces. We say that an isomorphism
$Z=\op{H}_1^{\textup{orb}}(\mathcal{O},\,\Z) \to Z'=\op{H}_1^{\textup{orb}}(\mathcal{O}',\,\Z)$
is \emph{symplectic} if the isomorphism $Z _{\textup{F}} \to Z'_{\textup{F}}$ 
respects the symplectic form, c.f. (\ref{if2}), i.e. the matrix of the isomorphism w.r.t. symplectic bases is in the integer symplectic group.
\end{definition}

\begin{definition}
Let $Z,\,Z'$ be respectively the first integral orbifold homology groups of
compact connected orbisurfaces $\mathcal{O},\,\mathcal{O}'$ of the
same Fuchsian signature.
We define the  following set of isomorphisms 
\begin{eqnarray} \label{Sn}
\textup{S}^{Z,Z'}:=\{ f \in \textup{Iso}(Z, \, Z') \, | \,  f \textup{ is torsion geometric}  \},
\end{eqnarray}
and we denote by
\begin{eqnarray} \label{Sp}
\textup{Sp}(Z,\,Z'):=\{h \in \textup{Iso}(Z, \, Z')    \, | \,  h \textup{ is  symplectic}\}.
 \end{eqnarray}
 \end{definition}

\begin{definition} \label{equivalent}
Let $Z_1,\,Z_2$ be respectively the first integral orbifold homology groups of
compact connected orbisurfaces $\mathcal{O}_1,\,\mathcal{O}_2$ of the
same Fuchsian signature. Let
$f_i \colon Z_i \to T$
be homomorphisms into a torus $T$. We say that $f_1$ is 
\emph{$\textup{Sp}(Z_1,\,Z_2)   \cap \textup{S}^{Z_1,Z_2}$\--equivalent to $f_2$} if there
exists an isomorphism $i \colon Z_1 \to Z_2$ such that there is an identity of maps
$f_2=f_1 \circ i$ and $i \in \textup{Sp}(Z_1,Z_2)   \cap \textup{S}^{Z_1,Z_2}$.
\end{definition}

\subsubsection{Geometric isomorphisms: characterization}

The main result of this section is Theorem \ref{veryhard},
where we give a characterization 
of geometric isomorphisms in terms of symplectic maps, c.f. Definition \ref{symiso},
and torsion geometric maps, c.f. Definition \ref{preiso}.  
Then we deduce from it Proposition \ref{conjecture0} which is a key ingredient of the proof 
of the classification theorem.

\begin{lemma} \label{preprop}
Let $\Sigma_i$, $i=1, \, 2$, be compact, connected, oriented smooth orbisurfaces
with the same Fuchsian signature $(g; \, \vec{o})$, and suppose that
$G$ is an isomorphism from the first integral orbifold homology group 
$\op{H}_1^{\textup{orb}}(\Sigma_1, \, \Z)$ onto
$\op{H}_1^{\textup{orb}}(\Sigma_2, \, \Z)$ which is symplectic and torsion geometric, c.f. Definition \ref{symiso}
and Definition \ref{preiso}. Then there
exists an orbifold diffeomorphism $g \colon \Sigma_1 \to \Sigma_2$
such that $G_{\op{T}}=(g_*)_{\op{T}}$ and $G_{\textup{F}}=(g_*)_{\textup{F}}$.
\end{lemma}

\begin{proof}
Because an orbifold is classified up to orbifold diffeomorphisms by its Fuchsian 
signature, c.f. Theorem 
\ref{2dorbifoldtheorem}, without loss of generality we may assume that 
$\Sigma=\Sigma_1=\Sigma_2$,
and let  $\{\gamma_k\}$ be a geometric torsion basis of 
$\op{H}_1^{\textup{orb}}(\Sigma, \, \Z)$, c.f Definition \ref{geometrictorsionbasis}.
Choose an orbifold atlas for $\Sigma$ such that the orbifold chart $U_k$ around the
$k^{\scriptop{th}}$ singular point $p_k$ is homeomorphic to a disk $D_k$
modulo a finite group of diffeomorphisms, and such that every singular
point is contained in precisely one chart.
The oriented boundary loop 
$\partial U_k$ represents
the class $\gamma_k \in \op{H}_1^{\textup{orb}}(\Sigma, \, \Z)_{\op{T}}$.
This in particular implies that there exists an orbifold diffeomorphism
$f_{\tau} \colon  U:=\bigcup_{k} U_k \to U$ 
such that $f_{\tau}(p_k)=p_{\tau(k)}$ and 
\begin{eqnarray} \label{ftu}
f_{\tau}(U_k)=U_{\tau(k)}.
\end{eqnarray}
Replace each orbifold chart around a singular point of $\Sigma$ by a manifold chart. This
gives rise to a manifold atlas, which defines a compact, connected, smooth orientable
surface
$\widehat{\Sigma}$ without boundary\footnote{recall that we are always assuming, unless otherwise
specified, that all manifolds
and orbifolds in this paper have no boundary.}
with the same underlying space as that of the orbisurface $\Sigma$.
Let
$\op{H}_1^{\textup{orb}}(\Sigma, \, \Z)_{\textup{F}} \to \op{H}_1(\widehat{\Sigma}, \,
\Z)$ be the natural intersection form preserving automorphism in (\ref{if3}). The automorphism obtained from $G_{\textup{F}}$ by conjugation with the 
aforementioned automorphism preserves the intersection 
form on  $\op{H}_1(\widehat{\Sigma}, \, \Z)$.
By Lemma \ref{triviallemma} there is a surface diffeomorphism 
which induces it, which sends
$U_k$ to $U_{\tau(k)}$.
Because each
diffeomorphism of a circle is isotopic to a rotation or a reflection, by (\ref{ftu})
the restriction of  the aforementioned surface diffeomorphism
to the punctured surface  $\widehat{\Sigma} \setminus \bigcup_{k} U_k$
may be glued to $f_{\tau}$, 
along the boundary circles $\partial U_k$,  to give rise
to an orbifold diffeomorphism of $\Sigma$ which satisfies the required properties.
\end{proof}

\begin{theorem} \label{veryhard}
Let $\Sigma_1, \, \Sigma_2$ be compact,
connected, orientable smooth orbisurfaces of the same Fuchsian signature. An isomorphism
$Z_1=\op{H}_1^{\textup{orb}}(\Sigma_1,\,\Z) \to Z_2=\op{H}_1^{\textup{orb}}(\Sigma_2,\,\Z)$
is orbisurface geometric if and only if it is symplectic and torsion geometric.
\end{theorem}

\begin{proof}
Suppose that $G \colon Z_1 \to Z_2$ is an orbisurface geometric isomorphism, c.f.
Definition \ref{geomiso}.
It follows from the definition of orbifold diffeomorphism that if the group isomorphism
$G$ is induced by an orbifold diffeomorphism, then 
the induced map $G_{\textup{F}}$ on the free quotient preserves the intersection form 
and $G$ sends geometric torsion basis to geometric torsion basis of the orbifold homology preserving the order of the orbifold singularities.

Conversely, suppose that $G$ is symplectic and torsion geometric, c.f. Definition \ref{symiso}, Definition \ref{preiso}.
By Theorem \ref{2dorbifoldtheorem} we may
assume without loss of generality that $\Sigma=\Sigma_1=\Sigma_2$, so
$G$ is an automorphism of $\op{H}_1^{\textup{orb}}(\Sigma,\, \Z)$.
By Lemma \ref{preprop} there exists an orbifold diffeomorphism $g \colon \Sigma \to \Sigma$
such that 
\begin{eqnarray} \label{identities}
(g_*)_{\textup{F}}=G_{\textup{F}}, \, \,  \, \, \, \, (g_*)_{\op{T}}=G_{\op{T}}.
\end{eqnarray}
Let us define the mapping
\begin{eqnarray} \label{K}
K:=g_* \circ G^{-1} \colon \op{H}_1^{\textup{orb}}(\Sigma, \, \Z) \to 
\op{H}_1^{\textup{orb}}(\Sigma, \, \Z).
\end{eqnarray}
$K$ given by (\ref{K}) is a group isomorphism because it is the composite of two group isomorphisms. Because
of the identities in expression (\ref{identities}),  $K$ satisfies that
\begin{eqnarray} \label{equations2}
K_{\op{T}}=\textup{Id}_{\op{H}_1^{\textup{orb}}(\Sigma, \, \Z)_{\op{T}}}, \, \, \, \, 
K_{\textup{F}}=\textup{Id}_{\op{H}_1^{\textup{orb}}(\Sigma, \, \Z)_{\textup{F}}}.
\end{eqnarray}
If the genus of the underlying surface $|\Sigma|$ is $0$, then $\op{H}_1^{\textup{orb}}(\Sigma, \, \Z)_{\textup{F}}$ is trivial and $\op{H}_1^{\textup{orb}}(\Sigma, \, \Z)_{\op{T}}=
\op{H}_1^{\textup{orb}}(\Sigma, \, \Z)$, $K=K_{\op{T}}=
\textup{Id}$,
so $g_*=G$, and we are done. If the genus of $|\Sigma|$ is strictly positive, then
there are two cases. First of all, if $K$ is the identity map, then $g_*=G$, and we are done. 
If otherwise $K$ is not the identity map, then
$g_* \neq G$, and let us choose a free subgroup $F$ of $\op{H}_1^{\textup{orb}}(\Sigma, \, \Z)$
such that $\op{H}_1^{\textup{orb}}(\Sigma, \, \Z)=F \oplus \op{H}_1^{\textup{orb}}(\Sigma, \, \Z)_{\op{T}}$, and a symplectic basis $\alpha_i, \, \beta_i$ of the free group
$F$. To make the forthcoming notation simpler rename $e_{2i-1}=\alpha_i$ and
$e_{2i}=\beta_i$, for all $i$ such that $1 \le i \le g$.
By the right hand side of equation (\ref{equations2}), we have that
there exist non--negative integers $a^i_k$ 
such that if $n$ is the number singular of points of $\Sigma$, then
\begin{eqnarray} \label{Kformula}
K(e_i)=e_i+ \sum_{k=1}^n  a^i_k \, \gamma_k.
\end{eqnarray}
By the left hand side of (\ref{equations2}),
\begin{eqnarray} \label{KT}
K(\gamma_k)=\gamma_k. 
\end{eqnarray}
Next we define a new isomorphism, which we call $K_{\scriptop{new}}$,  
by altering $K$ in the following fashion.

Choose an orbifold atlas for the orbisurface $\Sigma$ such that for each singular
point $p_j$ of $\Sigma$ there is a unique orbifold chart which contains it, and which
is homeomorphic to a disk $D_j$ modulo a finite group of diffeomorphisms. Let 
$\widehat{\Sigma}$ be the smooth surface whose underlying topological space
is the underlying surface $|\Sigma|$ to $\Sigma$, and whose smooth structure
is given by the atlas defined by replacing each chart which contains a singular point
by a manifold chart from the corresponding disk. 
Let $\widetilde{e}_i$ be a  loop
which represents the homology class $e_i \in \op{H}^{\scriptop{orb}}_1(\Sigma, \,\Z)$.
Take an annulus $\mathcal{A}$ 
in $\Sigma$ such that the loop $\widetilde{e}_i$ crosses $\mathcal{A}$ exactly once, in the sense
that it intersects the boundary $\partial \mathcal{A}$  of $\mathcal{A}$ twice, 
once at each of the two boundary components of $\partial \mathcal{A}$, and such that it contains 
the $j^{\scriptop{th}}$ singular point $p_j$
of $\Sigma$ which is enclosed by the oriented loop whose homology
class is $\gamma_j$, and no other singular point of $\Sigma$.
This can be done by choosing the annulus $\mathcal{A}$ in such a way that its boundary curves
represent the same class as the dual element $e_{i+1}$ to $e_i$
(instead of $e_{i+1}$ we may have $e_{i-1}$, depending on how the symplectic basis of
$e_i$ is arranged). 
Equip the topological space $|\mathcal{A}| \subset |\Sigma| =|\widehat{\Sigma}|$ with the
smooth structure which comes from restricting to $|\mathcal{A}|$ the charts of the smooth structure of
$\widehat{\Sigma}$, and let $\widehat{\mathcal{A}}$ be the smooth submanifold--with--boundary of 
$\widehat{\Sigma}$ which arises in such way.
The loop $\widetilde{e}_i$ intersects the boundary of the annulus
at an initial point $x_{ij}$ and at an end point $y_{ij}$, and intersects the annulus itself at a curved
segment path $[x_{ij},\, y_{ij}]$, which by possibly choosing
a different representative of the class $e_i$, may be assumed to be a smooth embedded
$1$\--dimensional submanifold--with--boundary of $\widehat{\mathcal{A}}$.
Replace the segment $[x_{ij},\, y_{ij}]$ of the loop $\widetilde{e}_i$ by a smoothly embedded
path $P(x_{ij},\, y_{ij})$ which starts at the point $x_{ij}$, goes towards the cone point enclosed by 
$\gamma_j$ while inside
of the aforementioned annulus, and goes around 
it precisely once, to finally come back to end up at the end point $y_{ij}$.
The replacement of the segment $[x_{ij},\, y_{ij}]$ by the path $P(x_{ij},\, y_{ij})$ gives rise to a new loop $\widetilde{e}'_i$, which agrees  with the loop $\widetilde{e}_i$ outside of $\mathcal{A}$.

We claim that there exists an orbifold 
diffeomorphism $\widehat{h}_{ij} \colon \widehat{\Sigma} \to \widehat{\Sigma}$ which 
is the identity map outside of the annulus and which
sends the segment $[x_{ij},\, y_{ij}]$ to the path $P(x_{ij},\, y_{ij})$,
and hence the loop $\widetilde{e}_i$ to the loop $\widetilde{e}'_i$. 
Indeed,
let  $\widehat{f}$  be a diffeomorphism of the aforementioned annulus $\widehat{\mathcal{A}}$
with 
$$
\widehat{f}([x_{ij}, \, y_{ij}]) = P(x_{ij},\, y_{ij}),
$$
 and with  $\widehat{f}$  being the 
identity on $\partial \widehat{\mathcal{A}}$.  Then let  $\widehat{k} \colon \widehat{\mathcal{A}}
\to \widehat{\mathcal{A}}$  
 be a diffeomorphism which is  the identity on  $P(x_{ij},\,y_{ij})$  and 
on $\partial \widehat{\mathcal{A}}$,  such that  $\widehat{k}(\widehat{f}(p_j))=p_j$.  Such a diffeomorphism 
$\widehat{k}$  exists since 
cutting  the annulus  along  $P(x_{ij},\,y_{ij})$  
 gives a disk, and in a disk there is a diffeomorphism taking any  
interior point  $x$  to any other interior point  $y$  and fixing a  neighborhood of the 
boundary. The composition 
$
\widehat{h}'_{ij}:=\widehat{k} \circ \widehat{f} \colon \widehat{\mathcal{A}} \to \widehat{\mathcal{A}}
$
is a diffeomorphism which is the identity on $\partial \widehat{\mathcal{A}}$, 
which sends the segment $[x_{ij},\, y_{ij}]$ to the path $P(x_{ij},\, y_{ij})$ and
such that
\begin{eqnarray} \label{hijp}
\widehat{h}'_{ij}(p_j)=p_j.
\end{eqnarray}
Because the mapping $\widehat{h}'_{ij}$ is the identity along $\partial \widehat{\mathcal{A}}$, 
it extends to a diffeomorphism $\widehat{h}_{ij}$ along such boundary, 
which satisfies the required properties.
Since $p_j$ is contained in a unique orbifold chart and (\ref{hijp}) holds,
the way in which we defined the smooth structure of $\widehat{\Sigma}$ from
the orbifold structure of $\Sigma$ gives that $\widehat{h}_{ij}$ defines an orbifold diffeomorphism
$\Sigma \to \Sigma$. To emphasize that $\widehat{h}_{ij}$ is a diffeomorphism at the level of
orbifolds, we denote it by $h_{ij} \colon \Sigma \to \Sigma$.

The isomorphism $h^*_{ij}$ induced on the orbifold homology 
by the orbifold diffeomorphism $h_{ij}$ is given by
\begin{eqnarray} \label{for1}
h^*_{ij}(e_k )=e_k, \, \, k \neq i, \, \, \, \, \,  h^*_{ij}(e_i)=e_i+\gamma_j,
\end{eqnarray}
\begin{eqnarray} \label{for2}
 h^*_{ij}(\gamma_k)=\gamma_k,
\end{eqnarray}
where $1 \le i \le 2g$ and $1 \le k \le n$, and notice that in the first equality we have used
that the boundary curve of the annulus  is in the class of $e_{i+1}$.
Define $f_{ij} \colon \Sigma \to \Sigma$ to be
the orbifold diffeomorphism obtained by composing $h_{ij}$ with itself
precisely $a^j_i$ times. It follows from (\ref{for1}) and (\ref{for2}) that
\begin{eqnarray} \label{f1}
f^*_{ij}(e_k)=e_k \, \, \textup{ if } k \neq i, \, \, \, f^*_{ij}(e_i)=e_i+a_i^j\, \gamma_j, 
\end{eqnarray}
and
\begin{eqnarray} \label{f2}
 f^*_{ij}(\gamma_k)=\gamma_k,
\end{eqnarray}
where $1 \le i \le 2g$ and $1 \le k \le n$.

The isomorphisms $f^*_{ij}$ commute with each other, because they are the identity on the torsion subgroup,
and only change one loop of the free part which does not affect the other loops\footnote{Observe that this is completely false
in the fundamental group.}. Therefore combining expressions (\ref{Kformula}), (\ref{f1}) and (\ref{f2})
we arrive at the identity
\begin{eqnarray} \label{fij}
((f_{ij}^*)^{-1} \circ K)(e_i)=e_i+ \sum_{k=1, \, k \neq j}^n a^i_k \, \gamma_k.
\end{eqnarray}
On the other hand, it follows from (\ref{KT}) and (\ref{f2}) that
\begin{eqnarray}
((f_{ij}^*)^{-1} \circ K)(\gamma_k)=\gamma_k. \label{KKT}
\end{eqnarray}
Now consider the isomorphism $K_{\scriptop{new}}$ of the orbifold homology group
defined by
\begin{eqnarray} \label{alter}
K_{\scriptop{new}}:= (\bigcirc_{1 \le i \le 2g, \, 1 \le j \le n} \, (f^*_{ij})^{-1}) \circ K,
\end{eqnarray}
where recall that in (\ref{alter}), $g$ is the genus of the surface underlying $\Sigma$, and $n$ is the number of singular points. 
It follows from expression (\ref{KKT}), and from (\ref{fij}), by induction on $i$ and $j$, that
$K_{\scriptop{new}}$ given by (\ref{alter}) is the identity map on the orbifold homology, which
then by formula (\ref{K}) implies that 
\begin{eqnarray} 
G=g_* \circ K^{-1} &=&  g_{*} \circ  (\bigcirc_{1 \le i \le 2g, \, 1 \le j \le n}\,  f^*_{ij}) \nonumber \\
                                 &=& (g \circ  (\bigcirc_{1 \le i \le 2g, \, 1 \le j \le n}\,  f_{ij}))_*, \nonumber
\end{eqnarray}
and hence $G$ is a geometric isomorphism induced by 
$$
g \circ  (\bigcirc_{1 \le i \le 2g, \, 1 \le j \le n}\,  f_{ij}),
$$
which is a composite of orbifold diffeomorphisms, and hence an orbifold diffeomorhism itself.
\end{proof}

Recall that the mappings $\mu_{\textup{h}},\,
\mu_{\textup{h}}'$ are, respectively, the homomorphisms induced on homology
by the monodromy homomorphisms $\mu,\,\mu'$ of the connections $\Omega,\, \Omega'$
of symplectically orthogonal complements to the tangent spaces to the $T$\--orbits in
$M,\,M'$, respectively (c.f. Proposition \ref{easynow}). Recall that $\nu$,
$\nu' $ are the unique $2$\--forms respectively on $M/T$ and $M'/T$
such that $\pi^*\nu|_{\Omega_x}=\sigma|_{\Omega_x}$ and 
$\pi'^*\nu'|_{\Omega'_{x'}}=\sigma'|_{\Omega'_{x'}}$, for every $x \in M$, $x' \in M'$.

\begin{prop} \label{conjecture0}
Let $(M, \,\sigma), \, (M', \, \sigma')$ be compact, connected $2n$\--dimensional symplectic
manifolds equipped with an effective  symplectic action of a $(2n-2)$\--dimensional torus $T$ for which 
at least one, and hence every $T$\--orbit is a $(2n-2)$\--dimensional
symplectic submanifold of $(M, \, \sigma)$ and $(M', \, \sigma')$, respectively. Let $K=\op{H}_1^{\textup{orb}}(M/T, \, \Z)$ and similarly $K'$. 
Suppose that the symplectic orbit spaces $(M/T, \, \nu)$ and $(M'/T, \, \nu')$ are orbifold symplectomorphic and that $\mu_{\textup{h}}$ is 
$\textup{Sp}(K',\,K)   \cap \textup{S}^{K',K}$\--equivalent 
to $\mu'_{\textup{h}}$ via an automorphism $G$ from the first integral orbifold homology group $\op{H}_1^{\textup{orb}}(M'/T, \, \Z)$
onto $\op{H}_1^{\textup{orb}}(M/T, \, \Z)$ (c.f. Corollary \ref{3.7} and 
Lemma \ref{orbifoldsymplectomorphic}
for the definitions of $\nu,\,\nu'$, or see the preceding paragraphs). Then there exists an orbifold
diffeomorphism $g \colon M'/T \to M/T$ such that $G=g_*$ and $\mu'_{\textup{h}}=\mu_{\textup{h}} \circ g_*$.
\end{prop}
\begin{proof}
It follows from Theorem \ref{veryhard} applied to the groups $Z', \, Z$, which respectively are the
first integral orbifold homology group of $\Sigma'$, and of $\Sigma$,
where  $\Sigma'=M'/T$, $\Sigma=M/T$.
\end{proof}

\subsubsection{Fuchsian signature space}

We define the invariant of $(M, \, \sigma)$
which encodes the monodromy of the connection
for $\pi \colon M \to M/T$ of symplectic orthogonal complements
to the tangent spaces to the $T$\--orbits, c.f. Definition \ref{nfmono}.

\begin{definition} \label{signature}
Let $\mathcal{O}$ be a smooth orbisurface with $n$ cone points $p_k$, $1 \le k \le n$.
The \emph{Fuchsian signature $\textup{sig}(\mathcal{O})$ 
of $\mathcal{O}$} is the $(n+1)$\--tuple $(g; \, \vec{o})$
where $g$ is the  genus of the underlying surface to the orbisurface $\mathcal{O}$,
$o_k$ is the order of the point $p_k$, which we require to be strictly positive, and $\vec{o}=(o_k)_{k=1}^n$, where
$o_k \le o_{k+1}$, for all $1 \le k \le n-1$.
\end{definition}

\begin{definition} \label{MSn}
Let $\vec{o}$ be an $n$\--dimensional tuple of strictly positive integers. We define
$$
\mathcal{M}\textup{S}^{\vec{o}}_n:=\{B \in \textup{GL}(n, \, \Z) \, | \, B \cdot \vec{o}=\vec{o}\},
$$
i.e. $\mathcal{M}\textup{S}^{\vec{o}}_n$ is
the group of $n$\--dimensional matrices which permute elements preserving the tuple of orders $\vec{o}$.
\end{definition}

\begin{definition} \label{key2}
Let $(g; \, \vec{o})$ be an $(n+1)$\--tuple of integers, where the $o_k$'s are strictly positive
and non--decreasingly ordered. Let $T$ be a torus.
The \emph{Fuchsian signature space associated to $(g;\, \vec{o})$} is the quotient space
\begin{eqnarray}
T^{2g+n}_{(g;\, \vec{o})}/\mathcal{G}_{(g, \, \vec{o})}
\end{eqnarray}
 where
$$
T^{2g+n}_{(g;\, \vec{o})}:=\{(t_i)_{i=1}^{2g+n} \in T^{2g+n}  \, | \, \textup{the order of } t_i \textup{ is a multiple of } o_i, \, \, 2g +1\le i \le n \}.
$$
and where
 $\mathcal{G}_{(g, \, \vec{o})}$ is the group of matrices
\begin{eqnarray} \label{matrixgroup}
\mathcal{G}_{(g; \, \vec{o})}:=\{  
\left( \begin{array}{cc}
A & \textup{0}  \\
C & D 
\end{array} \right) \in \textup{GL}(2g+n, \, \Z)
 \, | \, A \in \textup{Sp}(2g, \, \Z), \, \,  D \in \mathcal{M}\textup{S}^{\vec{o}}_n \},
 \end{eqnarray}
where $\textup{Sp}(2g, \, \Z)$ is the group of $2g$\--dimensional symplectic matrices with integer entries,
c.f. Section 3.1.1. and Definition \ref{Fkey}, and
$\mathcal{M}\textup{S}^{\vec{o}}_n$ is the group of $n$\--dimensional matrices which permute elements preserving
the tuple of orders $\vec{o}$, c.f. Definition \ref{MSn}.
\end{definition}

\begin{remark}
The lower left block $C$ in the definition of $\mathcal{G}_{(g; \, \vec{o})}$ in Definition \ref{key2}
of the description of item 4) is allowed to be any matrix.
Intuitively, this reflects that there are many free subgroups of
the first integral orbifold homology group with together with the torsion span the entire group. 
\end{remark}

\begin{definition} \label{nfmono}
Let $(M, \, \sigma)$ be a smooth compact and connected symplectic manifold of dimension
$2n$ and let $T$ be a $(2n-2)$\--dimensional torus which acts effectively  on $(M, \, \sigma)$ by means of
symplectomorphisms. We furthermore assume that at least one, and hence every
$T$\--orbit is $(2n-2)$\--dimensional symplectic
submanifold of $(M, \, \sigma)$. Let $(g; \, \vec{o}) \in \Z^{1+m}$ be the Fuchsian signature
of the orbit space $M/T$.
Let $\{\gamma_k\}_{k=1}^m$ be a geometric torsion basis, 
c.f. Definition \ref{geometrictorsionbasis}.
Let $\{\alpha_i, \, \beta_i\}_{i=1}^{g}$ be a symplectic basis of a free subgroup 
of our choice of $\op{H}_1^{\textup{orb}}(M/T, \, \Z)$, c.f. expression (\ref{Fcap}) and  \cite[Th.\,2.3]{MS} whose direct sum with the torsion subgroup is equal to  $\op{H}_1^{\textup{orb}}(M/T, \, \Z)$.
Let $\mu_{\textup{h}}$ be the homomomorphism induced on homology by
the monodromy homomorphism $\mu$ associated to the connection $\Omega$, c.f. 
Proposition \ref{easynow}.
The \emph{monodromy invariant of $(M, \sigma, \, T)$} is the  
$\mathcal{G}_{(g; \, \vec{o})}$\--orbit
\begin{eqnarray}
 \mathcal{G}_{(g, \, \vec{o})} \cdot ((\mu_{\textup{h}}(\alpha_i), \, \mu_{\textup{h}}(\beta_i))_{i=1}^{g}, \, (\mu_{\textup{h}}(\gamma_k))_{k=1}^m),
\end{eqnarray}
of the $(2g+n)$\--tuple $((\mu_{\textup{h}}(\alpha_i), \, \mu_{\textup{h}}(\beta_i))_{i=1}^{g}, \, (\mu_{\textup{h}}(\gamma_k))_{k=1}^n)$,
where $\mathcal{G}_{(g, \, \vec{o})}$ is the group of matrices
given in (\ref{matrixgroup}).
\end{definition}

Because the invariant in Definition \ref{nfmono} depends on choices, it is unclear whether it is well defined.

\begin{lemma} \label{welldefined}
Let $(M, \, \sigma)$ be a smooth compact and connected symplectic manifold of dimension
$2n$ and let $T$ be a $(2n-2)$\--dimensional torus which acts effectively on $(M, \, \sigma)$ by means of
symplectomorphisms. We furthermore assume that at least one, and hence every
$T$\--orbit is $(2n-2)$\--dimensional symplectic
submanifold of $(M, \, \sigma)$.
The monodromy invariant of $(M, \, \sigma, \, T)$ is well defined, in the following sense.

Suppose that the signature of $M/T$ is $(g; \, \vec{o}) \in \Z^{1+m}$.
Let $F$, $F'$ be any two free subgroups of $\op{H}_1^{\textup{orb}}(M/T, \, \Z)$
whose direct sum with the torsion subgroup is equal to  $\op{H}_1^{\textup{orb}}(M/T, \, \Z)$, and let
$\alpha_i, \, \beta_i$ and $\alpha'_i, \, \beta'_i$ be symplectic bases of $F, \, F'$, respectively.
Let  $\tau \in \textup{S}^{\vec{o}}_m$.
Let $\mu_{\textup{h}}$ be the homomorphism induced in homology by means of the Hurewicz map
from the monodromy homomorphism $\mu$ of the connection of symplectic orthogonal complements
to the tangent spaces to the $T$\--orbits, c.f. Proposition \ref{easynow}. Then
the tuples $((\mu_{\textup{h}}(\alpha_i), \, \mu_{\textup{h}}(\beta_i))_{i=1}^{g}, \, (\mu_{\textup{h}}(\gamma_k))_{k=1}^m)$
and  $((\mu_{\textup{h}}(\alpha'_i), \, \mu_{\textup{h}}(\beta'_i))_{i=1}^{g}, \, (\mu_{\textup{h}}(\gamma_{\tau(k)})_{k=1}^m)$
lie in the same $\mathcal{G}_{(g; \, \vec{o})}$\--orbit in the Fuchsian signature space.
\end{lemma}

\begin{proof}
A symplectic basis of the maximal free subgroup $F$ of $\op{H}_1^{\textup{orb}}(M/T, \, \Z)$
may be taken to a symplectic basis of the maximal free subgroup $F'$ of $\op{H}_1^{\textup{orb}}(M/T, \, \Z)$
by a matrix 
$$
X:=\left(\begin{array}{cc}
A & \textup{0}  \\
C & \textup{Id} 
\end{array} \right) \in  \textup{GL}(2g+m, \, \Z),
$$
where the upper block $A$ is a $2g$\--dimensional matrix in the integer symplectic linear group $\textup{Sp}(2g, \, \Z)$,
and the lower block $C$ is $(m \times 2g)$\--dimensional matrix with integer entries. Here $\textup{Id}$
denotes the $n$\--dimensional identity matrix, and $\textup{0}$ is 
the $(m \times 2g)$\--dimensional matrix
all the entries of which equal $0$.
A geometric torsion basis
can be taken to another geometric torsion basis by preserving the order $\vec{o}$
of the orbifold singularities by a matrix of the form
$$
Y:=\left( \begin{array}{cc}
\textup{Id} & \textup{0}  \\
\textup{0} & N 
\end{array} \right)  \in \textup{GL}(2g+m, \, \Z),
$$
for a certain matrix $N \in \mathcal{M}\textup{S}_m^{\vec{o}}$, and
the product matrix 
$$
X \, Y=\left( \begin{array}{cc}
A & \textup{0}  \\
C & N 
\end{array} \right) 
$$
lies in $\mathcal{G}_{(g; \, \vec{o})}$.
\end{proof}

With this matrix terminology, we can restate Proposition \ref{conjecture0} in the following terms.

See the paragraph preceding Proposition \ref{conjecture0} for a reminder of the terminology which we use next.

\begin{prop} \label{conjecture01}
Let $(M, \,\sigma), \, (M', \, \sigma')$ be compact, connected $2n$\--dimensional symplectic
manifolds equipped with an effective symplectic action of a $(2n-2)$\--dimensional torus $T$ for which 
at least one, and hence every $T$\--orbit is a $(2n-2)$\--dimensional
symplectic submanifold of $(M, \, \sigma)$ and $(M', \, \sigma')$, respectively. 
Suppose that the symplectic orbit spaces $(M/T, \, \nu)$ and $(M'/T, \, \nu')$ have Fuchsian signature $(g; \, \vec{o})$, that they are orbifold symplectomorphic and that 
the $\mathcal{G}_{(g;\, \vec{o})}$\--orbits of the $(2g+m)$\--tuples 
$((\mu_{\textup{h}}(\alpha'_i), \, \mu_{\textup{h}}(\beta'_i))_{i=1}^{g}, \, (\mu_{\textup{h}}(\gamma'_k))_{k=1}^m)$
and $((\mu'_{\textup{h}}(\alpha'_i), \, \mu'_{\textup{h}}(\beta'_i))_{i=1}^{g}, \, 
(\mu'_{\textup{h}}(\gamma'_k))_{k=1}^m)$ are equal, where
$\alpha_i, \, \beta_i$ and $\alpha'_i, \, \beta'_i$ are respectively symplectic bases of free 
homology subgroups of $\op{H}_1^{\textup{orb}}(M/T, \, \Z)$
and $\op{H}_1^{\textup{orb}}(M'/T, \, \Z)$ which together with the corresponding torsion subgroups span the entire group,
and $\gamma_k, \, \gamma'_k$ are corresponding geometric torsion bases.
Then there exists an orbifold
diffeomorphism $g \colon M'/T \to M/T$ such that $\mu'_{\textup{h}}=\mu_{\textup{h}} \circ g_*$.
\end{prop}

\begin{proof}
By assumption the 
$\mathcal{G}_{(g;\, \vec{o})}$\--orbits of the
$(2g+m)$\--tuples $((\mu_{\textup{h}}(\alpha'_i), \, \mu_{\textup{h}}(\beta'_i))_{i=1}^{g}, \, (\mu_{\textup{h}}(\gamma'_k))_{k=1}^m)$
and $((\mu'_{\textup{h}}(\alpha'_i), \, \mu'_{\textup{h}}(\beta'_i))_{i=1}^{g}, \, (\mu'_{\textup{h}}(\gamma_k))_{k=1}^m)$ are equal, and hence
$\mu'_{\textup{h}}=\mu_{\textup{h}} \circ G$, where $G$ is the isomorphism defined by $G(\alpha_i)=\alpha'_i$,
$G(\beta_i)=\beta'_i$ and $G(\gamma_i)=\gamma'_{\tau(i)}$. By its definition $G$ is symplectic and torsion geometric, i.e.
$G \in \textup{Sp}(K',\,K) \cap S^{K',\,K}$, where $K:=\op{H}_1^{\textup{orb}}(M/T, \, \Z)$
and $K':=\op{H}_1^{\textup{orb}}(M'/T, \, \Z)$.
Now the result follows from Proposition \ref{conjecture0}.
\end{proof}

\subsection{Uniqueness}

\subsubsection{List of ingredients of $(M,\,\sigma,\,T)$}

We start by assigning a list 
of invariants to $(M, \, \sigma)$.

\begin{definition} \label{def2}
Let $(M, \,\sigma)$ be a compact connected $2n$\--dimensional symplectic
manifold equipped
with an effective symplectic action of a $(2n-2)$\--dimensional torus $T$ for which at least
one, and hence every $T$\--orbit is a $(2n-2)$\--dimensional
symplectic submanifold of $M$. The \emph{list of ingredients of $(M,\,\sigma,\,T)$}
consists of the following items. 
\begin{itemize}
\item[1)]
The Fuchsian signature  $(g; \, \vec{o}) \in \Z^{1+m}$ 
of the orbisurface $M/T$ (c.f. Remark \ref{manifoldcase} and
Definition \ref{signature}).
\item[2)]
The total symplectic area of the symplectic orbisurface $(M/T,\, \nu)$, where 
the symplectic form $\nu$ is defined
by the condition 
$\pi^* \nu|_{\Omega_x}=\sigma|_{\Omega_x}$ for every $x \in M$,
where $\pi \colon M \to M/T$ is the projection map and
$\Omega_x=(\op{T}_x(T \cdot x))^{\sigma_x}$ (c.f. Lemma \ref{3.7}).
\item[3)]
The non--degenerate antisymmetric bilinear form 
$\sigma^{\mathfrak{t}} \colon \mathfrak{t} \times \mathfrak{t} \to \R$ on the Lie
algebra $\mathfrak{t}$ of $T$ such that for all $X, \, Y \in \mathfrak{t}$
and all $x \in M$ 
$\sigma_x(X_M(x),\, Y_M(x))=\sigma^{\mathfrak{t}}(X, \, Y)$  (c.f. Lemma \ref{sigmat}).
\item[4)]
The  monodromy invariant of $(M, \, \sigma, \, T)$, i.e. the $\mathcal{G}_{(g; \, \vec{o})}$\--orbit 
$$
 \mathcal{G}_{(g, \, \vec{o})} \cdot ((\mu_{\textup{h}}(\alpha_i), \, \mu_{\textup{h}}(\beta_i))_{i=1}^{g}, \, (\mu_{\textup{h}}(\gamma_k))_{k=1}^m),
$$
of the $(2g+m)$\--tuple $((\mu_{\textup{h}}(\alpha_i), \, \mu_{\textup{h}}(\beta_i))_{i=1}^{g}, \, (\mu_{\textup{h}}(\gamma_k))_{k=1}^m)$,
c.f. Definition \ref{nfmono}.
\end{itemize}
\end{definition}

\begin{theorem} \label{veryhardmatrix}
Suppose that $G$ is the first integral orbifold homology group of a compact, connected,
orientable smooth orbisurface of Fuchsian signature $(g; \, \vec{o}) \in \Z^{1+m}$. 
Choose a set of generators $\{\alpha_i,\, \beta_i\}_{i=1}^{g}$ of 
a maximal\footnote{This means that together with the torsion subgroup
span the entire group.} free subgroup and let $\{\gamma_k\}_{k=1}^m$
be a geometric torsion basis.
The group of geometric isomorphisms of $G$, c.f. Definition \ref{geomiso}, is equal to the group
of isomorphisms of $G$ induced by linear isomorphisms $f_B$ of 
$\Z^{2g+m}$ where $B \in \mathcal{G}_{(g; \, \vec{o})}$, 
and  $\mathcal{G}_{(g; \, \vec{o})}$ is given in Definition \ref{key2}. 
(Recall that we had to choose the generators $\alpha_i,\,
\beta_i, \, \gamma_k$ in order to define an endomorphism $f_B$ of $\Z^{2g+m}$ from the matrix $B$).
\end{theorem}

\begin{proof}
After fixing a group, the statement corresponds to that of Theorem \ref{veryhard} formulated in the language of matrices. 
\end{proof}

\subsubsection{Uniqueness Statement}

The next two results say that 
the list of ingredients of $(M,\,\sigma,\,T)$ as in Definition
\ref{def2} is a complete set of invariants of $(M,\,\sigma,\,T)$ .

\begin{lemma} \label{mt1part1}
Let $(M, \,\sigma)$ be a compact connected $2n$\--dimensional symplectic
manifold equipped
with an effective symplectic action of a $(2n-2)$\--dimensional torus $T$ for which 
at least one, and hence every
$T$\--orbit is a $(2n-2)$\--dimensional symplectic submanifold of $(M,\, \sigma)$.
Then if $(M',\,\sigma')$ is 
a compact connected $2n$\--dimensional symplectic
manifold equipped
with an effective symplectic action of $T$ for which at least one, and hence
every $T$\--orbit is a $(2n-2)$\--dimensional symplectic submanifold of $(M, \, \sigma')$,
and $(M', \, \sigma')$ is
$T$\--equivariantly symplectomorphic to $(M, \, \sigma)$, then 
the list of ingredients of $(M',\,\sigma',\,T)$ is equal to the list of ingredients of
$(M,\,\sigma,\,T)$.
\end{lemma}

The proof of Lemma \ref{mt1part1} is analogous to the proof of Lemma \ref{Fmt1part1}.

\begin{prop} \label{mt1part2}
Let $(M, \,\sigma)$ be a compact connected $2n$\--dimensional symplectic
manifold which is equipped
with an effective symplectic action of a $(2n-2)$\--dimensional torus $T$ for which 
at least one, and hence every
$T$\--orbit is a $(2n-2)$\--dimensional symplectic submanifold.
Then if $(M',\,\sigma')$ is 
a compact connected $2n$\--dimensional symplectic
manifold equipped
with an effective symplectic action of $T$ for which at least one, and hence
every $T$\--orbit is a $(2n-2)$\--dimensional
symplectic submanifold of $(M', \, \sigma')$
and the list of ingredients of $(M',\,\sigma',\,T)$ is 
equal to the list of ingredients of $(M,\,\sigma,\,T)$, then
$(M', \, \sigma')$ is $T$\--equivariantly symplectomorphic to $(M, \, \sigma)$. 
\end{prop}

\begin{proof}
Suppose that the list of ingredients of $(M, \, \sigma)$ equals the list of
ingredients of $(M', \, \sigma')$. Because $M/T$ and $M'/T$
have the same Fuchsian signature  and symplectic area, by the orbifold version of Moser's theorem \cite[Th.\,3.3]{MW}, the (compact, connected, smooth, orientable)
orbisurfaces $(M/T, \, \nu)$ and $(M'/T, \, \nu')$ are symplectomorphic,
where $\nu$ and $\nu'$ are the symplectic forms given by Lemma \ref{3.7}.

Let $\mu, \, \mu',\, \mu_{\textup{h}}, \, \mu_{\textup{h}}'$ be the 
monodromy homomorphisms from the orbifold fundamental groups
$\pi^{\textup{orb}}_1(M/T, \, p_0)$, $\pi^{\textup{orb}}_1(M'/T, \, p_0)$, and from the first integral orbifold homology groups 
$\op{H}_1^{\textup{orb}}(M/T,\, \Z)$, $\op{H}_1^{\textup{orb}}(M'/T,
\, \Z)$ into the torus $T$, respectively associated 
to the symplectic manifolds $(M, \, \sigma)$, $(M', \, \sigma')$ as
in Definition \ref{below}. Because ingredient 4) of $(M, \, \sigma)$ equals
ingredient 4) of $(M', \, \sigma')$, by Proposition \ref{conjecture01} and
Definition \ref{geomiso}
there exists an orbifold diffeomorphism
$F \colon M/T \to M'/T$ such that 
$\mu_{\textup{h}}=\mu'_{\textup{h}} \circ F_*$, and hence $\mu=\mu' \circ F_*$.

The remaining part of the proof is analogous to the proof in the free case, ``proof of
Proposition \ref{Fmt1part2}'', and the details of the arguments that follow may be found there.
If $\nu_0:=F^*\nu'$,
the symplectic orbit space $(M/T, \, \nu_0)$ is symplectomorphic to $(M'/T, \, \nu')$,
by means of $F$. Let $\widetilde{\nu_0}$
be the pullback of the $2$\--form $\nu_0$ by the orbifold universal cover $\psi \colon
\widetilde{M/T} \to M/T$
of $M/T$ based at $p_0=\psi(x_0)$, and similarly we define $\widetilde{\nu'}$ by means of the
universal cover $\psi' \colon \widetilde{M'/T} \to M'/T$ based at $F(p_0)$. 
As in the proof of Proposition \ref{Fmt1part2}, the 
orbifold symplectomorphism $F$ between $(M/T, \, \nu_0)$ and $(M'/T, \, \nu')$ lifts to a unique
symplectomorphism $\widetilde{F}$ between $(\widetilde{M/T}, \, \widetilde{\nu_0})$ 
and $(\widetilde{M'/T}, \, \widetilde{\nu'})$ such that $\widetilde{F}(x_0)=x'_0$.
Then the assignment 
$$
 [[\gamma], \, t]_{\pi^{\textup{orb}}_1(M/T, \, p_0)} \mapsto 
[\widetilde{F}([\gamma]), \, t]_{\pi^{\textup{orb}}_1(M'/T, \, p'_0)},
$$
is a $T$\--equivariant symplectomorphism between the associated
bundle $\widetilde{M/T} \times_{\pi^{\textup{orb}}_1(M/T, \,p_0)} T$ and the associated bundle
$\widetilde{M'/T} \times_{\pi^{\textup{orb}}_1(M'/T, \, p'_0)} T$, which by
by Theorem \ref{tt}, gives a $T$\--equivariant symplectomorphism 
between $(M, \, \sigma)$ and $(M', \, \sigma')$.
\end{proof}

\subsection{Existence}

\subsubsection{List of ingredients for $T$}

We assign to a torus $T$ a list of four ingredients. This is  
analogous to Definition \ref{Fdef1}.

\begin{definition} \label{def1}
Let $T$ be a torus. The \emph{list of ingredients for $T$} consists of the following.
\begin{itemize}
\item[i)]
An $(m+1)$\--tuple $(g; \, \vec{o})$ of integers,  where $\vec{o}$ is non--decreasingly ordered and consists of strictly positive integers and $m$ is a non--negative integer, and such that  
$(g; \, \vec{o})$ is not of the form $(0;\,o_1)$ or of the form $(0; \, o_1, \, o_2)$ with
$o_1 < o_2$.
\item[ii)]
A positive real number $\lambda>0$.
\item[iii)]
A non--degenerate antisymmetric bilinear form $\sigma^{\mathfrak{t}} \colon \mathfrak{t} \times
\mathfrak{t} \to \R$ on the Lie algebra $\mathfrak{t}$ of $T$.
\item[iv)]
An orbit $ \mathcal{G}_{(g, \, \vec{o})} \cdot \xi   \in 
T^{2g+m}_{(g;\, \vec{o})}/  \mathcal{G}_{(g, \, \vec{o})}$
in the Fuchsian signature space associated to $(g; \, \vec{o})$,
where
$$
T^{2g+m}_{(g;\, \vec{o})}=\{(t_i)_{i=1}^{2g+m}  \in T^{2g+m}\, | \,  \textup{the order of } t_i \textup{ is a multiple of  } o_i, 2g +1\le i \le 2g+m \},
$$
and where
 $\mathcal{G}_{(g, \, \vec{o})}$ is the group of matrices in Definition \ref{key2}.
\end{itemize}
\end{definition}

\subsubsection{Existence statement}

Any list of ingredients as in Definition \ref{def1} gives rise to one of our manifolds with symplectic
$T$\--action. 

\begin{prop} \label{existence}
Let $T$ be a  $(2n-2)$\--dimensional torus. Then given a list of ingredients for $T$, as in Definition \ref{def1}, there exists a $2n$\--dimensional symplectic manifold $(M,\, \sigma)$ with an effective
symplectic action of $T$ for which at least one, and hence every $T$\--orbit is a $(2n-2)$\--dimensional
symplectic submanifold of $(M, \, \sigma)$, and 
such that the list of ingredients of $(M,\,\sigma,\,T)$ is equal to the list of ingredients for $T$.
\end{prop}

\begin{proof} 
Let $\mathcal{I}$ be a list of ingredients for the torus $T$, as in Definition \ref{def1}.
Let the pair $(\Sigma, \, \sigma^{\Sigma})$ be a compact, connected symplectic
orbisurface of Fuchsian signature  $(g; \, \vec{o})$
given by ingredient i) of $\mathcal{I}$ in Definition \ref{def1}, and
with total symplectic area equal to the positive real number $\lambda$, 
where $\lambda$ is given by ingredient ii) of $\mathcal{I}$. By \cite[Th.\,13.3.6]{thurston3},
since $(g; \, \vec{o})$ is not of the form $(0;\,o_1), \, (0; \, o_1, \, o_2)$ 
with $o_1 <o_2$ by assumption, $\Sigma$ is a (very) good orbisurface.
Let the space $\widetilde{\Sigma}$
be the orbifold universal cover of $\Sigma$, which is a smooth manifold,
based at an arbitrary regular point $p_0 \in \Sigma$ which we fix for the rest of the proof,
c.f. the construction we gave prior to Definition \ref{themodel} and Theorem \ref{tt}.
Because $\Sigma$ is very good orbisurface, $\widetilde{\Sigma}$
is a smooth surface.
Let $\alpha_i, \, \beta_i$, for $1 \le i \le 2g$, be a symplectic basis of a free subgroup 
of $\op{H}_1^{\scriptop{orb}}(\Sigma,\,\Z)$
which together
with the torsion subgroup span the entire group, and let
$\gamma_k$, for $1 \le k \le m$, be a geometric torsion basis and 
define
$\mu_{\textup{h}} \colon \op{H}_1^{\scriptop{orb}}(\Sigma, \, \Z) \to T$ to be the unique homomorphism such that
$\mu_{\textup{h}}(\alpha_i)=a_i, \, \mu_{\textup{h}}(\beta_i)=b_i,
\, \mu_{\textup{h}}(\gamma_k)=c_k$, where
the tuple $((a_i,\, b_i), \, (\gamma_k))$ 
represents ingredient iv) of $\mathcal{I}$, and $h_1$
denotes the orbifold Hurewicz homomorphism from 
$\pi^{\textup{orb}}_1(\Sigma, \, p_0) $ to
$\op{H}_1^{\textup{orb}}(\Sigma, \, \Z)$.
Let
$\mu \colon \pi^{\textup{orb}}_1(\Sigma, \, p_0) \to T$ 
be the homomorphism defined as $\mu:=\mu_{\textup{h}} \circ h_1$. 
Let the orbifold fundamental group 
$\pi^{\textup{orb}}_1(\Sigma, \, p_0)$ act on the smooth manifold given as the
Cartesian product $\widetilde{\Sigma} \times T$ by 
the diagonal action  $[\delta] \cdot ([\gamma],\, t)=([\delta \, \gamma^{-1}],\, \mu([\delta]) \, t)$.
Define the bundle space 
\begin{eqnarray} \label{modell}
M^{\Sigma}_{\textup{model},\,p_0}:=
\widetilde{\Sigma} \times_{\pi^{\textup{orb}}_1(\Sigma, \, p_0)}T.
\end{eqnarray}
The symplectic form
and torus action on $M^{\Sigma}_{\textup{model}, \, p_0}$ are 
constructed in the exact same way as in the free case, c.f. proof of 
Proposition \ref{Fexistence}.

The proof that ingredients 1)--3) of $(M^{\Sigma}_{\textup{model},\,p_0}, \, \sigma^{\Sigma}_{\textup{model}})$ are equal to
ingredients 1)--3) of the list $\mathcal{I}$ is the same as in the free case (c.f. proof of
Proposition \ref{Fexistence}), with the observation
that in the non--free case we use the classification theorem of compact, connected, smooth
orientable orbisurfaces Theorem \ref{2dorbifoldtheorem} of Thurston's instead of the classical
classification theorem of compact, connected, smooth surfaces, to prove that the corresponding
ingredients 1) agree. We have left
to show that ingredient 4) of 
$(M^{\Sigma}_{\textup{model},\,p_0}, \, \sigma^{\Sigma}_{\textup{model}})$ equals
ingredient 4) of the list $\mathcal{I}$. 
Let $\Omega^{\Sigma}_{\textup{model}}$
stand for the flat connection on $M^{\Sigma}_{\textup{model}}$ given by the symplectic orthogonal
complements to the tangent spaces to the $T$\--orbits, see Proposition \ref{easynow}, and let
$\mu^{\Omega^{\Sigma}_{\textup{model}}}_{\textup{h}}$ stand for the induced homomorphism
$\mu^{\Omega^{\Sigma}_{\textup{model}}}_{\textup{h}} \colon \op{H}^{\textup{orb}}_1(M^{\Sigma}_{\textup{model}}/T, \, \Z)
\to T$ in homology, by the monodromy of such connection. 
If $f \colon \op{H}^{\textup{orb}}_1(M^{\Sigma}_{\textup{model}}/T,
\, \Z) \to \op{H}^{\textup{orb}}_1(\Sigma, \, \Z)$ is
the group isomorphism induced by the orbifold symplectomorphism
\begin{eqnarray} \label{compositemap2}
M^{\Sigma}_{\textup{model}}/T \to \widetilde{\Sigma}
/\pi^{\textup{orb}}_1(\Sigma, \, p_0) \to \Sigma,
\end{eqnarray}
where each arrow in (\ref{compositemap2}) represents the natural map,
\begin{eqnarray} \label{monoequation2}
\mu_{\textup{h}}^{\Omega^{\Sigma}_{\textup{model}}}=\mu_{\textup{h}} \circ f.
\end{eqnarray}
Because $f$ is induced by a diffeomorphism, 
by Theorem \ref{veryhard}, $f$ is symplectic and torsion geometric, c.f. Definition \ref{symiso}
and Definition \ref{preiso}. Therefore there exists a unique collection of elements
$\alpha'_i, \, \beta'_i$, $1 \le i \le g$, in $\op{H}^{\textup{orb}}_1(M^{\Sigma}_{\textup{model}}/T,
\, \Z)$ such that $f(\alpha'_i)=\alpha_i$ and
$f(\beta'_i)=\beta_i$, for all $1 \le i \le g$. The elements $\alpha'_i, \, \beta'_i$, $1 \le i \le g$,
form a symplectic basis of a free subgroup $F^{\Omega^{\Sigma}}$ of the orbifold homology group 
$\op{H}^{\textup{orb}}_1(M^{\Sigma}_{\textup{model}}/T, \, \Z)$,
which together with the torsion subgroup spans the entire group. Similarly let the collection
$\gamma'_k$, for $1 \le k \le m$, be such that $f(\gamma'_k)=\gamma_k$, for all $1 \le k \le m$.
The $\gamma'_k$, $1 \le k \le m$, form a geometric torsion basis, c.f. Definition \ref{geometrictorsionbasis} such that $o_k=o'_{\tau(k)}$ for all $k$, $1 \le k \le m$, 
for a permutation $\tau \in \textup{S}^{\vec{o}}_m$.
Let $\widehat{\xi}$ be the $(2g+m)$\--tuple of elements 
$\mu_{\textup{h}}^{\Omega^{\Sigma}_{\textup{model}}}(\alpha'_i)$,
$\mu_{\textup{h}}^{\Omega^{\Sigma}_{\textup{model}}}(\beta'_i),\, 
\mu_{\textup{h}}^{\Omega^{\Sigma}_{\textup{model}}}(\gamma'_k)$,
where $1 \le i \le g$ and $1 \le k \le m$.
Therefore by (\ref{monoequation2})
$$
\widehat{\xi}=((\mu_{\textup{h}}(\alpha_i),\,
\mu_{\textup{h}}(\beta_i))_{i=1}^g, (\mu_{\textup{h}}(\gamma_k))_{k=1}^m),
$$
which in particular implies that ingredient 4) of 
$(M^{\Sigma}_{\textup{model},\,p_0}, \, \sigma^{\Sigma}_{\textup{model}})$ equals
$((a_i,\,b_i),\, c_k)$.
\end{proof}

\subsection{Classification Theorem}

By putting
together the results of the previous subsections, we obtain
the main result of the section:

\begin{theorem} \label{TT}
Let $T$ be a $(2n-2)$\--dimensional torus. Let $(M,\,\sigma)$ be a compact connected 
$2n$\--dimensional symplectic manifold 
on which $T$ acts effectively and symplectically and such that at least one,
and hence every $T$\--orbit is a $(2n-2)$\--dimensional symplectic
submanifold of $(M, \, \sigma)$.

Then the list of ingredients of $(M,\,\sigma,\,T)$ as in Definition
\ref{def2} is a complete set of invariants of $(M,\,\sigma,\,T)$, 
in the sense that, if $(M',\,\sigma')$ is 
a compact connected $2n$\--dimensional symplectic
manifold equipped
with an effective symplectic action of $T$ for which at least one, and hence every
$T$\--orbit is a $(2n-2)$\--dimensional symplectic submanifold of $(M', \, \sigma')$,
$(M', \, \sigma')$ is
$T$\--equivariantly symplectomorphic to $(M, \, \sigma)$ if and only if 
the list of ingredients of $(M',\,\sigma',\,T)$ is equal to the list of ingredients of
$(M,\,\sigma,\,T)$.

And given a list of ingredients for $T$, as in Definition \ref{def1}, there 
exists a symplectic $2n$\--dimensional 
manifold $(M,\, \sigma)$ with an effective symplectic torus action of $T$ for which at least
one, and hence every 
$T$\--orbit is a $(2n-2)$\--dimensional symplectic submanifold of $(M, \, \sigma)$, 
such that the list of ingredients of $(M,\,\sigma,\,T)$ is equal to the
list of ingredients for $T$.
\end{theorem}
\begin{proof}
It follows by putting together Lemma \ref{mt1part1}, Proposition \ref{mt1part2} and Proposition 
\ref{existence}. The combination of Lemma \ref{mt1part1}, Proposition \ref{mt1part2}
gives the existence part of the theorem, while Proposition \ref{existence} gives the existence
part.
\end{proof}

\begin{remark}
The author is grateful to 
professor Deligne for pointing out an imprecision
in a earlier version of the following statements.
Let $T$ be a $(2n-2)$\--dimensional 
torus. Let $\mathcal{M}$ denote the category of which 
the objects are the compact connected symplectic $2n$\--dimensional manifolds 
$(M,\,\sigma )$  together with an effective symplectic 
$T$\--action on $(M,\,\sigma )$ such that at least one, and hence
every $T$\--orbit is a $(2n-2)$\--dimensional symplectic 
submanifold of $M$, and of which the 
morphisms are the $T$\--equivariant symplectomorphisms
of $(M, \, \sigma)$. 
Let $\mathcal{I}$ denote the set of all lists of 
ingredients as in Definition \ref{def1}, 
viewed as a category, and of which the identities are
the only endomorphisms of categories. 
Then the assignment $\iota$ in Definition \ref{def2} 
is a full functor categories from the
category $\mathcal{M}$ onto the category
$\mathcal{I}$.
In particular the proper class $\mathcal{M}/ \sim$ of isomorphism classes in $\mathcal{M}$ 
is a set, and the functor $\iota  \colon \mathcal{M} \to \mathcal{I}$ 
in Definition \ref{def2}  induces a bijective mapping $\iota/\sim$ 
from $\mathcal{M}/\sim$ onto $\mathcal{I}$. 
The fact that the mapping $\iota  \colon \mathcal{M} \to \mathcal{I}$ is a functor 
and the mapping $\iota/ \sim$ is injective 
follows from the uniqueness part of the statement of Theorem \ref{TT}. The surjectivity of $\iota$, 
follows from the existence part of the statement of Theorem \ref{TT}. 
\end{remark}

\section{The four dimensional classification}

We give a classification of effective symplectic actions of $2$\--tori
on compact, connected, symplectic $4$\--dimensional manifolds under
no additional assumption.

This generalizes the $4$\--dimensional case of Delzant's theorem
on the classification of symplectic--toric manifolds
to symplectic actions which are not Hamiltonian.

\subsection{Classification Statement}

We state the main theorem.

The following is a special case of \cite[Def. 9.1]{DuPe}. 
We call the list ``free and Lagrangian''
for obvious reasons, c.f. Theorem \ref{oldclassificationtheorem}.

\begin{definition}  \label{def3ing}
Let $T$ be a $2$\--torus. 
A {\em free and Lagrangian list of ingredients for $T$} consists of: 
\begin{itemize} 
\item[1)] A discrete cocompact subgroup $P$ of $\mathfrak{t}^*$.
\item[2)] An antisymmetric bilinear mapping 
$c: \mathfrak{t}^*\times \mathfrak{t}^*\to\mathfrak{t}$ 
with the following properties. 
\begin{itemize}
\item[2a)]  
For every $\zeta ,\,\zeta '\in P$, the element 
$c(\zeta ,\,\zeta ')\in \mathfrak{t}$ belongs to the integral lattice 
$T_{{\Z}}$ in $\mathfrak{t}$, the kernel of the exponential 
mapping $\textup{exp}:\mathfrak{t}\to T$. 
\item[2b)] For every $\zeta ,\,\zeta ',\,\zeta ''\in \mathfrak{t}^*$ 
we have that 
\[
\zeta (c(\zeta ',\,\zeta ''))
+\zeta '(c(\zeta '',\,\zeta ))
+\zeta ''(c(\zeta ,\,\zeta '))=0.
\]
\end{itemize}
\item[3)]  An element in $\mathcal{T}$, c.f. (\ref{calT}).
\end{itemize}
\label{ingredientdef}
\end{definition}

The following definition is a particular case of \cite[Def. 7.9]{DuPe}. 
%(Therein we take $\mathfrak{l}=\mathfrak{t}$ and $N=\mathfrak{t}^*$).

\begin{definition}
Let $c \colon \got{t}^* \times \got{t}^* \to \got{t}$ be as in item 2) of Definition \ref{def3ing}.
Let $\op{Hom}_c(P,\, T)$ denote the space of mappings 
$
\tau:P\to T, \, \zeta\mapsto \tau_{\zeta}, 
$
such that 
\begin{equation}
\tau _{\zeta '}\,\tau _{\zeta}=\tau _{\zeta +\zeta '}
\,\op{e}^{c(\zeta ',\,\zeta )/2},\quad \zeta,\,\zeta '\in P. 
\label{tauplushom}
\end{equation} 
For each $\zeta' \in \mathfrak{t}^*$, $\zeta\mapsto 
c(\zeta ,\,\zeta ')$ is a homomorphism from $P$ to 
$\got{t}$. Write 
$c(\cdot ,\, \mathfrak{t}^*)$ for the set of all $c(\cdot ,\,\zeta ')
\in\op{Hom}(P,\,\got{t})$ such that $\zeta '\in \mathfrak{t}^*$. 
$c(\cdot ,\, \mathfrak{t}^*)$ is a linear subspace of the Lie 
algebra $\op{Hom}(P,\,\got{t})$ of $\op{Hom}(P,\, T)$. 
If $h:\zeta\mapsto h_{\zeta}$ is a homomorphism from $P$ to $T$, 
then $h\cdot\tau :\zeta\mapsto \tau _{\zeta}\, h_{\zeta} 
\in\op{Hom}_c(P,\, T)$ for every $\tau\in\op{Hom}_c(P,\, T)$, 
and $(h,\,\tau )\mapsto h\cdot\tau$ defines a 
free, proper, and transitive action of $\op{Hom}(P,\, T)$ 
on $\op{Hom}_c(P,\, T)$. 
Let $\op{Sym}$ denote the space of all linear mappings 
$\alpha :\got{t}^*\to\got{t}$, $\xi \mapsto \alpha_{\xi}$, which are symmetric in the 
sense of 
\begin{equation}
\xi(\alpha_{\xi'})-\xi'(\alpha_{\xi}) = 0.
\label{alphasym}
\end{equation}
For each $\alpha\in\op{Sym}$, 
the restriction $\alpha |_P$ of $\alpha$ to $P$ is a 
homomorphism from $P$ to $\got{t}$. 
In this way the set $\op{Sym}|_P$ of all $\alpha |_P$ 
such that $\alpha\in\op{Sym}$ is another linear 
subspace of $\op{Hom}(P,\,\got{t})$. 
Write 
\begin{equation}
{\cal T}:=\op{Hom}_c(P,\, T)/\op{exp}{\cal A},
\quad {\cal A}:=c(\cdot ,\, \mathfrak{t}^*)+\op{Sym}|_P
\label{calT}
\end{equation} 
for the orbit space of the action of  
the Lie subgroup $\op{exp}{\cal A}$ 
of $\op{Hom}(P,\, T)$ on $\op{Hom}_c(P,\, T)$. 
\end{definition}

\begin{remark}
Because $c(P\times P)\subset T_{\Z}$, 
the factor $\op{e}^{c(\zeta ',\zeta )/2}$ in (\ref{tauplushom}) 
is an element of order two in $T$. 
The elements of $\op{Hom}_c(P,\, T)$ are homomorphisms from $P$ to $T$
if $c(P\times P)\subset 2T_{\Z}$.
$\op{Hom}(P,\, T)$ is 
a torus group, a compact, connected, and commutative 
Lie group with Lie algebra equal to the 
vector space $\op{Hom}(P,\,\got{t})$ of dimension 
$(\op{dim}T)^2$. 
$\op{Hom}_c(P,\, T)$ 
is diffeomorphic to a torus of dimension $(\op{dim}T)^2$. 
\end{remark}

The following is our main result of this paper, a classification of symplectic actions of $2$\--tori on compact, connected symplectic $4$\--manifolds, up to $T$\--equivariant symplectomorphisms. 
Theorem \ref{520} consists of two separate statements. 
The first statement says that the examples in i) and ii) are symplectic $T$\--manifolds with certain characteristics.  Moreover it says that any symplectic $4$\--manifold with a $2$\--torus action is of the form given in i), ii), 1) or 2). The second statement
says that for any two different lists of ingredients (a discrete cocompact subgroup, a bilinear form etc.) 
as in i) or ii) give
rise to symplectic manifolds which are not $T$\--equivariantly symplectomorphic. Hence Theorem \ref{520} is a full classification.

\begin{theorem} \label{520}
Let $T$ be a $2$\--dimensional torus. Then the following hold:
 \begin{itemize}
 \item[i)]
For any choice of a discrete cocompact subgroup $P$ 
of $\mathfrak{t}^*$,  a bilinear form $c \colon\mathfrak{t}^* \times \mathfrak{t}^* \to
\mathfrak{t}$ and an equivalence 
class of a mapping $[\tau]_{\textup{exp}(\mathcal{A})}$, 
as in Definition \ref{ingredientdef}, let  $\iota \colon P \to T \times \mathfrak{t}^*$
be the mapping $\zeta \mapsto (\tau_{\zeta}^{-1},\,\zeta)$,
which is a homomorphism onto a discrete cocompact subgroup of $T\times \got{t}^*$
with respect to the non--abelian group structure given by 
$$
(t,\,\zeta )\, (t',\,\zeta ')=(t\, t'\,\op{e}^{-c(\zeta ,\,\zeta ')/2},\,\zeta +\zeta ').
$$
Equip $T \times \mathfrak{t}^*$ with the standard cotangent bundle symplectic form.
Then 
$(T \times \mathfrak{t}^*)/\iota(P)$ equipped with the action of $T$
which comes from the action of $T$ by translations on the left factor of $T \times \mathfrak{t}^*$,
and where the symplectic form on $(T \times \mathfrak{t}^*)/\iota(P)$ is the $T$\--invariant form
induced by the symplectic form on $T \times \mathfrak{t}^*$, 
is a compact, connected symplectic $4$\--manifold on which $T$ acts freely and for which the
$T$\--orbits are Lagrangian $2$\--tori. 
\item[ii)]
For any choice of
an $(m+1)$\--tuple $(g; \, \vec{o})$ of integers, a positive real number $\lambda>0$, 
a non--degenerate antisymmetric bilinear form $\sigma^{\mathfrak{t}}$ on $\mathfrak{t}$,
and an element $\mathcal{G}_{(g, \, \vec{o})} \cdot \xi   \in 
T^{2g+m}_{(g;\, \vec{o})}/  \mathcal{G}_{(g, \, \vec{o})}$ where
$\xi=((a_i,\, b_i), \, (c_k))$ as in Definition \ref{def1}, let $\Sigma$ be an orbisurface with
Fuchsian signature $(g;\, \vec{o})$, and total symplectic area $\lambda$, and 
let $p_0 \in \Sigma$. 
Let $\alpha_i, \, \beta_i$, for $1 \le i \le 2g$, be a symplectic basis of a free subgroup 
of $\op{H}_1^{\op{orb}}(\Sigma, \, \Z)$
which together
with the torsion subgroup spans $\op{H}_1^{\op{orb}}(\Sigma, \, \Z)$, and let
$\gamma_k$, for $1 \le k \le m$, be a geometric torsion basis (c.f. Definition \ref{geometrictorsionbasis}).
Let $h_1$ be the Hurewicz homomorphism.
Let $f_{\textup{h}}$ to be the unique homomorphism such that
$f_{\textup{h}}(\alpha_i)=a_i, \, f_{\textup{h}}(\beta_i)=b_i,
\, f_{\textup{h}}(\gamma_k)=c_k$.
Let $f:=f_{\textup{h}} \circ h_1$. 
Let $\pi^{\textup{orb}}_1(\Sigma, \, p_0)$ act on $\widetilde{\Sigma} \times T$ by 
$$
[\delta] \cdot ([\gamma],\, t)=([\delta \, \gamma^{-1}],\, f([\delta]) \, t).
$$
Equip the universal cover
$\widetilde{\Sigma}$ with the symplectic form pullback from $\Sigma$,
and $\widetilde{\Sigma} \times T$ with the product symplectic form.
Let
$T$ act by translations on the right factor of $\widetilde{\Sigma} \times T$. 
Then the space $\widetilde{\Sigma} \times_{\pi^{\textup{orb}}_1(\Sigma, \, p_0)}T$
endowed with the unique 
symplectic form and $T$\--action induced by the product ones
is a compact, connected 
symplectic $4$\--manifold on which $T$ acts effectively and for which the 
$T$\--orbits are symplectic $2$\--tori.
\end{itemize}
Moreover, if $(M, \, \sigma)$ is a compact connected symplectic $4$\--dimensional manifold equipped with an effective symplectic action of a $2$\--torus $T$, then one and only one of the following
cases occurs:
\begin{itemize}
\item[1)]
$(M, \, \sigma)$ is a $4$\--dimensional symplectic toric manifold, hence determined
up to $T$\--equivariant symplectomorphisms by its Delzant polygon $\mu(M)$ centered
at the origin, where $\mu \colon M \to \mathfrak{t}^*$ is the momentum map for the $T$\--action. 
\item[2)]
$(M, \, \sigma)$ is equivariantly symplectomorphic to
a product $\T^2 \times S^2$, where
$\T^2=(\R/\Z)^2$ and
the first factor of $\mathbb{T}^2$ acts on the left factor by translations
on one component, and the second factor acts on $S^2$ by rotations about the
vertical axis of $S^2$. The symplectic form is
a positive linear combination of the
standard translation invariant form on $\mathbb{T}^2$
and the standard rotation invariant form on $S^2$.
\item[3)]
$(M, \, \sigma)$ is $T$\--equivariantly symplectomorphic to one and only one of the
symplectic $T$\--manifolds in part i), given by a unique list of ingredients.
\item[4)]
 $(M, \, \sigma)$ is $T$\--equivariantly symplectomorphic to one and only one of the
symplectic $T$\--manifolds in part ii), given by a unique list of ingredients.
\end{itemize}
\end{theorem}

\subsection{Proof of Theorem \ref{520}}

 We make use of results  in sections 3, 4 and ideas/methods in 
proofs of \cite[Prop. 5.5, Lem. 7.1, Lem. 7.5]{DuPe}. Because we do not reprove anything
which appeared therein, it is unrealistic to expect understanding every detail of the proof below 
without prior knowledge of them. References are given to specific parts of the previous sections and
to \cite{DuPe} when appropiate. We divide the proof into three steps.
\\
\\
\\
{\large {\bf Step 1.}}
First suppose that the $2$\--dimensional $T$\--orbits which are Lagrangian submanifolds
of $(M, \, \sigma)$. 
The ingredients of the construction
of the model $(M_{\scriptop{model}}, \, \sigma_{\scriptop{model}})$
in Claim 1.1 below,
which we obtain as a particular case of  \cite[Prop.\,7.2, Prop.\,7.4]{DuPe},
are the following.
$(M_{\textup{h}}, \, \sigma_{\scriptop{h}})$ is a Delzant submanifold of $(M, \, \sigma)$,
where $\sigma_{\scriptop{h}}$ equals $\sigma$ restricted to $M_{\scriptop{h}}$;
$T_{\textup{h}}$ is the maximal subtori of $T$ which acts on $M$ in a Hamiltonian fashion.
 Because $T$ is $2$\--dimensional, it follows from \cite[Lem. 2.6, Lem. 3.11]{DuPe} 
 that $T=T_{\textup{h}}$ if the $T$\--action is Hamiltonian, $T=\{1\}$ if the action is free, and
 otherwise $T_x=T_{\textup{h}}$ for any $x$ such that $T_x \neq \{1\}$;
$N=(\mathfrak{t}/\mathfrak{t}_{\textup{h}})^*$; $G$ is the Lie group
\begin{eqnarray} \label{G}
G:=T \times N, 
\end{eqnarray}
which is two--step nilpotent with respect to the non--abelian product
\begin{equation}
(t,\,\zeta )\, (t',\,\zeta ')=
(t\, t'\,\op{e}^{-c(\zeta ,\,\zeta ')/2},\,\zeta +\zeta '),
\label{Hproduct}
\end{equation}
where 
\begin{eqnarray} \label{thec}
c \colon N \times N \to \mathfrak{t},
\end{eqnarray}
is a certain antisymmetric bilinear
form which represents the Chern class of the principal bundle
$M_{\textup{reg}} \to M_{\textup{reg}}/T$; the so called period lattice
\begin{eqnarray} \label{P}
P \le N,
\end{eqnarray}
which is a discrete, cocompact, subgroup;
$H \le G$ is a commutative, closed, Lie subgroup, defined in terms of the 
holonomy $\tau$, which lives in certain exponential quotient of the, roughly speaking,
homomorphism group $\textup{Hom}_c(P,\,T)$, c.f. (\ref{calT});
$H$ acts on $G \times M_{\textup{h}}$ by 
\begin{eqnarray} \label{holo}
((t,\,\zeta),\,x) \mapsto (t \, \tau_{\zeta}) \cdot x.
\end{eqnarray}
In a formula,
\begin{eqnarray} \label{H}
H=\{(t,\,\zeta) \, |\,  \xi \in P, \, t \, \tau_{\zeta} \in T_{\textup{h}} \  \} \le G.
\end{eqnarray}
The map 
\begin{equation}
((t,\,\zeta ),\, x)\mapsto (t\,\tau _{\zeta})\cdot x:
H\times M_{\scriptop{h}}\to M_{\scriptop{h}}
\label{Hact}
\end{equation}
defines a smooth action of $H$ on the Delzant manifold 
$M_{\scriptop{h}}$. 
The symplectic form $\sigma_{\textup{model}}$ on $M_{\textup{model}}:=G \times_H M_{\scriptop{h}}$
is given by a explicit formula as follows, c.f. \cite[Prop.\,7.4]{DuPe}, where
$X \mapsto X_{\textup{h}}$ is a
projection onto $\mathfrak{t}_{\textup{h}}$ (this depends
on a choice of freely acting torus $T_{\textup{f}}$, which induces a
decomposition of $\mathfrak{t}$, c.f. \cite[Sec.\,5.2]{DuPe}). 
First consider the form on $G \times M_{\textup{h}}$
\begin{eqnarray} \label{symplecto}
\omega_a(\delta a,\,\delta'a)&=&\delta \zeta(X')-\delta' \zeta(X) \nonumber \\
& & - \mu(x)(c_{\textup{h}}(\delta \zeta,\, \delta \zeta'))  
+ (\sigma_{\textup{h}})_x(\delta x,\,(X'_{\textup{h}})_{M_{\textup{h}}}(x))  \nonumber \\
& &-(\sigma_{\textup{h}})_x(\delta' x,(X_{\textup{h}})_{M_{\textup{h}}}(x))
+(\sigma_{\textup{h}})_x(\delta x,\,\delta'x),
\end{eqnarray}
where $c_{\scriptop{h}}$ is the $\got{t}_{\scriptop{h}}$\--component of $c$
 and $\mu$ is the momentum
map for the Hamiltonian $T_{\textup{h}}$\--action, c.f. \cite[Lem. 3.11]{DuPe}.
Here $\delta a=((\delta t,\,\delta\zeta ),\,\delta x)$ 
and $\delta' a=((\delta 't,\,\delta '\zeta ),\,\delta ' x)$ 
are tangent vectors to $G\times M_{\scriptop{h}}$ at 
$a=((t,\,\zeta ),\, x)$, 
where we identify each tangent space of the torus $T$ with 
$\got{t}$, and $X=\delta t+c(\delta\zeta ,\,\zeta )/2$ 
and $X'=\delta 't+c(\delta '\zeta ,\,\zeta )/2$.
Let $\pi_{M_{\textup{model}}} \colon G \times M_{\textup{h}} \to M_{\textup{model}}$ be the
canonical projection. Then there exists a unique $2$\--form $\sigma_{\textup{model}}$ such
that 
\begin{eqnarray} \label{sm}
\omega=\pi_{M_{\textup{model}}}^* \, \sigma_{\textup{model}}.
\end{eqnarray}
The translational $T$\--action on $M_{\textup{model}}$ by left multiplications
on $T$ is so that 
$T_{\textup{h}}$ acts on the Delzant manifold
$M_{\textup{h}}$, and its complement
$T_{\textup{f}}$ in $T$, which acts freely, permutes the Delzant submanifolds.
\\
\\
{\bf Claim 1.1.} There exists a $T$\--equivariant symplectomorphism from
$(M_{\textup{model}}, \sigma_{\textup{model}})$ onto $(M, \, \sigma)$.
\\
\\
\emph{Proof of Claim 1.1.}(Sketch).
The first observation is that the orbit space $M/T$ is a  polyhedral $\mathfrak{t}^*$\--parallel space.
A \emph{$\mathfrak{t}^*$\--parallel space} is a
Hausdorff, topological space modelled on a corner of $\mathfrak{t}^*$, c.f. \cite[Def.\,10.1]{DuPe}. 
The local charts $\phi_{\alpha}$ into $\mathfrak{t}^*$,
satisfy that the mapping
$x \mapsto \phi_{\alpha}(x)-\phi_{\beta}(x)$ has to be locally constant for all
values of $\alpha$ and $\beta$. Consider the form
$\hat{\sigma}(X):=-\textup{i}_{X_M} \sigma \in \Omega^1(M)$,
which is a closed, basic form, if $X \in \mathfrak{t}$.  
The assignment 
$\hat{\sigma} \colon x \mapsto (\hat{\sigma}_x \colon \op{T}_xM \to \mathfrak{t}^*)$, where 
we use the identification $\hat{\sigma}_x \simeq \op{T}_x\pi$,
induces and isomorphism $\hat{\sigma}_p \colon \op{T}_p(M_{\textup{reg}}/T) \to \mathfrak{t}^*$. 
This implies that 
a ``constant vector field on $\mathcal{X}^{\infty}(M_{\textup{reg}}/T)$''
can be thought of as an element ``$\xi \in \mathfrak{t}^*$''.
$L_{\xi}$ is a \emph{lift of $\xi$} if $\hat{\sigma}_x(L_{\xi})=\xi$.
Secondly, as $\mathfrak{t}^*$\--parallel spaces, 
there is an isomorphism
$ M/T \simeq \Delta \times S$, where $\Delta$ is a Delzant polytope, and $S$ is a torus, c.f. \cite[Prop.\,3.8, Th.\,10.12]{DuPe}.
Proving this involves the classification of $V$\--parallel spaces.
Moreover, it involves generalizing the Tietze--Nakajima theorem 
in \cite{tietze}, \cite{nakajima}.
Here the Delzant polytope $\Delta$ captures the Hamiltonian part of the action, while 
the torus $S$ captures the free part.

Recall that $M_{\textup{reg}}$ is the subset of $M$ where the $T$\-action is free.
In \cite[Prop.\,5.5]{DuPe} we showed that there exists a  ``nice'' admissible connection 
$
\xi \in \mathfrak{t}^* \mapsto L_{\xi} \in \mathcal{X}^{\infty}(M_{\textup{reg}}),
$
\cite[Def.\,5.3]{DuPe}, for the principal $T$\--bundle $\pi \colon M_{\textup{reg}} \to M_{\textup{reg}}/T$,
The lifts $L_{\xi}$, $\xi \in N$, have smooth extensions to $M$.
Here ``nice'' means with simple Lie brackets $[L_{\xi}, \,L_{\eta}]$, such as 
$[L_{\xi}, \,L_{\eta}]=c(\xi,\,\eta)_M$, $\xi,\,\eta \in N$, where $c$ as in (\ref{thec})
represents the Chern class of $\pi \colon M_{\textup{reg}} \to M_{\textup{reg}}/T$, and
zero otherwise.
We also achieve simple symplectic pairings $\sigma(L_{\xi}, \, L_{\eta})$. In the
particular case that $c_{\scriptop{h}}=0$, we have that 
$
\sigma(L_{\xi},\,L_{\eta})=0, \, \, \, \textup{for all } \xi, \, \eta \in \mathfrak{t}^*.
$
The lifts $L_{\xi}$, $\xi \notin N$, are singular on $M \setminus M_{\textup{reg}}$.
The singularities are also required to be simple. This step
involves computations in the cohomology ring $\op{H}^*(M/T)$.

Then in \cite[Prop.\,6.1]{DuPe} we define the integrable distribution$\{D_x\}_{x \in M}$on $M$
where
$
D_x
$
is the span of $L_{\eta}(x),\,Y_M(x)$ as $Y \in \mathfrak{t}_{\textup{h}}, \, \eta \in C$,
where
$
C \oplus N = \mathfrak{t}^*
$. The integral manifolds of this distribution are $T$\--equivariantly symplectomorphic 
to the Delzant manifold $(M_{\textup{h}},\,\sigma_{\textup{h}}, \,T_{\textup{h}})$.
Assuming this, we make the definition of $G$ in (\ref{G}), 
$H$ in (\ref{H}), and $M_{\textup{model}}$ from the connection $\xi \mapsto L_{\xi}$, and the
definition of $\sigma_{\textup{model}}$ in (\ref{sm}).
The definition of $H$ involves the holonomy of the connection $\xi \mapsto L_{\xi}$.
And the $T$\--equivariant symplectomorphism between the model $G \times_H M_{\textup{h}}$
and our symplectic $T$\--manifold $(M, \, \sigma)$ is induced by
\begin{eqnarray} \label{tsts}
((t,\,\xi),\,x) \mapsto t \cdot \textup{e}^{L_{\xi}}(x) \colon G \times M_{\scriptop{h}} \to M.
\end{eqnarray}
$P$ is rigorously introduced in Lemma 10.12, 
Proposition 3.8 in \cite{DuPe}; $c$ is introduced in
in Proposition 5.5 therein.
\\
\\
\\
{\large {\bf Step 2.}}
In Step 1 we introduced ingredients i) and ii) below.
The definition of ingredient iii) takes Section 7.5 in \cite{DuPe}.\footnote{Recall: the 
proof of Theorem \ref{520}
depends on Definition \ref{invariantdef} below, and in particular
in item iii), but its statement did not.}

\begin{definition} \label{invariantdef}
Let $(M, \, \sigma)$ be a compact, connected symplectic manifold equipped with
a free symplectic action of a $2$\--torus $T$ for which the
$T$\--orbits are Lagrangian submanifolds of $(M, \, \sigma)$. 
The \emph{list of ingredients of $(M,\,\sigma ,\, T)$}, 
as in Definition \ref{ingredientdef}, consists of:   
\begin{itemize}
\item[i)] The period lattice group $P$ in $\mathfrak{t}^*$, where $\mathfrak{t}^*$ 
acts on the orbit space $M/T$, c.f. (\ref{P}).
\item[ii)] The
antisymmetric bilinear mapping $c \colon \mathfrak{t}^*\times \mathfrak{t}^* \to\mathfrak{t}$,
c.f. (\ref{thec}) with $\got{t}_{\scriptop{h}}=\{0\}$ and so $N=\got{t}^*$.
\item[iii)] The so called \emph{holonomy invariant of $(M,\,\sigma ,\, T)$},
which is an element in $\mathcal{T}$, c.f. (\ref{calT}).
\end{itemize}
\end{definition}

\cite[Th.\,9.4, Th.\,9.6]{DuPe} imply the following.

\begin{theorem}
 \label{oldclassificationtheorem}
Let $T$ be a $2$\--torus. 
The list of ingredients of $(M,\,\sigma ,\, T)$ 
is a complete set of invariants for the compact 
connected symplectic manifold $(M,\,\sigma )$ with 
free symplectic $T$\--action with Lagrangian orbits, 
in the following sense.

If $(M',\,\sigma ')$ is another 
compact connected symplectic manifold with a free symplectic $T$\--action 
with Lagrangian orbits, then there exists a $T$\--equivariant 
symplectomorphism from $(M,\,\sigma)$ onto 
$(M',\,\sigma ')$ if and only if the list of ingredients 
of $(M,\,\sigma ,\, T)$ is equal to the 
list of ingredients of $(M',\,\sigma ',\, T)$.  

Moreover,
every list of ingredients as in Definition 
\ref{ingredientdef} 
is equal to the list of 
invariants of a compact connected symplectic manifold $(M,\,\sigma )$ 
with free symplectic 
$T$\--action with Lagrangian orbits.
\end{theorem}
The next step combines Step 1 and Step 2, together with the results in sections 3, 4.
\\
\\
\\
{\large {\bf Step 3.}} Throughout we refer to the items in the statement of Theorem \ref{520}.
It follows from the definitions that
the mapping $\iota :\zeta\mapsto ({\tau_{\zeta}}^{-1},\,\zeta )$ in part i)
is a homomorphism from $P$ onto a discrete cocompact subgroup 
of $T \times \got{t}^*$. 
Proving that the space defined in part i) is a compact, connected
symplectic $4$\--manifold equipped with a effective symplectic action is also an exercise
using the definitions. Similarly, it follows from 
the pointwise expression for the symplectic form on $T \times \mathfrak{t}^*$ that 
the symplectic form on $(T \times \mathfrak{t}^*)/\iota(P)$
vanishes along the $T$\--orbits, which hence
are isotropic submanifolds of $(T \times \mathfrak{t}^*)/\iota(P)$. 
The action on $T \times \mathfrak{t}^*$ is free, and passes
to a free action on $(T \times \mathfrak{t}^*)/\iota(P)$, and hence all the $T$\--orbits 
$(T \times \{t\})/\iota(P)$ are
$2$\--dimensional Lagrangian submanifolds of $(T \times \mathfrak{t}^*)/\iota(P)$,
diffeomorphic to $T$ itself. 
The fact that the space in case ii) is a symplectic manifold equipped with an effective
action, and that all the elements involved in the definition are well defined was
checked in the proof of Proposition \ref{existence} and the
references therein given. The fact that the $T$\--orbits are $2$\--tori follows by the
same reasoning as in i).

To prove the second part, suppose that $(M, \, \sigma)$ is a compact connected symplectic $4$\--dimensional manifold equipped with an effective symplectic action of a $2$\--torus $T$.
If there exists a $2$\--dimensional symplectic $T$\--orbit then, by Lemma \ref{finitestabilizers00}, every $T$\--orbit
is a $2$\--dimensional symplectic submanifold of $(M, \, \sigma)$.
Assume that none of the $T$\--orbits
is a symplectic $2$\--dimensional submanifold of $(M, \, \sigma)$. Then the antisymmetric bilinear
form $\sigma^{\mathfrak{t}}$ in (\ref{abf}) is degenerate,
and hence it has a one or two--dimensional kernel $\got{l} \subset \mathfrak{t}$. 
If $\op{dim} \got{l}=2$, then $\got{l}=\mathfrak{t}$ and $\sigma^{\mathfrak{t}}=0$, so
every $T$\--orbit is an isotropic submanifold of $(M, \, \sigma)$. Hence
the $2$\--dimensional $T$\--orbits are Lagrangian submanifolds of $(M, \, \sigma)$.
If $\op{dim} \got{l}=1$ there
exists a one--dimensional complement $V$ to $\got{l}$ in $\mathfrak{t}$,
such that the restriction of $\sigma^{\mathfrak{t}}$ to $V$
is a non--degenerate, antisymmetric bilinear form, and hence
identically equal to zero, a contradiction.
Hence either every $T$\--orbit is a $2$\--dimensional symplectic submanifold of $(M, \, \sigma)$, 
or the $2$\--dimensional $T$\--orbits are Lagrangian
submanifolds of $(M, \, \sigma)$. We distinguish three cases according to this.
\\
\\
{\bf Case 3.1.} 
\emph{Suppose that the action of $T$ on $M$ is Hamiltonian}. Because $T$ is $2$\--dimensional,
$(M, \, \sigma)$ is a symplectic toric manifold, and hence by Delzant's theorem \cite{Delzant}, the image
$\mu(M)$ of $M$ under the momentum map $\mu \colon M \to \mathfrak{t}^*$ determines $(M,\, \sigma)$
up to $T$\--equivariant symplectomorphisms (the explicit construction of $M$ from $\mu(M)$ is given
in Delzant's article).
\\
\\
In the next three cases we assume that the $T$\--action on $M$ is not Hamiltonian, so 
that $(M, \, \sigma)$ is not a symplectic toric manifold.
\\
\\
{\bf Case 3.2.}\emph{
Suppose that $(M, \, \sigma)$ has Lagrangian $2$\--dimensional orbits, and that
$T$ does not act freely on $(M, \, \sigma)$ with all $T$\--orbits being Lagrangian $2$\--tori}.
By \cite[Lem. 2.6]{DuPe} the stabilizer subgroups are connected, and hence subtori of $T$.
Since the $T$ action is not free, there exists a $1$\--dimensional stabilizer
subgroup, and hence the Hamiltonian torus $T_{\textup{h}}$ is a $1$\--dimensional subtori of $T$. 
Let $T_{\scriptop{f}}$ be a complementary torus to the 
Hamiltonian torus $T_{\scriptop{h}}$ in $T$,
and let $\mathfrak{t}_{\op{f}}$ be its Lie algebra.
Since $\mathfrak{t}$ 
is $2$\--dimensional and $N=(\mathfrak{t}/\mathfrak{t}_{\textup{h}})^*$, $\op{dim}N=1$. 
Therefore the antisymmetric bilinear form
$c \colon N \times N \to \mathfrak{t}$ in (\ref{thec}) is identically zero.
Then the Lie group 
$G$ in (\ref{G}) is the Cartesian product 
$T_{\scriptop{h}}\times G_{\scriptop{f}}$, 
in which $G_{\scriptop{f}}:=T_{\scriptop{f}}\times N$, where the product in 
$G_{\scriptop{f}}$ is defined by
$
(t_{\scriptop{f}},\, \zeta)(t'_{\scriptop{f}},\, \zeta')=(t_{\scriptop{f}}\, t'_{\scriptop{f}},\, \zeta+\zeta')
$, and by \cite[Lem. 7.1, iii) form. (7.8)]{DuPe} we have that 
$\tau _{\zeta}\in T_{\scriptop{f}}$ for every $\zeta\in P$.

Then the mapping $\iota :\zeta\mapsto ({\tau_{\zeta}}^{-1},\,\zeta )$ 
is a homomorphism from $P$ onto a discrete cocompact subgroup 
of $G_{\scriptop{f}}$.
 With the same proof as in \cite[Prop. 7.2]{DuPe} that (\ref{tsts}) induces a 
 $T$\--equivariant diffeomorphism
$G \times_H M_{\scriptop{h}} \to M$, we obtain that the mapping
$
((t_{\scriptop{f}},\,\zeta ), x)\mapsto 
t_{\scriptop{f}}\cdot\op{e}^{L_{\zeta}}(x)
:G_{\scriptop{f}}\times M_{\scriptop{h}}\to M 
$
induces a $T$\--equivariant diffeomorphism $\alpha _{\scriptop{f}}$ from 
$M_{\scriptop{f}}\times M_{\scriptop{h}}$ onto 
$M$. 
Let $c_{\textup{h}}$
be the $\mathfrak{t}_{\textup{h}}$\--component of $c$ in 
$\mathfrak{t}=\mathfrak{t}_{\textup{h}} \oplus \mathfrak{t}_{\textup{f}}$. 
Let $\pi _{\scriptop{f}}$, 
$\pi _{\scriptop{h}}$ be the projection form 
$M_{\scriptop{f}}\times M_{\scriptop{h}}$ onto the first and the 
second factor, respectively.
The symplectic form 
${\alpha _{\scriptop{f}}}^*\,\sigma$ 
on $M_{\scriptop{f}}\times M_{\scriptop{h}}$ is equal to 
${\pi _{\scriptop{f}}}^*\,\sigma _{\scriptop{f}}+
{\pi  _{\scriptop{h}}}^*\,\sigma _{\scriptop{h}}$,
and the symplectic form 
$\sigma _{\scriptop{f}}$ on $M_{\scriptop{f}}$ is given,
according to (\ref{symplecto}) with $c_{\textup{h}}=0$,
by 
\begin{eqnarray} \label{sstt}
(\sigma _{\scriptop{f}})_{b}(\delta b,\,\delta 'b)=
\delta\zeta (\delta' t)-\delta '\zeta (\delta t).
\end{eqnarray}
Here $b=(t,\,\zeta )\,\iota (P)\in G_{\scriptop{f}}/\iota (P)$, the tangent vectors 
$\delta b=(\delta t,\,\delta\zeta )$ 
and $\delta 'b=(\delta 't,\,\delta '\zeta )$ 
are elements of $\got{t}_{\scriptop{f}}\times N$.
It follows that  $(M,\,\sigma ,\, T)$ is $T$\--equivariantly 
symplectomorphic to 
$(M_{\scriptop{f}},\,\sigma _{\scriptop{f}},\, T_{\scriptop{f}})
\times (M_{\scriptop{h}},\,\sigma _{\scriptop{h}},\, 
T_{\scriptop{h}})$, in which  
$(M_{\scriptop{f}},\,\sigma _{\scriptop{f}},\, T_{\scriptop{f}})$ is a 
compact connected symplectic manifold with a free 
symplectic action $T_{\scriptop{f}}$\--action. 
Here $t\in T$ acts on $M_{\scriptop{f}}\times 
M_{\scriptop{h}}$ by sending $(x_{\scriptop{f}},\, x_{\scriptop{h}})$ 
to $(t_{\scriptop{f}}\cdot x_{\scriptop{f}}
,\, t_{\scriptop{h}}\cdot x_{\scriptop{h}})$, 
if $t=t_{\scriptop{f}}\, t_{\scriptop{h}}$ with 
$t_{\scriptop{f}}\in T_{\scriptop{f}}$  
and $t_{\scriptop{h}}\in T_{\scriptop{h}}$. 

Then the classification of symplectic
toric manifolds \cite{Delzant} of Delzant implies that $M_{\textup{h}}$ is a
$2$\--sphere equipped with an $S^1$\--action by rotations about the vertical axis and endowed with a 
rotationally invariant symplectic form. On the other hand $M_{\scriptop{f}}$
is diffeomorphic to $T \simeq \mathbb{T}^2$, since $S^1$ does not act freely on a surface 
non\--diffeomorphic to $T$,
which in turn implies that $M$ is of the form given in part 2) of the statement.
\\
\\
{\bf Case 3.3.\footnote{
The approach to the proofs of case 3.2 and 3.3 are different. In Case 3.3,
the Delzant manifold $M_{\scriptop{h}}$ is trivial, so the proof is easily obtained as a particular case
of the model $G \times_H M_{\scriptop{h}}$. In Case 3.3 this approach does not lead
to a product $M_{\scriptop{f}} \times M_{\scriptop{h}}$ directly, and hence why we used the same
proof method as in Claim 1.1 but did not specialize.  Although 
if $M_{\scriptop{f}}:=G_{\scriptop{f}}/\iota (P)$ and $H^o:=T_{\textup{h}} \times \{1\}$, the mapping 
$g\mapsto g\, H^o$ defines an isomorphism from 
$G_{\scriptop{f}}$ onto the group $G/H^o$, and an isomorphism from 
$\iota (P)$ onto $H/H^o$, which leads to an identification 
of $M_{\scriptop{f}}$ with $G/H$,
$G \times_H M_{\scriptop{h}}$ is not even $T$\--equivariantly diffeomorphic to 
$G/H \times M_{\scriptop{h}}$ with the induced action.}}
\emph{Suppose that all $T$\--orbits are Lagrangian $2$\--tori and that $T$ acts freely.}
That $(M, \, \sigma)$
is as in part i) of the statement for a unique choice of the ingredients therein,
follows from Theorem \ref{oldclassificationtheorem}
once we show that the model $G \times_H M_{\scriptop{h}}$ of $(M, \, \sigma)$
in Claim 1.1 is of the form $(T \times \got{t}^*)/\iota(P)$
with the $T$\--actions and symplectic form in part i), which we do next. This unique choice of ingredients is given by Definition
\ref{invariantdef}.
Since the $T$\--action is free, $T_{\textup{h}}=\{e\}$ and hence, by (\ref{H}),
$
H=\{(t, \tau_{\zeta}) \,|\, \xi \in P, \, t \tau_{\zeta}=0\}=\iota(P) 
$
and $\mathfrak{t}_{\textup{h}}=\{0\}$, where $\iota$ was given in part i) of the statement.
Also $G=T \times \mathfrak{t}^*$, c.f. (\ref{G}). Because $T_{\textup{h}}=\{1\}$, the
Delzant submanifold $M_{\textup{h}}$ may be chosen to be any point $x \in M$, and there
is a natural $T$\--equivariant symplectomorphism
$$
G \times_H M_{\scriptop{h}} \to G/H \times \{x\} \to G/H = (T \times \got{t}^*)/\iota(P),
$$
where $(T \times \mathfrak{t}^*)/\iota(P)$
is equipped with the symplectic form (\ref{sstt}), but considering
that
$\delta b=(\delta t, \, \delta \zeta), \, \delta' b=
(\delta' t, \, \delta' \zeta) \in \mathfrak{t} \times \mathfrak{t}^*$. 
This expression for the symplectic form can be obtained
by simplifying expression (\ref{symplecto}) according to
$\sigma_{\textup{h}}$, $c_{\textup{h}}$ trivial.\footnote{It follows that
for any $x \in M$, the mapping
$$
(t, \, \zeta) \mapsto t \, \textup{e}^{L_{\zeta}}(x) \colon T \times \mathfrak{t}^* \to M,
$$
induced by the mapping (\ref{tsts}),
gives a $T$\--equivariant symplectomorphism between $(T \times \mathfrak{t}^*)/\iota(P)$
and $(M, \, \sigma)$.}
\\
\\
{\bf Case 3.4.}
\emph{Suppose that all $T$\--orbits are symplectic $2$\--tori.}
Then it follows from Theorem \ref{TT} that $(M, \, \sigma)$ is $T$\--equivariantly
symplectomorphic to the symplectic $T$\--manifold given in part ii) of the statement
for a unique choice of ingredients, and this unique list is given in Definition \ref{def2}
with $f=\mu$ and $\mu_{\scriptop{h}}=f_{\scriptop{h}}$. 
\\
\\
\\
This concludes the proof of Theorem \ref{520}.
$\Box$

\subsection{Corollaries of Theorem \ref{520}}

In the statement of Theorem \ref{520}, 
$(T \times \got{t}^*)/\iota(P)$ is a principal $T$\--bundle 
over the torus $\got{t}^*/P$. Palais and Stewart \cite{ps} 
showed that every principal torus bundle over a torus is diffeomorphic to 
a nilmanifold for a two\--step nilpotent Lie group. We have given an explicit
description of this nilmanifold structure in Theorem \ref{520} i).
Theorem \ref{520} also implies the following results.
\begin{theorem}
The only symplectic $4$\--dimensional manifold equipped with
a non--locally--free and non--Hamiltonian effective symplectic action of a $2$\--torus
is, up to equivariant symplectomorphisms, the product $\T^2 \times S^2$,
where
$\T^2=(\R/\Z)^2$ and
the first factor of $\mathbb{T}^2$ acts on the left factor by translations
on one component, and the second factor acts on $S^2$ by rotations about the
vertical axis of $S^2$. The symplectic form is
a positive linear combination of the
standard translation invariant form on $\mathbb{T}^2$
and the standard rotation invariant form on $S^2$.
\end{theorem}

\begin{proof}
Since the $T$\--action is not Hamiltonian, case 1) in the statement of 
Theorem \ref{520} cannot occur. Since
the action is non--locally--free, there are one--dimensional or two--dimensional stabilizer subgroups, 
$(M, \, \sigma)$ cannot be as in item 3) or 4): in item 3) the stabilizers are all trivial, and in item 4) the stabilizers are finite groups.
\end{proof}

\begin{theorem}
Let $(M, \, \sigma)$ be a non--simply connected 
symplectic $4$\--manifold equipped with a symplectic non--free
action of a $2$\--torus $T$ and such that $M$ is not homeomorphic to 
$T \times S^2$, then $(M, \, \sigma)$ is of the form given
in part ii) of the statement of Theorem \ref{520}, for some $m \ge 1$. 
\end{theorem}

\begin{proof}
This follows by Theorem \ref{520} by the fact that Delzant manifolds
are simply connected. Indeed, every Delzant manifold can be provided with the 
structure of a toric variety defined by a complete fan, 
cf. Delzant \cite{Delzant} and Guillemin \cite[App. 1]{guillemin}, 
and Danilov \cite[Th. 9.1]{danilov} observed that such a 
toric variety is simply connected. The argument is that 
the toric variety has an open cell which is isomorphic to 
$\C ^n$, of which the complement is a complex subvariety 
of complex codimension one. Therefore any loop can be deformed 
into the cell and contracted within the cell to a point. 
\end{proof}

\begin{remark}
The reasons
because of which 
we have imposed that the torus
$T$ is $2$\--dimensional and $(M, \, \sigma)$
is $4$\--dimensional in Theorem \ref{520} are 
\begin{itemize}
\item[i)] There does not exist a classification of $n$\--dimensional smooth orbifolds
if $n>1$.
\item[ii)] In dimensions greater than $2$, the symplectic form is not determined
by a single number (Moser's theorem).
\item[iii)]
If $M$ is not $4$\--dimensional, and $T$ is not $2$\--dimensional, then there are many
cases where not all of the torus orbits are symplectic, and not all of the torus
orbits are isotropic. Other than that Theorem \ref{520} may
be generalized.
Let $(M, \,\sigma)$ be a compact connected symplectic $2n$\--dimensional manifold equipped
with an effective symplectic action of  a torus $T$ and suppose that
one of the following two conditions hold:
\begin{enumerate}
\item
There exists  a $T$\--orbit of dimension $\op{dim}T$ which is a symplectic submanifold
of $(M, \, \sigma)$.
\item
There exists a principal $T$\--orbit which is a coisotropic submanifold of $(M, \, \sigma)$.
\end{enumerate}
Then such symplectic manifolds with $T$\--actions are classified analogously
to Theorem \ref{520}, but in weaker terms (e.g. involving an $n$\--dimensional
orbifold as in item i) above
instead of the first two items of Definition \ref{def1}).
\end{itemize}
\end{remark}

\bigskip\noindent
Alvaro Pelayo\\
Department of Mathematics, University of Michigan\\
2074 East Hall, 530 Church Street, Ann Arbor, MI 48109--1043, USA\\
e\--mail: apelayo@umich.edu

\end{document}